\documentclass[]{article}
\usepackage{amssymb,amsmath,hyperref,verbatim,graphicx,amsthm}
\makeatletter \let\cl@chapter\relax \makeatother
\usepackage{cleveref, xurl, xcolor}
\usepackage{multirow, longtable}
\usepackage{booktabs}
\usepackage{diagbox}
\usepackage{authblk}
\newtheorem{theorem}{Theorem}[section]
\newtheorem{proposition}[theorem]{Proposition}
\newtheorem{corollary}[theorem]{Corollary}

\newtheorem*{claim}{Claim}
\newenvironment{claimproof}[1]{\par\noindent\underline{Proof:}\space#1}{\hfill $\blacksquare$}
\usepackage{pdflscape}
\usepackage{mathtools}

\usepackage{longtable}
\usepackage{enumerate}
\usepackage[ruled,vlined,linesnumbered]{algorithm2e}
\hypersetup{
    colorlinks=true,
    linkcolor=blue,
    filecolor=magenta,      
    urlcolor=cyan,
}
\usepackage{natbib}
 \bibpunct[, ]{(}{)}{,}{a}{}{,}%
\begin{document}

\newcommand{\red}{\color{red}}
\newcommand{\blue}{\color{blue}}
\newcommand{\green}{\color{green}}
\allowdisplaybreaks

\title{Sparsity-Exploiting Distributed Projections onto a Simplex}
\author[]{Yongzheng Dai}
\author[]{Chen Chen}
\affil[]{ISE, The Ohio State University, Columbus, OH, USA}

\maketitle

\begin{abstract}
Projecting a vector onto a simplex is a well-studied problem that arises in a wide range of optimization problems.  Numerous algorithms have been proposed for determining the projection; however, the primary focus of the literature has been on serial algorithms. We present a parallel method that decomposes the input vector and distributes it across multiple processors for local projection. Our method is especially effective when the resulting projection is highly sparse; which is the case, for instance, in large-scale problems with i.i.d. entries. Moreover, the method can be adapted to parallelize a broad range of serial algorithms from the literature. We fill in theoretical gaps in serial algorithm analysis, and develop similar results for our parallel analogues. Numerical experiments conducted on a wide range of large-scale instances, both real-world and simulated, demonstrate the practical effectiveness of the method.
\end{abstract}

\section{Introduction}
\label{sec:intro}
  Given a vector $d\in\mathbb{R}^n$, we consider the following projection of $d$:
\begin{equation}
\label{eq:simplex}
  \mbox{proj}_{\Delta_{b}}(d) := \mbox{argmin}_{v\in \Delta_{b}} \|v-d\|_2,
\end{equation}
where $\Delta_b$ is a scaled standard simplex parameterized by some scaling factor $b>0$,

\[\Delta_{b} := \{v\in\mathbb{R}^n \ |\  \sum_{i=1}^n v_i = b, v\geq 0 \}.\]

\subsection{Applications}
Projection onto a simplex can be leveraged as a subroutine to determine projections onto more complex polyhedra. Such projections arise in numerous settings such as: image processing, e.g. labeling \citep{segmentation2009}, or multispectral unmixing \citep{unmixing2012}; portfolio optimization \citep{portfolio2009}; and machine learning \citep{svm2014}.  As a particular example, projection onto a simplex can be used to project onto the parity polytope \citep{toronto2019}:

\begin{equation} \label{eq:pp}
    \mbox{proj}_{\mathbb{PP}_{n}}(d):=\mbox{argmin}_{v\in \mathbb{PP}_{n}} \|v-d\|_{2},
\end{equation}
where $\mathbb{PP}_{n}$ is a $n$-dimensional parity polytope:
\[\mathbb{PP}_{n}:=\mbox{conv}(\{v\in\{0,1\}^{n}\ |\ \sum_{i=1}^{n}v_{i}=0\ (\mbox{mod} 2)\}).\]

Projection onto the parity polytope arises in linear programming (LP) decoding \citep{l1_LP2016,ADMMLP2011,PP2013,LargeLP2013,reduceLP2015}, which is used for signal processing.

Another example is projection onto a $\ell_{1}$ ball:
\begin{equation}\label{eq: l1 ball}
    \mathcal{B}_{b} := \{v\in\mathbb{R}^n \ |\  \sum_{i=1}^n |v_i| \leq b \}.
\end{equation}
\cite{Duchi2008} demonstrate that the solution to this problem can be easily recovered from projection onto a simplex. Furthermore, projection onto a $\ell_{1}$ ball can, in turn, be used as a subroutine in gradient-projection methods (see e.g. \cite{hybrid2020}) for a variety of machine learning problems that use $\ell_{1}$ penalty, such as: Lasso \citep{lasso1996}; basis-pursuit denoising \citep{bp1998,solveBP2009,hybrid2020}; sparse representation in dictionaries \citep{sparserepresent2003}; variable selection \citep{lassomodel1997}; and classification \citep{classification2017}.

Finally, we note that methods for projection onto the scaled standard simplex and $\ell_{1}$ ball can be extended to projection onto the weighted simplex and weighted $\ell_{1}$ ball \citep{wsimplex2020} (see Section~\ref{sec:weightedprob}). Projection onto the weighted simplex can, in turn, be used to solve the continuous quadratic knapsack problem \citep{qknapscak1992}. Moreover, $\ell_{p}$ regularization can be handled by iteratively solving weighted $\ell_{1}$ projections \citep{wl1lp2008,wl1lp20082,wl1lp2014}. 

\subsection{Contributions}
This paper presents a method to decompose the projection problem and distribute work across (up to $n$) processors. The key insight to our approach is captured by Proposition~\ref{prop: subvector to vector}: the projection of any subvector of $d$ onto a simplex (in the corresponding space with scale factor $b$) will have zero-valued entries only if the full-dimension projection has corresponding zero-valued entries. The method can be interpreted as a sparsity-exploiting distributed preprocessing method, and thus it can be adapted to parallelize a broad range of serial projection algorithms.  We furthermore provide extensive theoretical and empirical analyses of several such adaptations. We also fill in gaps in the literature on serial algorithm complexity. Our computational results demonstrate significant speedups from our distributed method compared to the state-of-the-art over a wide range of large-scale problems involving both real-world and simulated data.  

Our paper contributes to the limited literature on parallel computation for large-scale projection onto a simplex. Most of the algorithms for projection onto a simplex are for serial computing. Indeed, to our knowledge, there is only one published parallel method, and one distributed method for projection problem~(\ref{eq:simplex}). \cite{toronto2019} parallelize a basic sort and scan (specifically prefix sum) approach--a method that we use as a benchmark in our experiments. We also develop a modest but practically significant enhancement to their approach. \cite{iutzeler2018distributed} propose a gossip-based distributed ADMM algorithm for projection onto a Simplex. In this gossip-based setup, one entry of $d$ is given to each agent (e.g. processor), and communication is restricted according to a particular network topology. This differs fundamentally from our approach both in context and intended use, as we aim to solve large-scale problems and moreover our method can accommodate any number of processors up to $n$.
 
The remainder of the paper is organized as follows. Section~\ref{sec:serial} describes serial algorithms from the literature and develops new complexity results to fill in gaps in the literature. Section~\ref{sec:parallel} develops parallel analogues of the aforementioned algorithms using our novel distributed method. Section~\ref{sec:extensions} extends these parallelized algorithms to various applications of projection onto a simplex.  Section~\ref{sec:experiment} describes computational experiments. Section~\ref{sec:conclusion} concludes. Note that all appendices, mathematical proofs, as well as our code and data can be found in the online supplement.

\section{Background and Serial Algorithms}
\label{sec:serial}
 This section begins with a presentation of some fundamental results regarding projection onto a simplex, followed by analysis of serial algorithms for the problem, filling in various gaps in the literature. The final subsection, Section~\ref{sec:summarycomplex}, provides a summary. Note that, for the purposes of average case analysis, we assume uniformly distributed inputs, $d_{1},\dots,d_{n}$ are $\mathrm{i.i.d} \sim U[l,u]$)---a typical choice of distribution in the literature (e.g. \cite{condat2016}).
 
\subsection{Properties of Simplex Projection}
KKT conditions characterize the unique optimal solution $v^{*}$ to problem~(\ref{eq:simplex}):
\begin{proposition}[\cite{basic1974}] \label{prop: define tau}
    For a vector $d \in \mathbb{R}^{n}$ and a scaled standard simplex $\Delta_{b} \in \mathbb{R}^{n}$, there exists a unique $\tau \in \mathbb{R}$ such that
    \begin{equation} \label{eq: def_tau}
        v^{*}_{i}=\max\{d_{i}-\tau,0\}, \ \forall i=1,\cdots,n,
    \end{equation}
    where $v^{*}:=\mathrm{proj}_{\Delta_{b}}(d) \in \mathbb{R}^{n}$.
\end{proposition}
Hence, (\ref{eq:simplex}) can be reduced to a univariate search problem for the optimal \emph{pivot} $\tau$. Note that the nonzero (positive) entries of $v^*$ correspond to entries where $d_i >\tau$.  So for a given value $t\in \mathbb{R}$ and the index set $\mathcal{I} := \{1,...,n\}$, we denote the \emph{active index set} 
\[\mathcal{I}_t := \{i\in\mathcal{I}_t\ |\ d_i>t\},\] 
as the set containing all indices of \emph{active elements} where $d_i>t$. Now consider the following function, which will be used to provide an alternate characterization of $\tau$:
\begin{equation}\label{eq: def f}
    f(t) := \begin{cases}
      \frac{\sum_{i\in \mathcal{I}_t}d_{i}-b}{|\mathcal{I}_t|}-t, & t< \max_i \{d_i\}\\
        -b, & t \geq \max_i \{d_i\} 
    \end{cases}
\end{equation}

\begin{corollary} \label{cor: tau & t}
    For any $t_1,t_2 \in \mathbb{R}$ such that $t_1 < \tau < t_2$, we have
    \[f(t_1) > f(\tau) = 0 > f(t_2). \]
\end{corollary}

The sign of $f$ only changes once, and $\tau$ is its unique root. These results can be leveraged to develop search algorithms for $\tau$, which are presented next. This corollary and the use of $f$ is our own contribution, as we have found it a convenient organizing principle for the sake of exposition  We note, however, that the root finding framework has been in use in the more general constrained quadratic programming literature (see \cite[Equation 5]{reviewer2} and \cite[Section 2, Paragraph 2]{reviewer1}).

\subsection{Sort and Scan}
\label{sec:ss}
Observe that only the greatest $|\mathcal{I}_\tau|$ terms of $d$ are indexed in $\mathcal{I}_\tau$. Now suppose we sort $d$ in non-increasing order:
\[d_{\pi_1}\geq d_{\pi_2}\geq ...\geq d_{\pi_n}. \]
We can sequentially test these values in descending order, $f(d_{\pi_1}), f(d_{\pi_2}), $ etc. to determine $|\mathcal{I}_\tau|$. In particular, from Corollary~\ref{cor: tau & t} we know there exists some $\kappa := |\mathcal{I}_\tau|$ such that $f(d_{\pi_\kappa}) < 0 \leq f(d_{\pi_{\kappa+1}})$. Thus the projection must have $\kappa$ active elements, and since $f(\tau)=0$, we have $\tau = \frac{\sum_{i=1}^\kappa d_{\pi_i}-b }{\kappa}$. We also note that, rather than recalculating $f$ at each iteration, one can keep a running cumulative/prefix sum or \emph{scan} of $\sum_{i=1}^j d_{\pi_i}$ as $j$ increments.

\begin{algorithm}[!ht]
\SetAlgoLined
\LinesNumbered
\SetKwInput{Input}{Input}
\SetKwInput{Output}{Output}
\Input{vector $d=(d_{1},\cdots,d_{n})$, scaling factor $b$.}
\Output{projection $v^{*}$.}
 Sort $d$ as $d_{\pi_1}\geq \cdots \geq d_{\pi_n}$ \;
 Set $\kappa:=\mbox{max}_{1\leq j \leq n}\{j:\frac{\sum_{i=1}^{j}d_{\pi_i}-b}{j}< d_{\pi_j}\}$ (set $\kappa = \pi_n$ if maximum does not exist) \;
 Set $\tau:=\frac{\sum_{i=1}^{\kappa}d_{\pi_i}-b}{\kappa}$\;
 Set $\mathcal{I}_\tau:=\{i\ |\ d_i>\tau\}$, $V_\tau := \{d_i - \tau\ |\ i\in\mathcal{I}_\tau\}$\;
 \textbf{return} $SparseVector(\mathcal{I}_\tau, V_\tau)$.
 \caption{Sort and Scan (\cite{basic1974})}
 \label{Alg:basic}
\end{algorithm}

The bottleneck is sorting, as all other operations are linear time; for instance, \texttt{QuickSort} executes the sort with average complexity $O(n\log n)$ and worst-case $O(n^{2})$, while \texttt{MergeSort} has worst-case $O(n\log n)$ (see, e.g. \citep{sort1993}). Moreover, non-comparison sorting methods can achieve $O(n)$ (see, e.g. \citep{mahmoud2000sorting}), albeit typically with a high constant factor as well as dependence on the bit-size of $d$.

\subsection{Pivot and Partition}
Sort and Scan begins by sorting all elements of $d$, but only the greatest $|\mathcal{I}_\tau|$ terms are actually needed to calculate $\tau$.  Pivot and Partition, proposed by \cite{Duchi2008}, can be interpreted as a hybrid sort-and-scan that attempts to avoid sorting all elements. We present as Algorithm~\ref{Alg:pivot} a variant of this method approach given by \cite{condat2016}.

\begin{algorithm}[!ht]
\renewcommand{\arraystretch}{0.7}
\SetAlgoLined
\LinesNumbered
\SetKwInput{Input}{Input}
\SetKwInput{Output}{Output}
\Input{vector $d=(d_{1},\cdots,d_{n})$, scaling factor $b$.}
\Output{projection $v^{*}$.}
 Set $\mathcal{I}:=\{1,...,n\}$, $\mathcal{I}_\tau:=\emptyset$, $\mathcal{I}_p:=\emptyset$\;
 \While{$\mathcal{I} \neq \emptyset$}{
 	Select a pivot $p\in [\min_{i\in\mathcal{I}}\{d_i\},\max_{i\in\mathcal{I}}\{d_i\}]$\;
 	Set $\mathcal{I}_p := \{i\ |\ d_i>p, i \in \mathcal{I}\}$\;
 	 	\eIf{$(\sum_{i \in {\mathcal{I}_p\cup \mathcal{I}_\tau}}d_i -b) / (|\mathcal{I}_p|+|\mathcal{I}_\tau|)>p$}{
 		Set $\mathcal{I} := \mathcal{I}_p$;
 	}{
 		Set $\mathcal{I}_\tau:=\mathcal{I}_\tau\cup \mathcal{I}_p$, $\mathcal{I} := \mathcal{I}\setminus \mathcal{I}_p$\;
 	}
 }
 Set $\tau:=\frac{\sum_{i\in\mathcal{I}_\tau}d_i-b}{|\mathcal{I}_\tau|}$\;
 Set $V_\tau:=\{d_i - \tau\ |\ i\in \mathcal{I}_\tau\}$\;
 \textbf{return} $SparseVector(\mathcal{I}_\tau, V_\tau)$.
 \caption{Pivot and Partition}
 \label{Alg:pivot}
\end{algorithm}

The algorithm selects a pivot $p\in [\min_i\{d_i\},\max_i\{d_{i}\}]$, which is intended as a candidate value for $\tau$; the corresponding value $f(p)$ is calculated. From Corollary~\ref{cor: tau & t}, if $f(p)>0$, then $p < \tau$ and so $\mathcal{I}_p \supset \mathcal{I}_\tau$; consequently, a new pivot is chosen in the (tighter) interval $ [\min_{i\in\mathcal{I}_p} \{d_i\},\max_{i\in\mathcal{I}_p}\{d_{i}\}]$, which is known to contain $\tau$.  Otherwise, if $f(p) \leq 0$, then $p \geq \tau$, and so we can find a new pivot $p \in [\min_{i\in \bar{\mathcal{I}}_p} \{d_i\},\max_{i\in \bar{ \mathcal{I}}_p}\{d_{i}\}]$, where $\bar{\mathcal{I}}_p:= \{1,...,n\} \setminus \mathcal{I}_p$ is the complement set.  Repeatedly selecting new pivots and creating partitions in this manner results in a binary search to determine the correct active set $\mathcal{I}_\tau$, and consequently $\tau$.

Several strategies have been proposed for selecting a pivot within a given interval. \cite{Duchi2008} choose a random value in the interval, while \cite{pivot2008} uses the median value. The classical approach of \cite{michelot1986} can be interpreted as a special case that sets the initial pivot as $p^{(1)}= (\sum_{i\in \mathcal{I}}d_i-b)/|\mathcal{I}|$, and subsequently $p^{(i+1)} = f(p^{(i)})+p^{(i)}$. This ensures that $p^{(i)}\leq \tau$ which avoids extraneous re-evaluation of sums in the \textbf{if} condition. Note that $p^{(i)}$ generates an increasing sequence converging to $\tau$ \citep[Page 579, Paragraph 2]{condat2016}. Michelot's algorithm is presented separately as Algorithm~\ref{Alg:michelot}.

\begin{algorithm}[ht]
\SetAlgoLined
\LinesNumbered
\SetKwRepeat{Do}{do}{while}
\SetKwInput{Input}{Input}
\SetKwInput{Output}{Output}
\Input{vector $d=(d_{1},\cdots,d_{n})$, scaling factor $b$.}
\Output{projection $v^{*}$.}
 Set $\mathcal{I}_p:=\{1,...,n\}$, $\mathcal{I}:=\emptyset$\;
 \Do{$|\mathcal{I}|>|\mathcal{I}_p|$}{
    Set $\mathcal{I}:=\mathcal{I}_p$\;
    Set $p:=(\sum_{i\in \mathcal{I}}d_i-b)/|\mathcal{I}|$\;
    Set $\mathcal{I}_p := \{i\in\mathcal{I}\ | \ d_{i} > p\}$\;
 }
 Set $V_\tau:=\{d_i - \tau\ |\ i\in \mathcal{I}_\tau\}$\;
 \textbf{return} $SparseVector(\mathcal{I}_\tau, V_\tau)$. 
 \caption{Michelot's method}
 \label{Alg:michelot}
\end{algorithm}

\cite{condat2016} provided worst-case runtimes for each of the aforementioned pivot rules, as well as average case complexity (over the uniform distribution) for the random pivot rule (see Table~\ref{tab: complexity summary}). We fill in the gaps here and establish $O(n)$ runtimes for the median rule as well as Michelot's method. We note that the median pivot method is a linear-time algorithm, but relies on a median-of-medians subroutine \citep{blum1973time}, which has a high constant factor. For Michelot's method, we assume uniformly distributed inputs, $d_{1},\dots,d_{n}$ are $\mathrm{i.i.d} \sim U[l,u]$, and we have
\begin{proposition} \label{prop: Michelot Time Complexity}
	 Michelot's method has an average runtime of $O(n)$. 
\end{proposition}

The same argument holds for the median pivot rule, as half the elements are guaranteed to be removed each iteration, and the operations per iteration are within a constant factor of Michelot's; we omit a formal proof of its $O(n)$ average runtime for brevity.

\subsection{Condat's Method}
\label{sec:condat}
Condat's method \citep{condat2016}, presented as Algorithm~\ref{Alg:condat}, can be seen as a modification of Michelot's method in two ways. First, Condat replaces the initial scan with a Filter to find an initial pivot, presented as Algorithm~\ref{Alg:filter}. Lemma 2 (see Appendix. A) establishes that the Filter provides a greater (or equal) initial starting pivot compared to Michelot's initialization; furthermore, since Michelot approaches $\tau$ from below, this results in fewer iterations (see proof of Proposition~\ref{prop: time complexity of condat}). Second, Condat's method dynamically updates the pivot value whenever an inactive entry is removed from $\mathcal{I}_p$, whereas Michelot's method updates the pivot every iteration by summing over all entries. 

\begin{algorithm}[ht]
\SetAlgoLined
\LinesNumbered
\SetKwInput{Input}{Input}
\SetKwInput{Output}{Output}
\Input{vector $d=(d_{1},\cdots,d_{n})$, scaling factor $b$.}
\Output{$\mathcal{I}_{t}$.}
 Set $\mathcal{I}_{p}:=\{1\}$, $\mathcal{I}_{w}:=\emptyset$, $p:=d_1-b$\;
 \For{$i=2:n$}{
    \If{$d_i>p$}{
        Set $p:=p+\frac{d_i-p}{|\mathcal{I}_{p}|+1}$\;
        \eIf{$p>d_i-b$}{
            Set $\mathcal{I}_p:=\mathcal{I}_p\cup \{i\}$\;
        }{
            Set $\mathcal{I}_w:=\mathcal{I}_w\cup \mathcal{I}_p$, $\mathcal{I}_p:=\{i\}$, $p:=d_i-b$\;
        }
    }
 }
 
 \For{$i \in \mathcal{I}_w$}{
    \If{$d_i>p$}{
        Set $\mathcal{I}_p:=\mathcal{I}_p\cup \{i\}$, $p:=p+\frac{d_i-p}{|\mathcal{I}_p|}$\;
    }
 }
 \textbf{return} $\mathcal{I}_{p}$. 
 \caption{Filter}
 \label{Alg:filter}
\end{algorithm}

\begin{algorithm}[ht]
\SetAlgoLined
\LinesNumbered
\SetKwRepeat{Do}{do}{while}
\SetKwInput{Input}{Input}
\SetKwInput{Output}{Output}
\Input{vector $d=(d_{1},\cdots,d_{n})$, scaling factor $b$.}
\Output{projection $v^{*}$.}
 Set $\mathcal{I}_{p}:=\texttt{Filter}(d,b)$, $p:=\frac{\sum_{i\in\mathcal{I}_p}d_i-b}{|\mathcal{I}_{p}|}$, $\mathcal{I}:=\emptyset$\;
 \Do{$|\mathcal{I}|> |\mathcal{I}_{p}|$}{
    Set $\mathcal{I}:=\mathcal{I}_{p}$\;
 	\For{$i\in\mathcal{I}: d_{i}\leq p$}{
 		Set $\mathcal{I}_p:=\mathcal{I}_p\backslash\{i\}$\;
 		Set $p:=p+\frac{p-d_{i}}{|\mathcal{I}_p|}$\;
 	}
 }
 Set $V_\tau:=\{d_i - \tau\ |\ i\in \mathcal{I}_\tau\}$\;
 \textbf{return} $SparseVector(\mathcal{I}_\tau, V_\tau)$.
 \caption{Condat's method}
 \label{Alg:condat}
\end{algorithm}

\cite{condat2016} supplies a worst-case complexity of $O(n^{2})$. We supplement this with average-case analysis under uniformly distributed inputs, e.g.$d_{1},\dots,d_{n}$ are $\mathrm{i.i.d} \sim U[l,u]$.

\begin{proposition} \label{prop: time complexity of condat}
	Condat's method has an average runtime of $O(n)$.
\end{proposition}

\subsection{Summary of Results}
\label{sec:summarycomplex}

\begin{table}[!ht]
\renewcommand{\arraystretch}{0.7}
\centering
\caption{Time complexity of serial algorithms for projection onto a simplex (new results bolded).}
\setlength{\tabcolsep}{1mm}{
\begin{tabular}{lcc}
\toprule[1pt]
Pivot Rule & Worst Case & Average Case \\
\hline
(Quick)Sort and Scan &  $O(n^2)$ & $O(n)$\\
Michelot's method    & $O(n^2)$   & $\mathbf{O(n)}$\\
Pivot and Partition (Median)     & $O(n)$     & $\mathbf{O(n)}$       \\
Pivot and Partition (Random)     & $O(n^2)$   & $O(n)$      \\
Condat's method & $O(n^2)$ & $\mathbf{O(n)}$\\
Bucket method & $O(cn)$ & $\mathbf{O(cn)}$\\
\bottomrule[1pt]
\end{tabular}
}
\label{tab:pp}
\end{table}

Table~\ref{tab:pp} shows that all presented algorithms attain $O(n)$ performance on average given uniformly i.i.d. inputs. The methods are ordered by publication date, starting from the oldest result. As described in Section~\ref{sec:ss}, Sort and Scan can be implemented with non-comparison sorting to achieve $O(n)$ worst-case performance.  However, as with the linear-time median pivot rule, there are tradeoffs: increased memory, overhead, dependence on factors such as input bit-size, etc.  

Both sorting and scanning are (separately) well-studied in parallel algorithm design, so the Sort and Scan idea lends itself to a natural decomposition for parallelism (discussed in Section~\ref{sec:pss}).  The other methods integrate sorting and scanning in each iterate, and it is no longer clear how best to exploit parallelism directly. We develop in the next section a distributed preprocessing scheme that works around this issue in the case of sparse projections. Note that the table includes the Bucket Method; details on the algorithm are provided in Appendix B.1.

\section{Parallel Algorithms}
\label{sec:parallel}
In Section~\ref{sec:pss} we consider the parallel method proposed by \cite{toronto2019} and propose a modification. In Section~\ref{sec:dist} we develop a novel distributed scheme that can be used to preprocess and reduce the input vector $d$. The remainder of this section analyzes how our method can be used to enhance Pivot and Partition, as well as Condat's method via parallelization of the Filter method. Results are summarized in Section~\ref{sec:parsumm}. We note that the parallel time complexities presented are all unaffected by the underlying PRAM model (e.g. EREW vs CRCW, see  \citep[Chapter 1.4]{xavier1998introduction} for further exposition). This is well-known for parallel mergesort and parallel scan; moreover, our distributed scheme (see Section~\ref{sec:dist}) distributes work for the projection such that memory read/writes of each core are exclusive to that core's partition of $d$.

\subsection{Parallel Sort and Parallel Scan}
\label{sec:pss}
\cite{toronto2019} parallelize Sort and Scan in a natural way: first applying a parallel merge sort (see e.g. \cite[p. 797]{parallelmergesort}) and then a parallel scan \citep{parallelprefxisum} on the input vector. However, their scan calculates $\sum_{i=1}^{j}d_{\pi_i}$ for all $j\in\mathcal{I}$, but only $\sum_{i=1}^{\kappa}d_{\pi_\kappa}$ is needed to calculate $\tau$. We modify the algorithm accordingly, presented as Algorithm~\ref{Alg:IPscan}: checks are added (lines 7 and 14) in the for-loops to allow for possible early termination of scans.  As we are adding constant operations per loop, Algorithm~\ref{Alg:IPscan} has the same complexity as the original Parallel Sort and Scan. We combine this with parallel mergesort in Algorithm~\ref{Alg:PsortIPscan} and empirically benchmark this method with the original (parallel) version in Section~\ref{sec:experiment}.

\begin{algorithm}[ht]
\SetAlgoLined
\LinesNumbered
\SetKwInput{Input}{Input}
\SetKwInput{Output}{Output}
\Input{sorted vector $d_{\pi_1},\cdots, d_{\pi_n}$, scaling factor $b$}
\Output{$\tau$}
 Set $T:=\lceil\log_{2}n\rceil$, $s[1],...,s[n]=d_{\pi_1},...,d_{\pi_n}$\;
 \For{$j=1:T$}{
    \For(Parallel){$i=2^{j}:2^{j}:\min(n,2^{T})$}{
        Set $s[i]:=s[i]+s[i-2^{j-1}]$\;
    }
    Set $\kappa:=\min(n,2^{j})$\;
    \If{$\frac{s[\kappa]-a}{\kappa} \geq d_{\pi_\kappa}$}{
        break loop\;
    }
 }
 Set $p:=2^{j-1}$\;
 \For{$i=j-1:-1:1$}{
    Set $\kappa:=\min(p+2^{i-1},n)$, $s[\kappa]:=s[\kappa]+s[p]$\;
    \If{$\frac{s[\kappa]-a}{\kappa}<d_{\pi_\kappa}$}{
        break loop\;
    }
 }
 Set $\tau:=\frac{s[\kappa]-b}{\kappa}$\;
 \textbf{return} $\tau$.
 \caption{Parallel Partial Scan}
 \label{Alg:IPscan}
\end{algorithm}

\begin{algorithm}[ht]
\SetAlgoLined
\LinesNumbered
\SetKwInput{Input}{Input}
\SetKwInput{Output}{Output}
\Input{vector $d=(d_{1},\cdots,d_{n})$, scaling factor $b$.}
\Output{projection $v^{*}$.}
 Parallel mergesort $d$ so that $d_{\pi_1}\geq \cdots \geq d_{\pi_n}$\;
 Set $\tau=\texttt{PPScan}(\{d_{\pi_i}\}_{1\leq i \leq n},b)$\;
 Set $\mathcal{I}_\tau := \{i\ |\ d_i > \tau\}$, $V_\tau:=\{d_i - \tau\ |\ i\in \mathcal{I}_\tau\}$\;
 \textbf{return} $SparseVector(\mathcal{I}_\tau, V_\tau)$.
 \caption{Parallel Mergesort and Partial Scan}
 \label{Alg:PsortIPscan}
\end{algorithm}

\subsection{Sparsity-Exploiting Distributed Projections}
\label{sec:dist}
Our main idea is motivated by the following two theorems that establish that projections with i.i.d. inputs and fixed right-hand-side $b$ become increasingly sparse as the problem size $n$ increases. 

\begin{theorem} \label{thm: average active terms}
    $E[|\mathcal{I}_\tau|]  < \sqrt{\frac{2b(n+1)}{u-l}+\frac{1}{4}}+\frac{1}{2}$.
\end{theorem}
Theorem~\ref{thm: average active terms} establishes that, for i.i.d. uniformly distributed inputs, the projection has $O(\sqrt{n})$ active entries in expectation and thus has considerable sparsity as $n$ grows; we also show this in the computational experiments of Appendix E.1.

\begin{theorem} \label{thm: sparsity}
	Suppose $d_1,...,d_n$ are $\mathrm{i.i.d.}$ from an arbitrary distribution $X$, with PDF $f_X$ and CDF $F_X$. Then, for any $\epsilon>0$, $P(\frac{|\mathcal{I}_\tau|}{n}\leq \epsilon)=1$ as $n\rightarrow \infty$.
\end{theorem}

 Theorem~\ref{thm: sparsity} establishes arbitrarily sparse projections over arbitrary i.i.d. distributions given fixed $b$. Note that if $b$ is sufficiently large with respect to $n$ (rather than fixed), then the resulting projection could be too dense to attain the theorem result. However, sparsity can be assured provided $b$ does not grow too quickly with respect to $n$, namely:
\begin{corollary}
     Theorem~\ref{thm: sparsity} holds true if $b\in o(n)$.
\end{corollary}

We apply Theorem~\ref{thm: sparsity} to example distributions in Appendix C, and test the bounds empirically in Appendix E.1.

\begin{proposition} \label{prop: subvector to vector}
	Let $\hat{d}$ be a subvector of $d$ with $m\leq n$ entries; moreover, without loss of generality suppose the subvector contains the first $m$ entries. Let $\hat v^*$ be the projection of $\hat{d}$ onto the simplex $\hat \Delta := \{v\in\mathbb{R}^m \ |\  \sum_{i=1}^m v_i = b, v\geq 0 \}$, and $\hat{\tau}$ be the corresponding pivot value. Then, $\tau\geq \hat{\tau}$. Consequently, for $1\leq i\leq m$ we have that $\hat v_i^*=0 \implies v_i^* =0$. 
\end{proposition}
 Proposition~\ref{prop: subvector to vector} tells us that if we project a subvector of some length $m \leq n$ onto the same $b$-scaled simplex in the corresponding $\mathbb{R}^m$ space, the zero entries in the projected subvector must also be zero entries in the projected full vector.

\begin{figure}[!htbp]
\centering
\includegraphics[width=0.8\linewidth]{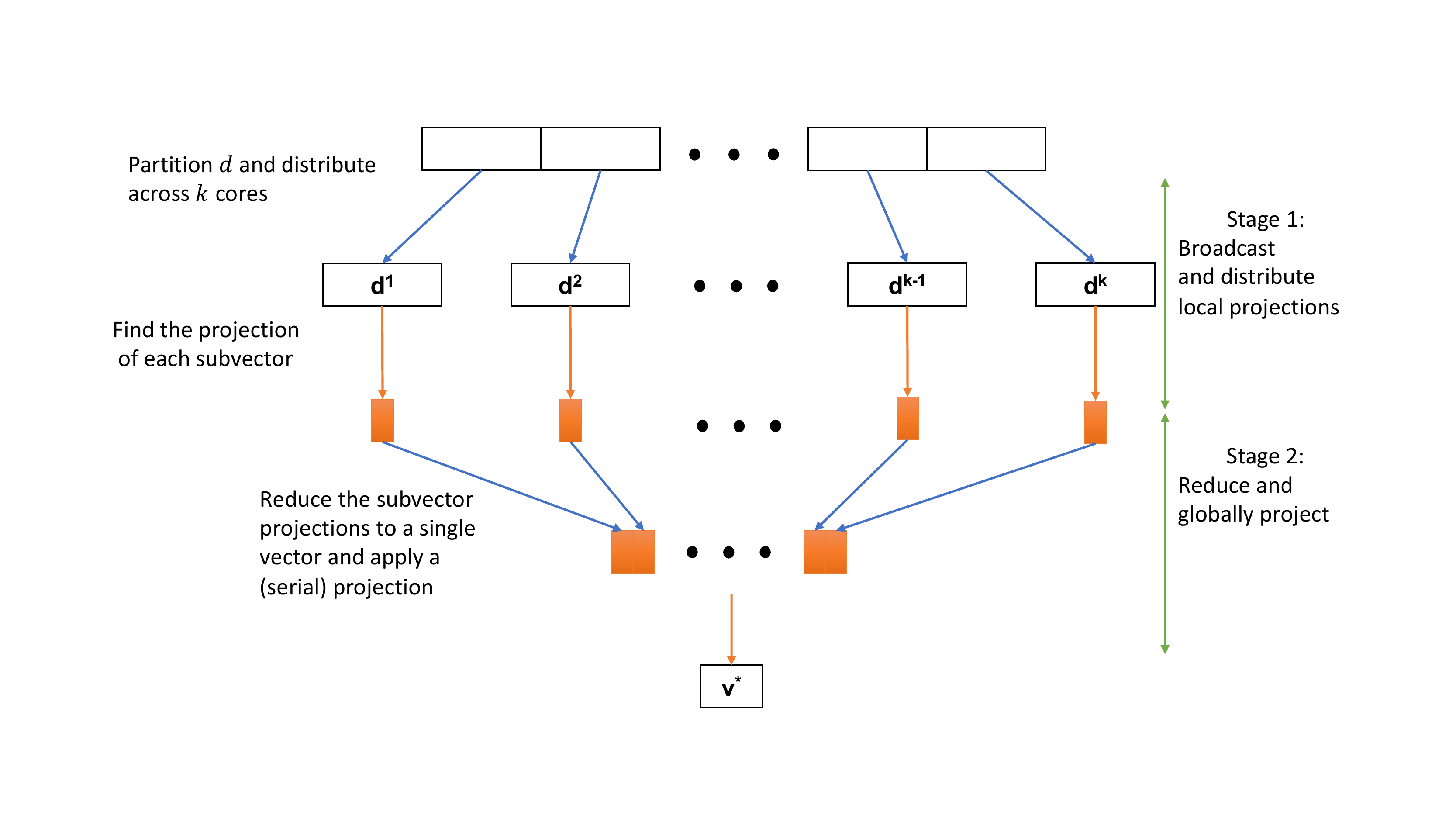}
\caption{Distributed Projection Algorithm}
\label{fig: distributed alg}
\end{figure}

Our idea is to partition and distribute the vector $d$ across cores (broadcast); have each core find the projection of its subvector (local projection); and combine the nonzero entries from all local projections to form a vector $\hat v$ (reduce), and apply a final (global) projection to $\hat v$. The method is outlined in Figure~\ref{fig: distributed alg}. Provided the projection $v^*$ is sufficiently sparse, which (for instance) we have established is the case for i.i.d. distributed large-scale problems, we can expect $\hat v$ to have been far less than $n$ entries.  We demonstrate the practical advantages of this procedure with various computational experiments in Section~\ref{sec:experiment}.

\subsection{Parallel Pivot and Partition}
The distributed method outlined in Figure~\ref{fig: distributed alg} can be applied directly to Pivot and Partition, as described in Algorithm~\ref{Alg: Dpivot}.  Note that, as presented, $v^{*}$ is a sparse vector: entries not processed in the final Pivot and Project iteration are set to zero (recall Proposition~\ref{prop: subvector to vector}).  


\begin{algorithm}[ht]
\SetAlgoLined
\LinesNumbered
\SetKwInput{Input}{Input}
\SetKwInput{Output}{Output}
\Input{vector $d=(d_{1},\cdots,d_{n})$, scaling factor $b$, cores number $k$.}
\Output{projection $v^{*}$.}
 Partition $d$ into subvectors $d^1,...,d^k$ of dimension $\leq \frac{n}{k}$ \;
 Set $\mathcal{I}_i$ be the active set from $\texttt{Pivot\_Project}(d^i,b)$ \ \ (distributed across cores $i=1,...,k$)\;
 Set $\hat{\mathcal{I}}:=\cup_{i=1}^{k}\mathcal{I}_i$\;
 \textbf{return} Sparse vector from the projection of $\{v_i\}_{i\in\hat{\mathcal{I}}}$. 
 \caption{Parallel Pivot and Partition}
 \label{Alg: Dpivot}
\end{algorithm}
 
We assume uniformly distributed inputs, $d_{1},\dots,d_{n}$ are $\mathrm{i.i.d} \sim U[l,u]$, and we have
\begin{proposition} \label{prop: Time Complexity of DMichelot}
	Parallel Pivot and Partition with either the median, random, or Michelot's pivot rule, has an average runtime of $O(\frac{n}{k}+\sqrt{kn})$. 
\end{proposition}

In the worst case we may assume the distributed projections are ineffective, $\mbox{dim}(\hat v) \in O(n)$, and so the final projection is bounded above by $O(n^2)$ with random pivots and Michelot's method, and $O(n)$ with the median pivot rule.

\subsection{Parallel Condat's Method}
We could apply the distributed sparsity idea as a preprocessing step for Condat's method. However, due to Proposition~\ref{prop: filter result} (and confirmed via computational experiments) we have found that Filter tends to discard many non-active elements. Therefore, we propose instead to apply our distributed method to parallelize the Filter itself. Our Distributed Filter is presented in Algorithm~\ref{Alg: Dfilter}: we partition $d$ and broadcast it to the cores, and in each core we apply (serial) Filter on its subvector. Condat's method with the distributed Filter is presented as Algorithm~\ref{Alg: Dcondat}.

\begin{algorithm}[ht]
\SetAlgoLined
\LinesNumbered
\SetKwInput{Input}{Input}
\SetKwInput{Output}{Output}
\Input{vector $d=(d_{1},\cdots,d_{n})$, scaling factor $b$, $k$ cores.}
\Output{Index set $\mathcal{I}$ of Stage $1$.}
 Partition $\mathcal{I}$ into index sets $\{\mathcal{I}_{1},\cdots,\mathcal{I}_{k}\}$ such that $\mathcal{I}_i \leq \frac{n}{k}, i=1,,,.k$\;
  \For(parallel){$i=1:k$}{
   Update $\mathcal{I}_{i}$ with (serial) \texttt{Filter}($d_{\mathcal{I}_i},b$)\;
    Set $p^i:=\frac{\sum_{j\in \mathcal{I}_{i}}d_j-b}{|\mathcal{I}_{i}|}$\;
    \For{$j\in \mathcal{I}_i$}{
        \If{$d_j\leq p^i$}{
            Set $p^i:=p^i+\frac{p^i-d_j}{|\mathcal{I}_i|-1}$, $\mathcal{I}_i:=\mathcal{I}_i\backslash \{j\}$
        }
    }
 }
 \textbf{return} $\mathcal{I}:=\cup_{i=1}^{k}\mathcal{I}_i$. 
 \caption{Distributed Filter (Dfilter)}
 \label{Alg: Dfilter}
\end{algorithm}

\begin{algorithm}[ht]
\SetAlgoLined
\LinesNumbered
\SetKwInput{Input}{Input}
\SetKwInput{Output}{Output}
\SetKwRepeat{Do}{do}{while}
\Input{vector $d=(d_{1},\cdots,d_{n})$, scaling factor $b$, $k$ cores.}
\Output{projection $v^{*}$.}
 Set $\mathcal{I}_p:=\texttt{Dfilter}(d,b,k)$, $\mathcal{I}:=\emptyset$, $p:=\frac{\sum_{i\in\mathcal{I}_p}d_i-b}{|\mathcal{I}_p|}$\;
 \Do{$|\mathcal{I}|> |\mathcal{I}_p|$}{
 	Set $\mathcal{I}:=\mathcal{I}_p$\;
 	\For{$i\in\mathcal{I}: d_{i}\leq p$}{
 		Set $\mathcal{I}_p:=\mathcal{I}_p\backslash\{i\}$, $p:=p+\frac{p-d_{i}}{|\mathcal{I}_p|}$\;
 	}
 }
 Set $V_\tau:=\{d_i - \tau\ |\ i\in \mathcal{I}_\tau\}$\;
 \textbf{return} $SparseVector(\mathcal{I}_\tau, V_\tau)$.
 \caption{Parallel Condat's method}
 \label{Alg: Dcondat}
\end{algorithm}

We assume uniformly distributed inputs, $d_{1},\dots,d_{n}$ are $\mathrm{i.i.d} \sim U[l,u]$, and we have
\begin{proposition} \label{prop: filter result}
 Let $\mathcal{I}_p$ be the output of $\mathrm{Filter}(d,b)$. Then $E[|\mathcal{I}_p|] \in O(n^{\frac{2}{3}})$.
\end{proposition}

Under the same assumption (uniformly distributed inputs) to Proposition~\ref{prop: filter result}, we have
\begin{proposition} \label{pro: time complexity of Dcondat}
	Parallel Condat's method has an average complexity $O(\frac{n}{k}+\sqrt[3]{kn^{2}})$.
\end{proposition}

In the worst-case we can assume Distributed Filter is ineffective ($|\mathcal{I}_p|\in O(n)$), and so the complexity of Parallel Condat's method is $O(n^2)$, same as the serial method.

\subsection{Summary of Results}
\label{sec:parsumm}
\begin{table}[!htbp]
\renewcommand{\arraystretch}{0.7}
\centering
\caption{\label{tab: complexity summary} Time complexity of serial vs parallel algorithms with problem dimension $n$ and $k$ cores}
\begin{tabular}{lcc}
\toprule[1pt]
Method & Worst case complexity  & Average complexity\\
\hline
Quicksort + Scan &$O(n^2)$ &$O(n\log n)$\\
(P)Mergesort + Scan &$O(\frac{n}{k}\log n)$ &$O(\frac{n}{k}\log n)$\\
(P)Mergesort + Partial Scan &$O(\frac{n}{k}\log n)$ &$O(\frac{n}{k}\log n)$\\
Michelot &$O(n^{2})$ &$O(n)$\\
(P)Michelot &$O(n^2)$ &$O(\frac{n}{k}+\sqrt{kn})$\\
Condat &$O(n^{2})$ &$O(n)$\\
(P)Condat &$O(n^{2})$ &$O(\frac{n}{k}+\sqrt[3]{kn^{2}})$\\
\bottomrule[1pt]
\end{tabular}
\end{table}

Complexity results for parallel algorithms developed throughout this section, as well as their serial counterparts, are presented in Table~\ref{tab: complexity summary}. Parallelized Sort and Scan has a dependence on $\frac{1}{k}$; indeed, sorting (and scanning) are very well-studied problems from the perspective of parallel computing. Parallel (Bitonic) Merge Sort has an average-case (and worst case) complexity in $O((n\log n)/k)$ \citep{sortcomplex}, and Parallel Scan has an average-case (and worst case) complexity in $\Theta(n/k+\log k)$ \citep{parallelscan}; thus, running parallel sort followed by parallel scan is $O(\frac{n}{k}\log n)$. Now, Michelot's and Condat's serial methods are observed to have favorable practical performance in our computational experiments; this is expected as these more modern approaches were explicitly developed to gain practical advantages in e.g. the constant runtime factor. Conversely, our distributed method does not improve upon the \emph{worst-case} complexity of Michelot's and Condat's methods, but is able to attain a $\frac{1}{k}$ factor for \emph{average} complexity for $n\gg k$, which is the case on large-scale instances, i.e. for all practical purposes.  The average case analyses were conducted under the admittedly limited (typical) assumption of uniform i.i.d. entries, but our computational experiments over other distributions and real-world data confirm favorable practical speedups from our parallel algorithms.

\section{Parallelization for Extensions of Projection onto a Simplex}
\label{sec:extensions}
This section develops extensions involving projection onto a simplex, to be used for experiments in Section~\ref{sec:experiment}. 

\subsection{Projection onto the $\ell_{1}$ Ball}
\label{sec:ell1desc}
Consider projection onto an $\ell_1$ ball:
\begin{equation}
    \label{eq: project l1 ball}
    \mbox{Proj}_{\mathcal{B}_b}(d):=\arg\min_{v\in\mathcal{B}_b}\|v-d\|_2,
\end{equation}
where $\mathcal{B}_b$ is given by Equation~(\ref{eq: l1 ball}). \cite{Duchi2008} show Problem~(\ref{eq: project l1 ball}) is linear-time reducible to projection onto simplex (see \cite[Section 4]{Duchi2008}). Hence, any parallel method for projection onto a simplex can be applied to Problem~(\ref{eq: project l1 ball}). 


As mentioned in Section~\ref{sec:intro}, Problem~(\ref{eq: project l1 ball}) can itself be used as a subroutine in solving the Lasso problem, via (e.g.) Projected Gradient Descent (PGD) (see e.g. \cite[Exercise 10.2]{boyd2004convex}). To handle large-scale datasets, we instead use the mini-batch gradient descent method \citep[Sec. 12.5]{minibatch} in PGD in Sec~\ref{sec:experiment}. 


\subsection{Centered Parity Polytope Projection}
\label{sec:cpppdesc}
Leveraging the solution to Problem~(\ref{eq:simplex}), Wasson et al. \cite[Algorithm 2]{toronto2019} develop a method to project a vector onto the centered parity polytope, $\mathbb{PP}_{n}-\frac{1}{2}$ (recall Problem~\ref{eq:pp}); we present a slightly modified version as Algorithm 2 in Appendix D. The modification is on line $11$, where we determine whether a simplex projection is required to avoid unnecessary operations; the original method executes line $14$ before line $11$.

\section{Numerical Experiments} \label{sec:experiment}
All algorithms were implemented in Julia 1.5.3 and run on a single node from the Ohio Supercomputer Center \citep{osc}. This node includes 4 sockets, with each socket containing a 20 Intel Xeon Gold 6148 CPUs; thus there are 80 cores at this node. The node has 3 TB of memory and runs 64-bit Red Hat Enterprise Linux 3.10.0. The code and data are available at All code and data can be found at: \href{https://github.com/foreverdyz/Parallel_Projection}{Github}\footnote{https://github.com/foreverdyz/Parallel\_Projection} or the IJOC repository \citep{code}

\subsection{Testing Algorithms}
In this subsection, we compare runtime results for serial methods and their parallel versions with two measures: absolute speedup and relative speedup. Absolute speedup is the parallel method's runtime vs the fastest serial method's runtime, e.g. serial Condat time divided by parallel Sort and Scan time; relative speedup is the parallel method's runtime vs its serial equivalent's runtime, e.g. serial Sort and Scan time divided by parallel Sort and Scan time. We test parallel implementations using $8, 16, 24, ..., 80$ cores. Note that we have verified that each parallel algorithm run on a single core is slower than its serial equivalent (see Appendix E.4).

\subsubsection{Projection onto Simplex}
\begin{figure}[!htbp]
    \centering
    \begin{minipage}[b]{0.45\textwidth}
         \centering
         \includegraphics[width=\textwidth]{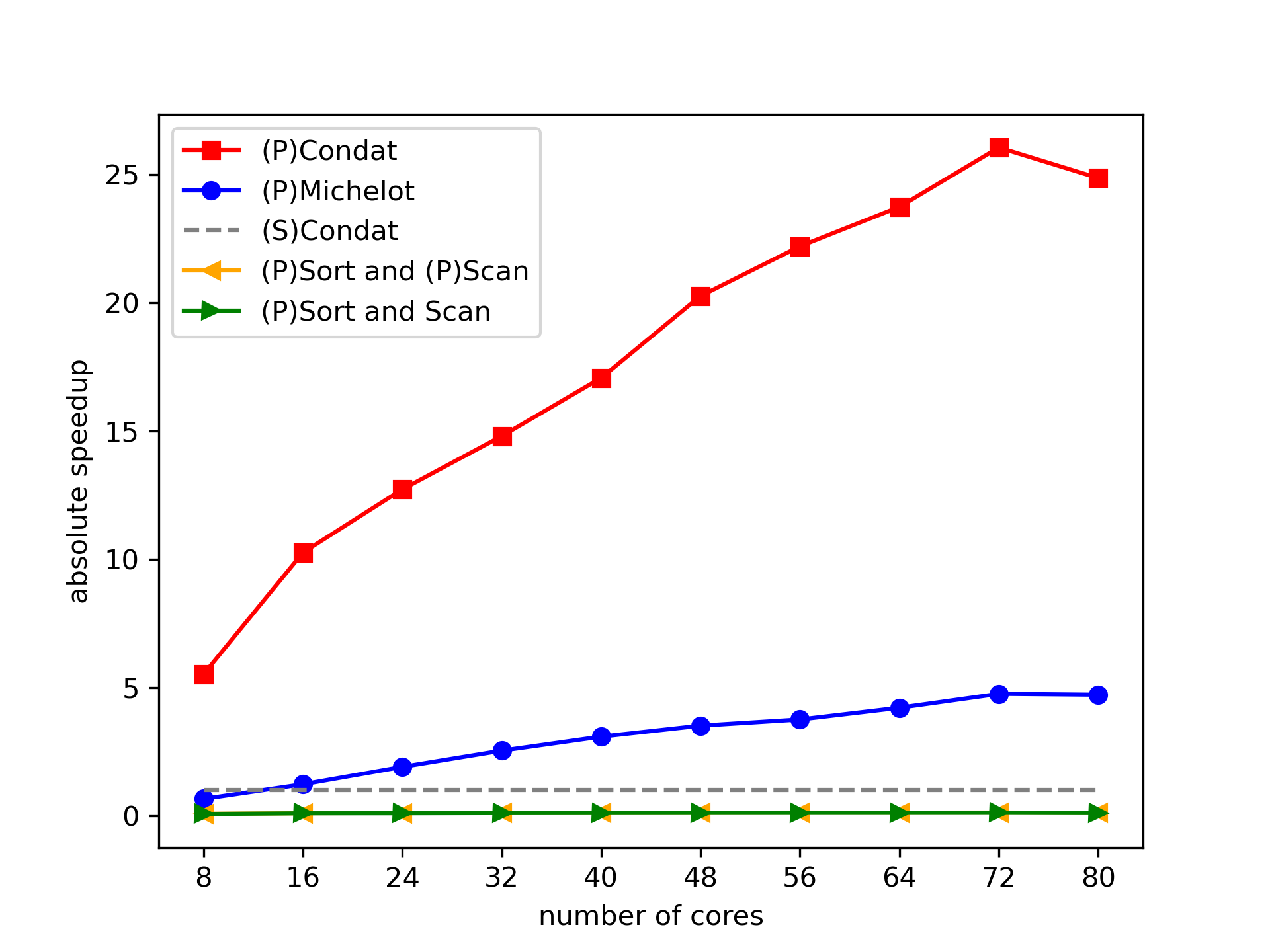}
    \end{minipage}
    \hfill
    \begin{minipage}[b]{0.45\textwidth}
         \centering
         \includegraphics[width=\textwidth]{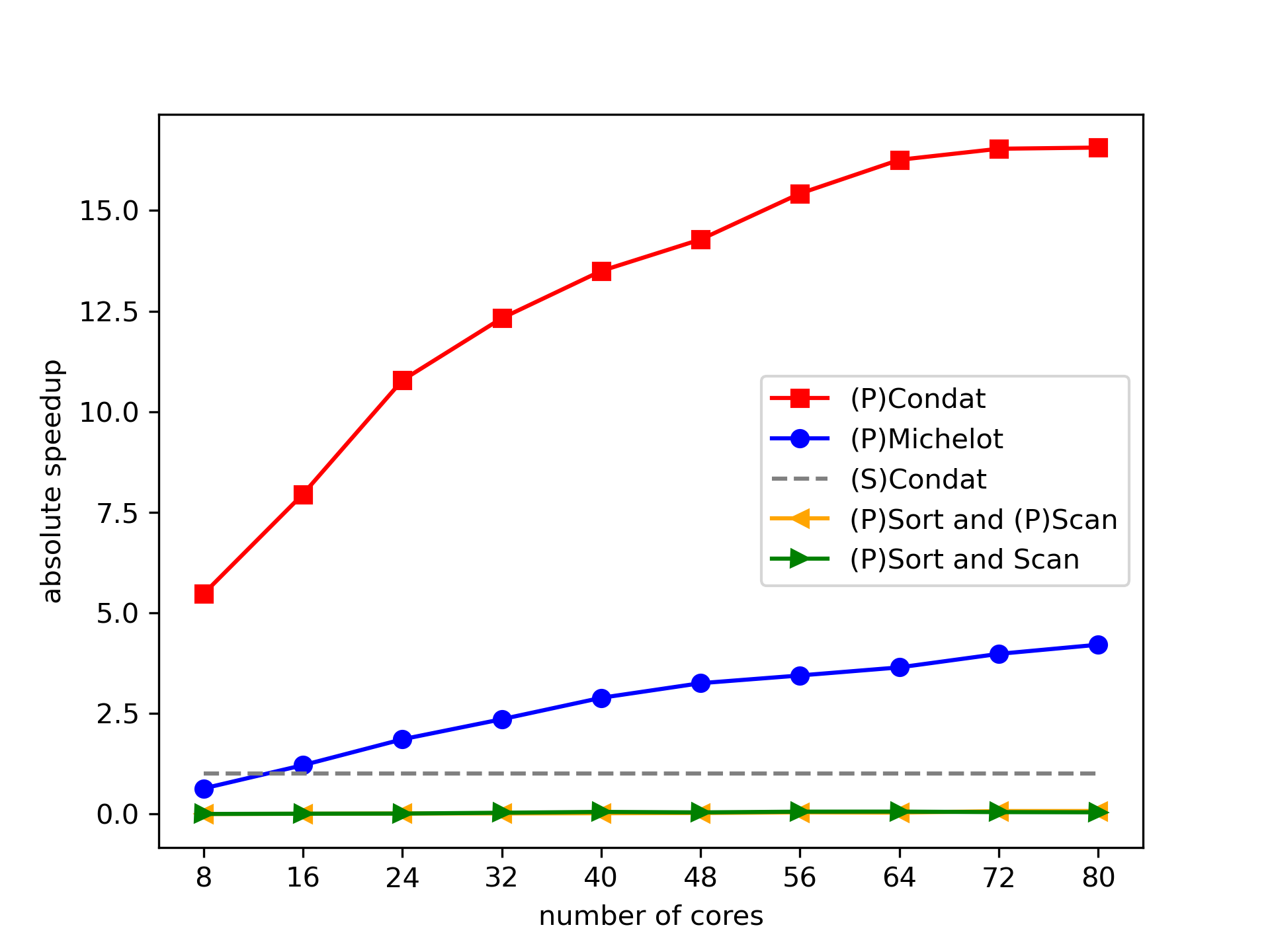}
    \end{minipage}
    \hfill
    \begin{minipage}[b]{0.45\textwidth}
         \centering
         \includegraphics[width=\textwidth]{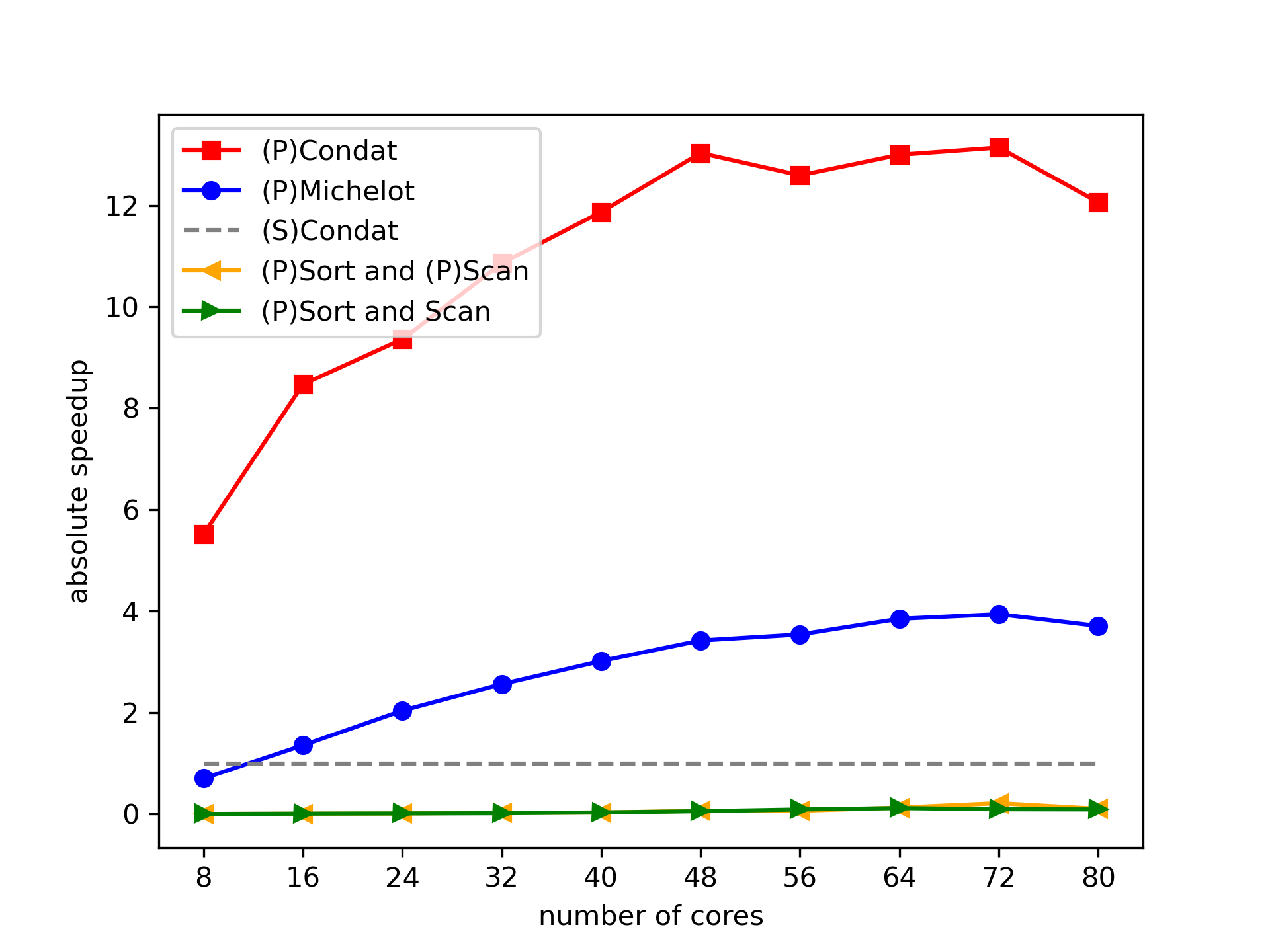}
    \end{minipage}
    \caption{Absolute speedup vs cores in simplex projection. Each line represents a different algorithm, and each graph represents different input distributions: $N(0,1)$ (left top), $U[0,1]$ (right top), and $N(0, 10^{-3})$ (below).}
    \label{fig:simplex_comp}
\end{figure}

\begin{figure}[!htbp]
    \centering
    \begin{minipage}[b]{0.45\textwidth}
         \centering
         \includegraphics[width=\textwidth]{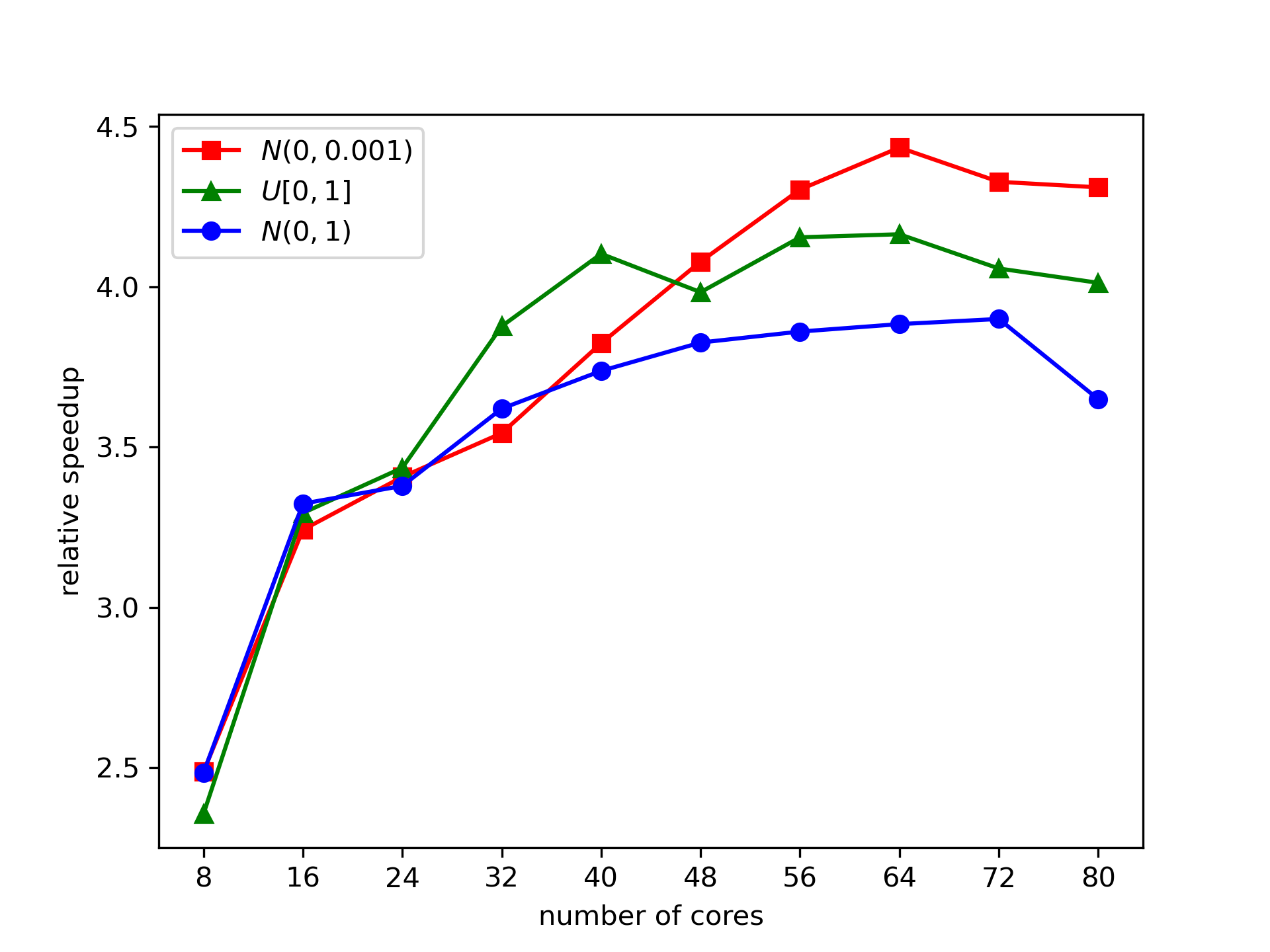}
    \end{minipage}
    \hfill
    \begin{minipage}[b]{0.45\textwidth}
         \centering
         \includegraphics[width=\textwidth]{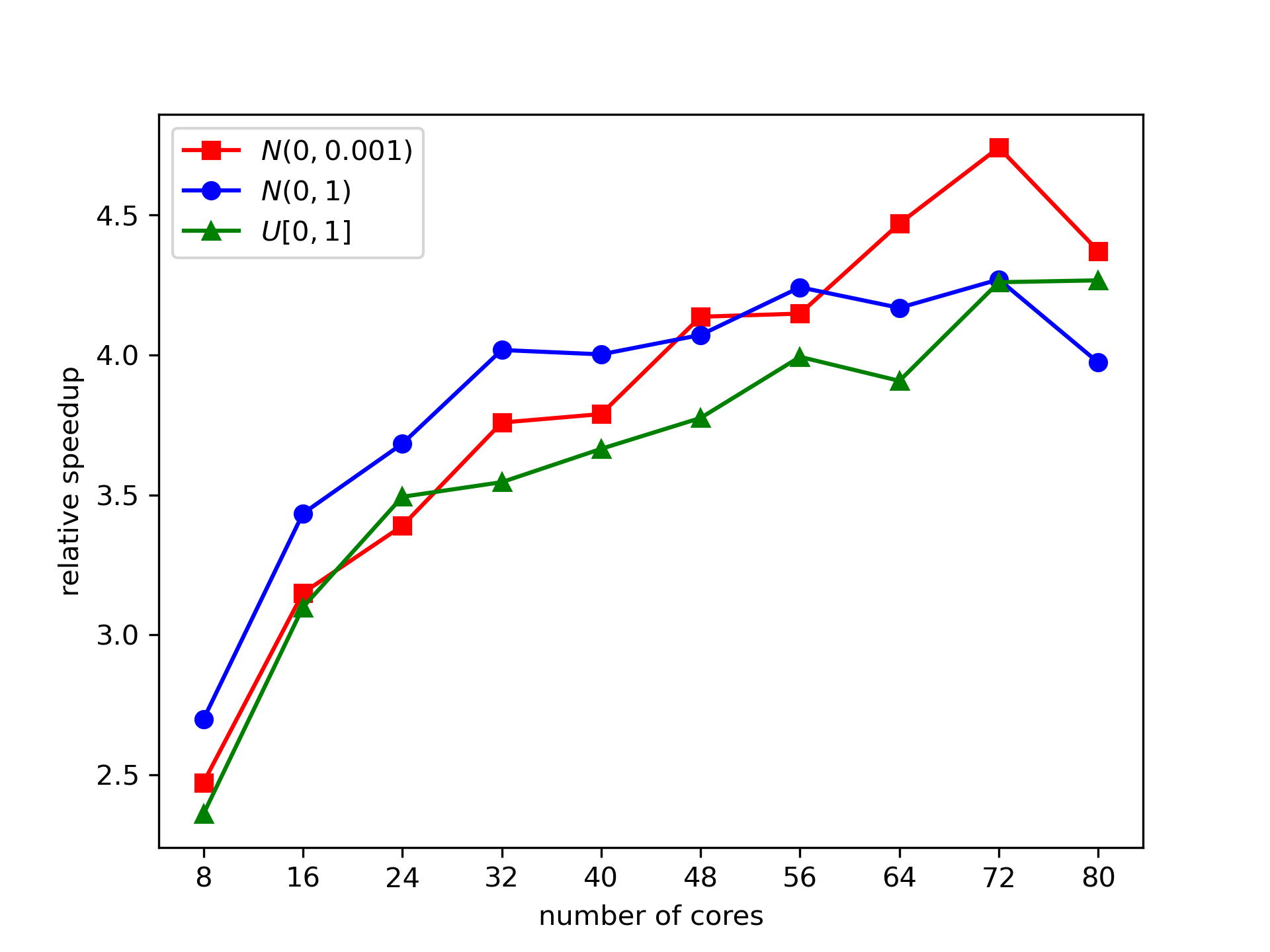}
    \end{minipage}
    \hfill
    \begin{minipage}[b]{0.45\textwidth}
         \centering
         \includegraphics[width=\textwidth]{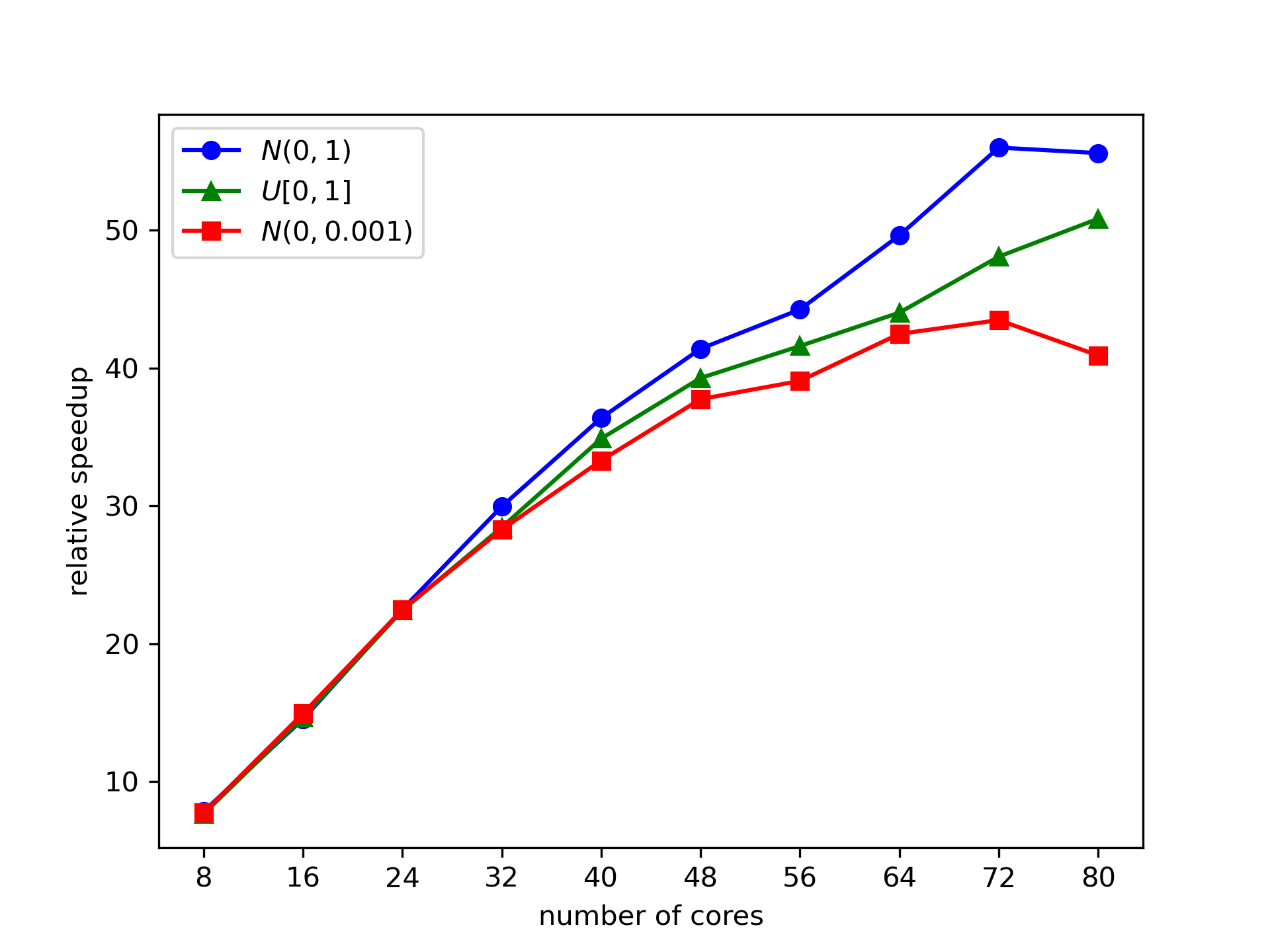}
    \end{minipage}
    \hfill
    \begin{minipage}[b]{0.45\textwidth}
         \centering
         \includegraphics[width=\textwidth]{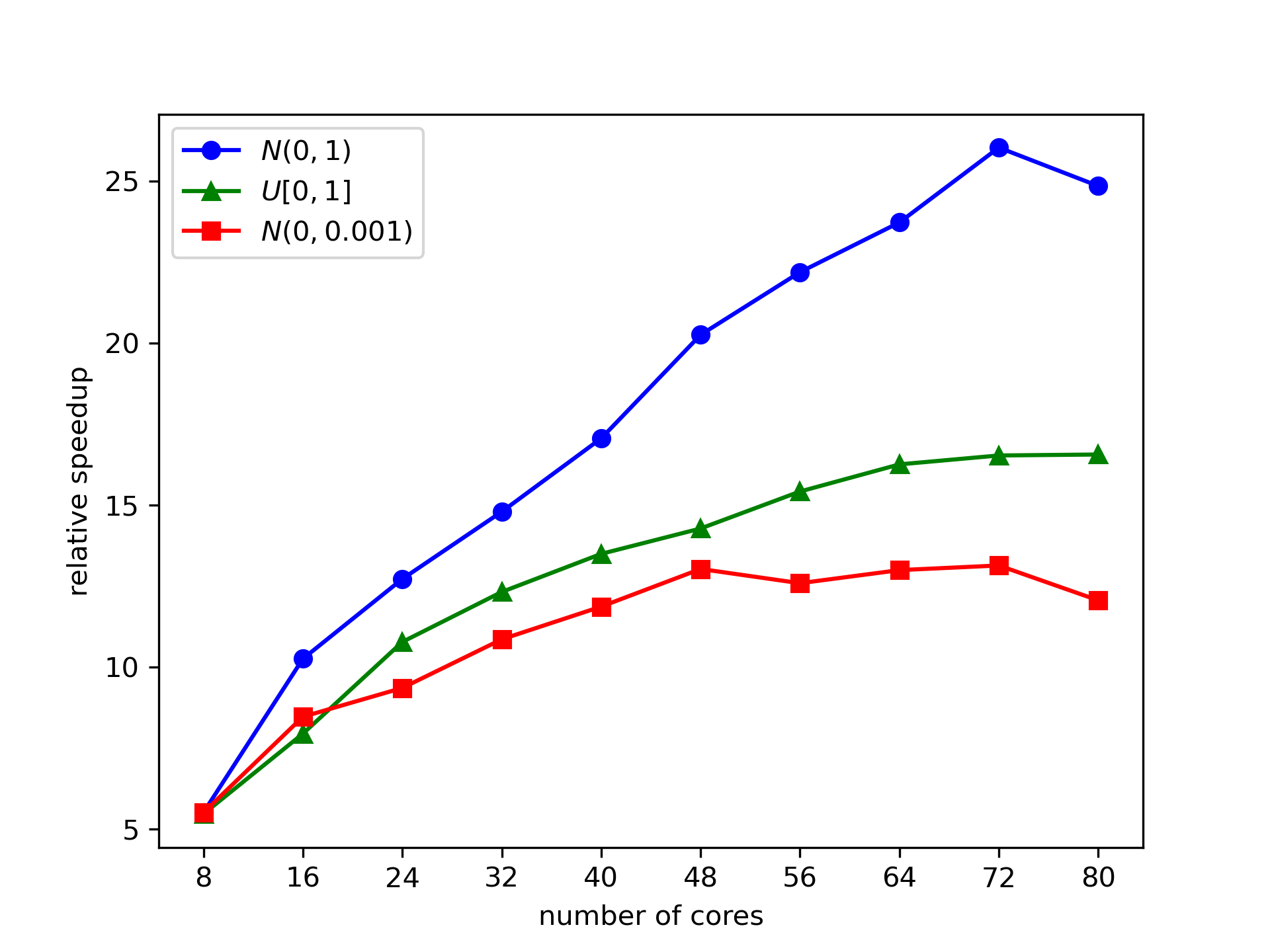}
    \end{minipage}
    \caption{Relative speedup vs cores in simplex projection. Each line represents a different input distribution, and 4 graphs represent 4 different  projection methods, Parallel Sort and Scan (left top), Parallel Sort and Partial Scan (right top), Parallel Michelot (left below), and Parallel Condat (right below).}
    \label{fig:simplex_projection}
\end{figure}

Instances in Figures~\ref{fig:simplex_comp} and~\ref{fig:simplex_projection} are generated with a size of $n=10^{8}$ and scaling factor $b = 1$. Inputs $d_i$ are drawn i.i.d. from three benchmark distributions: $U[0,1]$, $N(0,1)$, and $N(0,10^{-3})$.  This is a common benchmark used in previous works e.g. \citep{Duchi2008, condat2016}. Serial Condat is the benchmark serial algorithm (i.e. with the fastest performance), with the dotted line representing a speedup of 1.  Our parallel Condat algorithm achieves up to 25x absolute speedup over this state-of-the-art serial method.  Parallel Sort and Scan ran slower than the serial Condat's method, due to the dramatic relative slowdown in the serial version. In terms of relative speedup, our method offers superior performance compared to the Sort and Scan approach.  Although not visible on the absolute speedup graph, we note that in relative speedup it can be seen that our partial Scan technique offers some modest improvements over the standard parallel Sort and Scan.

\begin{figure}[!htbp]
    \centering
    \begin{minipage}[b]{0.45\textwidth}
        \centering
        \includegraphics[width=\textwidth]{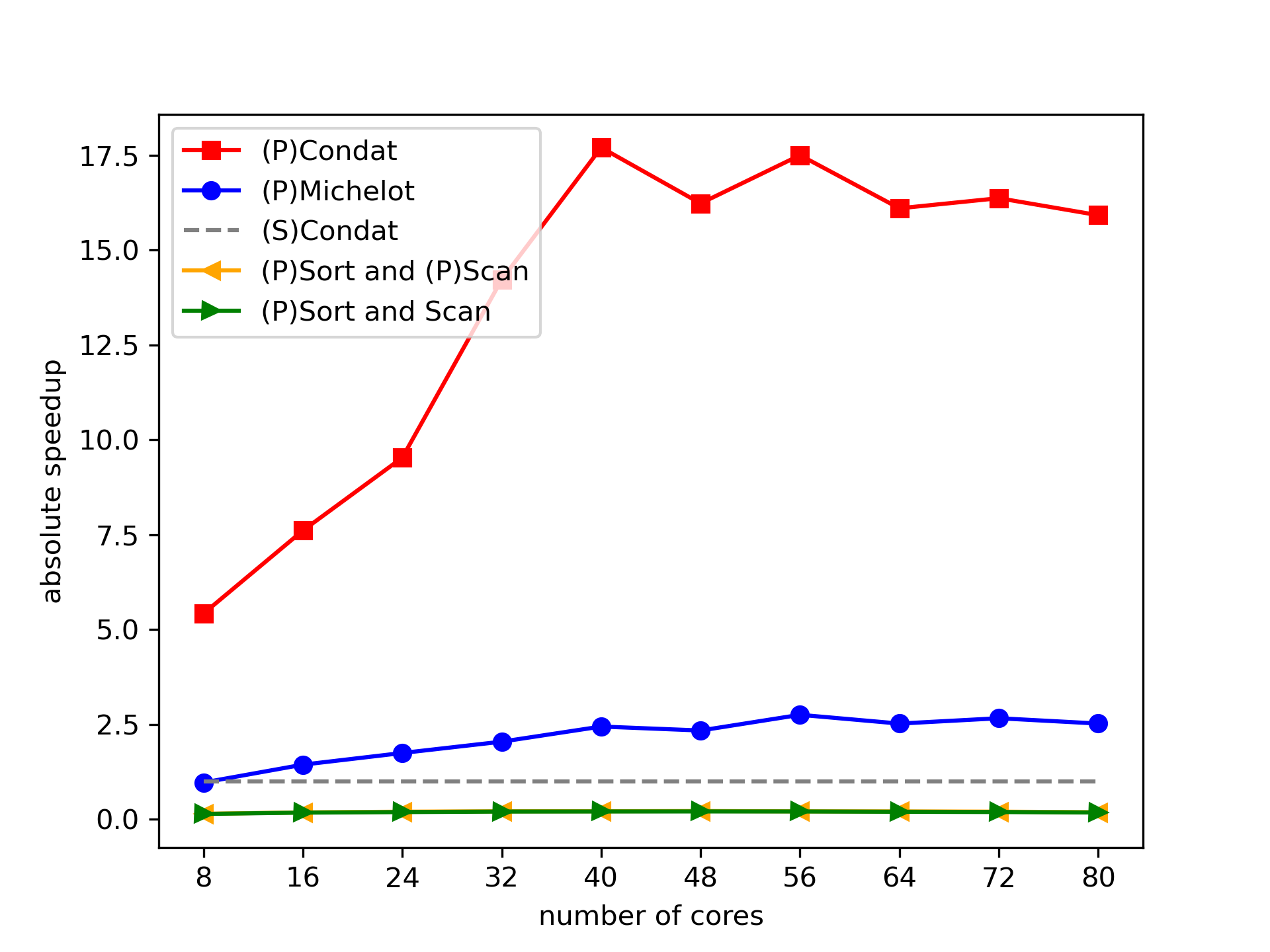}
    \end{minipage}
    \hfill
    \begin{minipage}[b]{0.45\textwidth}
        \centering
        \includegraphics[width=\textwidth]{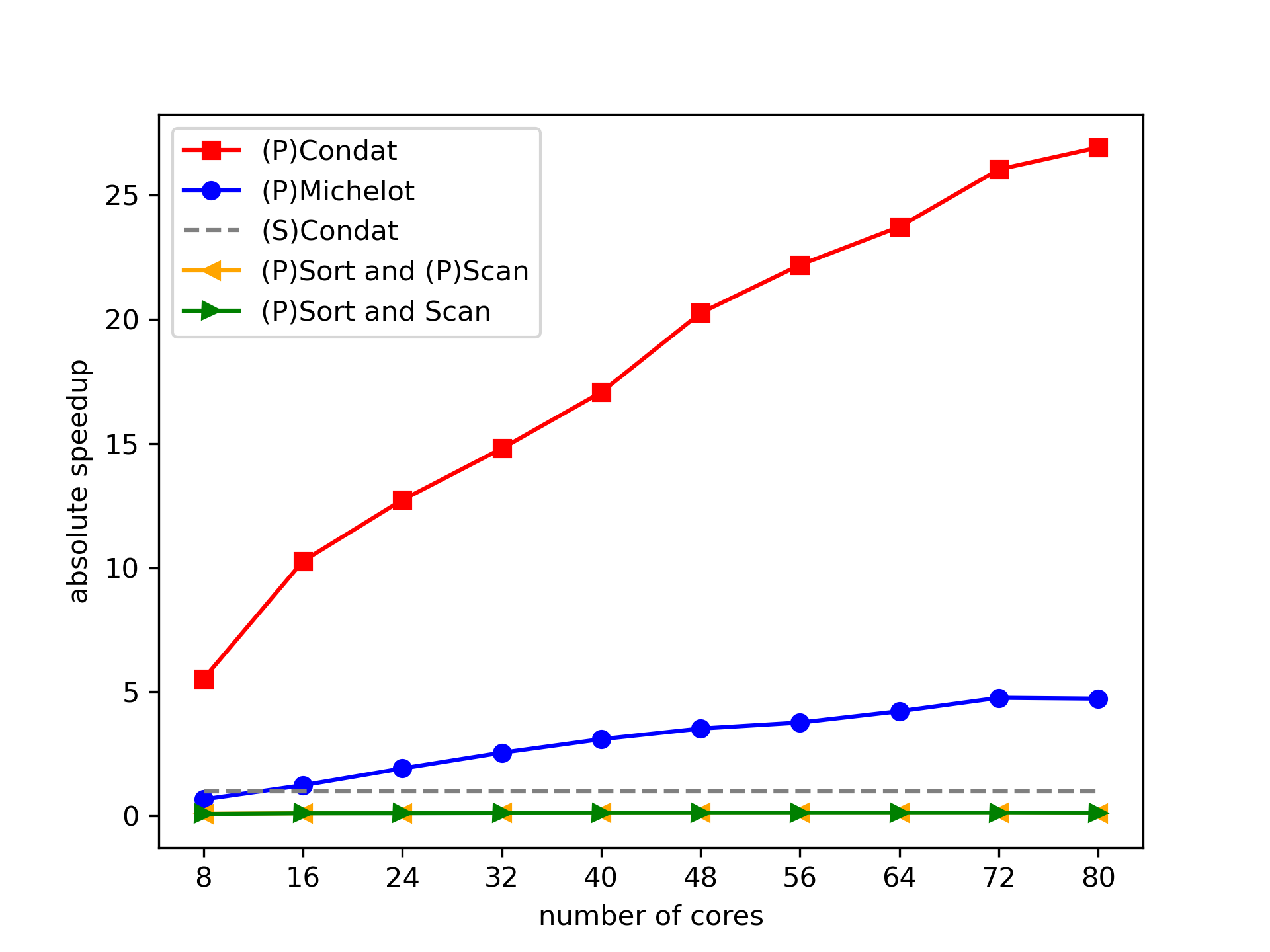}
    \end{minipage}
    \hfill
    \begin{minipage}[b]{0.45\textwidth}
        \centering
        \includegraphics[width=\textwidth]{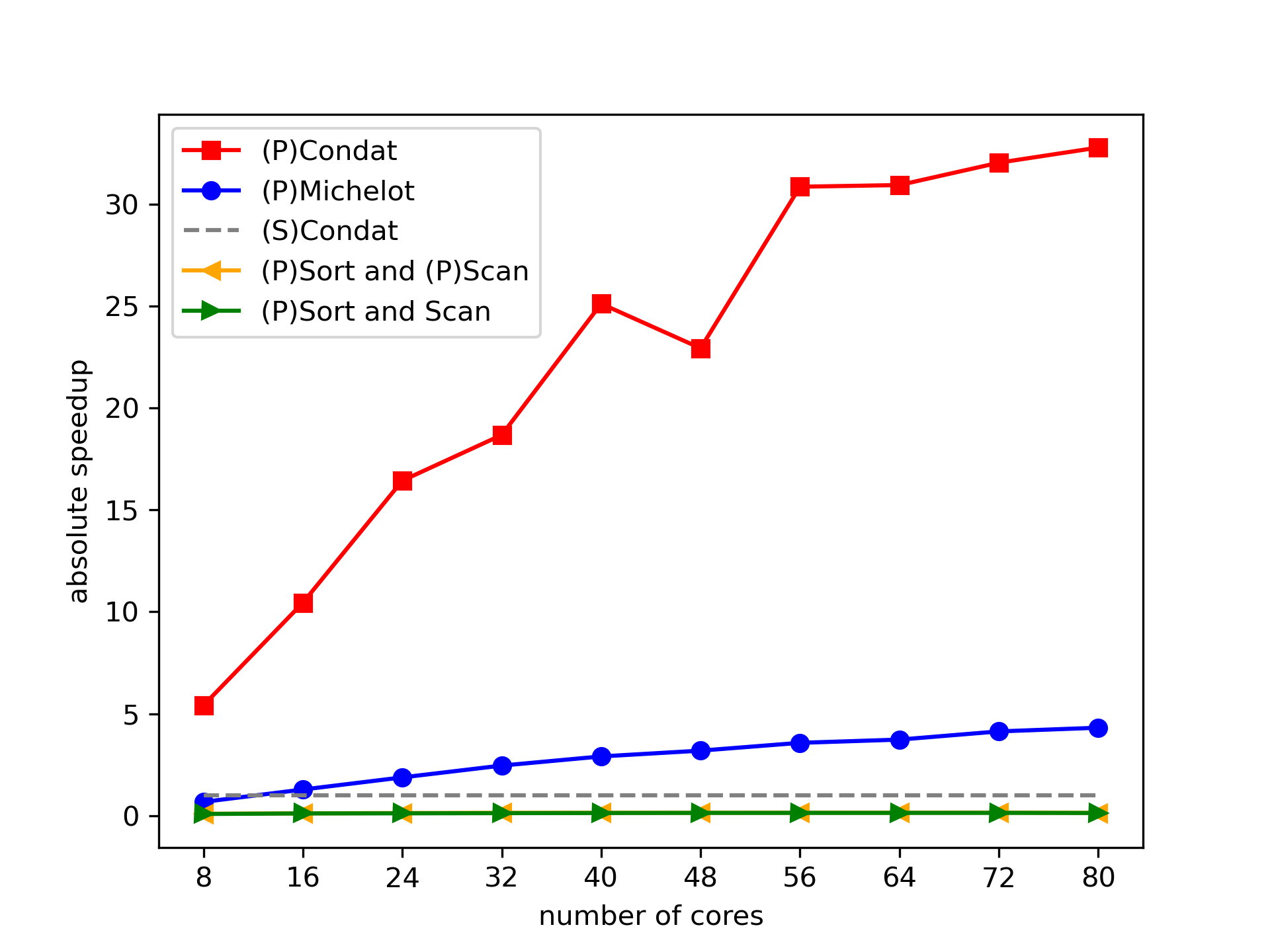}
    \end{minipage}
    \caption{Absolute speedup vs cores in simplex projection. Each line represents a different projection method, and 3 graphs represent 3 different sizes of input vector: $10^7$ (left top), $10^8$ (right top), and $10^9$ (below).}
    \label{fig:mlen_comp}
\end{figure}

\begin{figure}[!htbp]
    \centering
    \begin{minipage}[b]{0.45\textwidth}
        \centering
        \includegraphics[width=\textwidth]{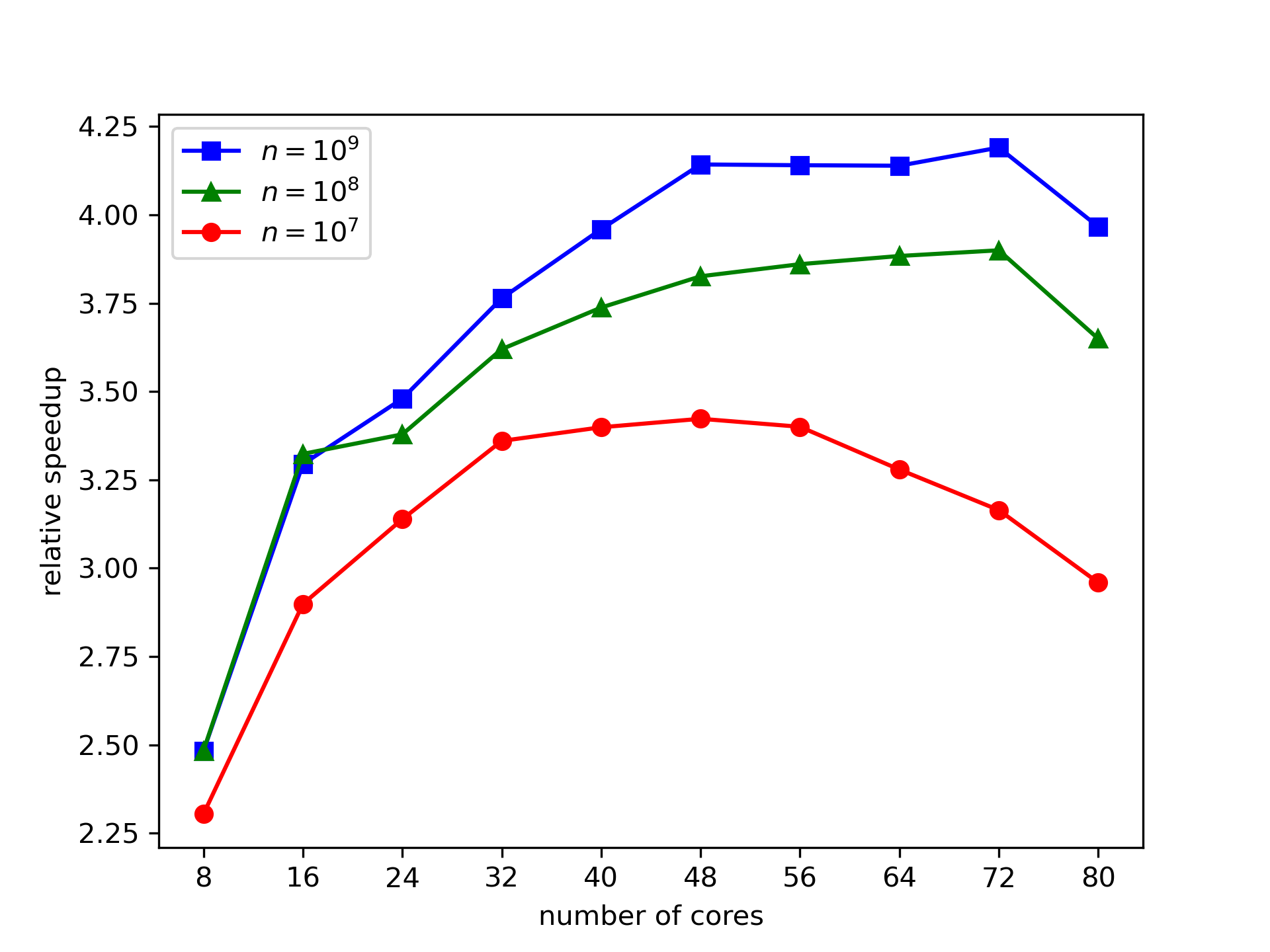}
    \end{minipage}
    \hfill
    \begin{minipage}[b]{0.45\textwidth}
        \centering
        \includegraphics[width=\textwidth]{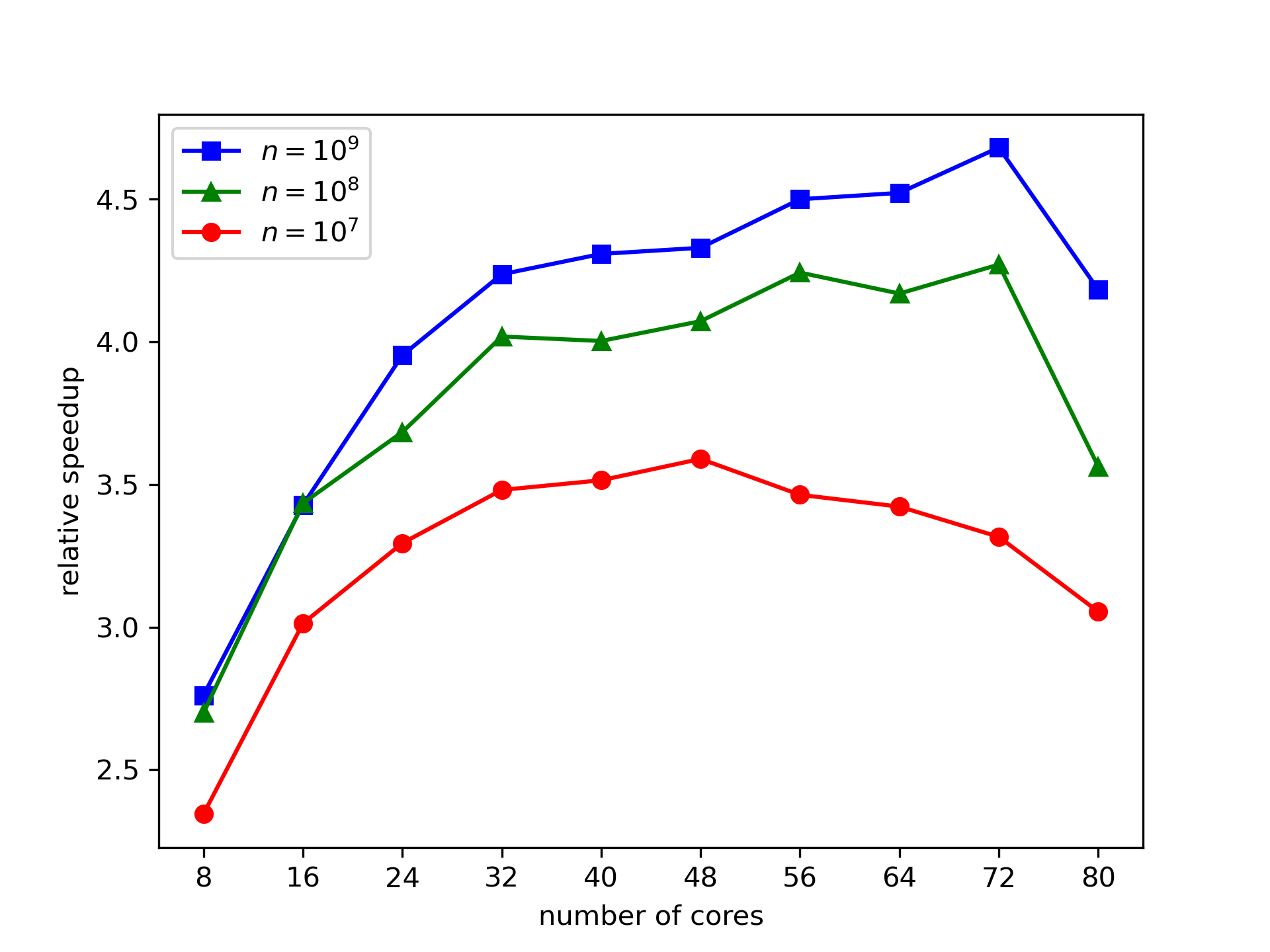}
    \end{minipage}
    \hfill
    \begin{minipage}[b]{0.45\textwidth}
        \centering
        \includegraphics[width=\textwidth]{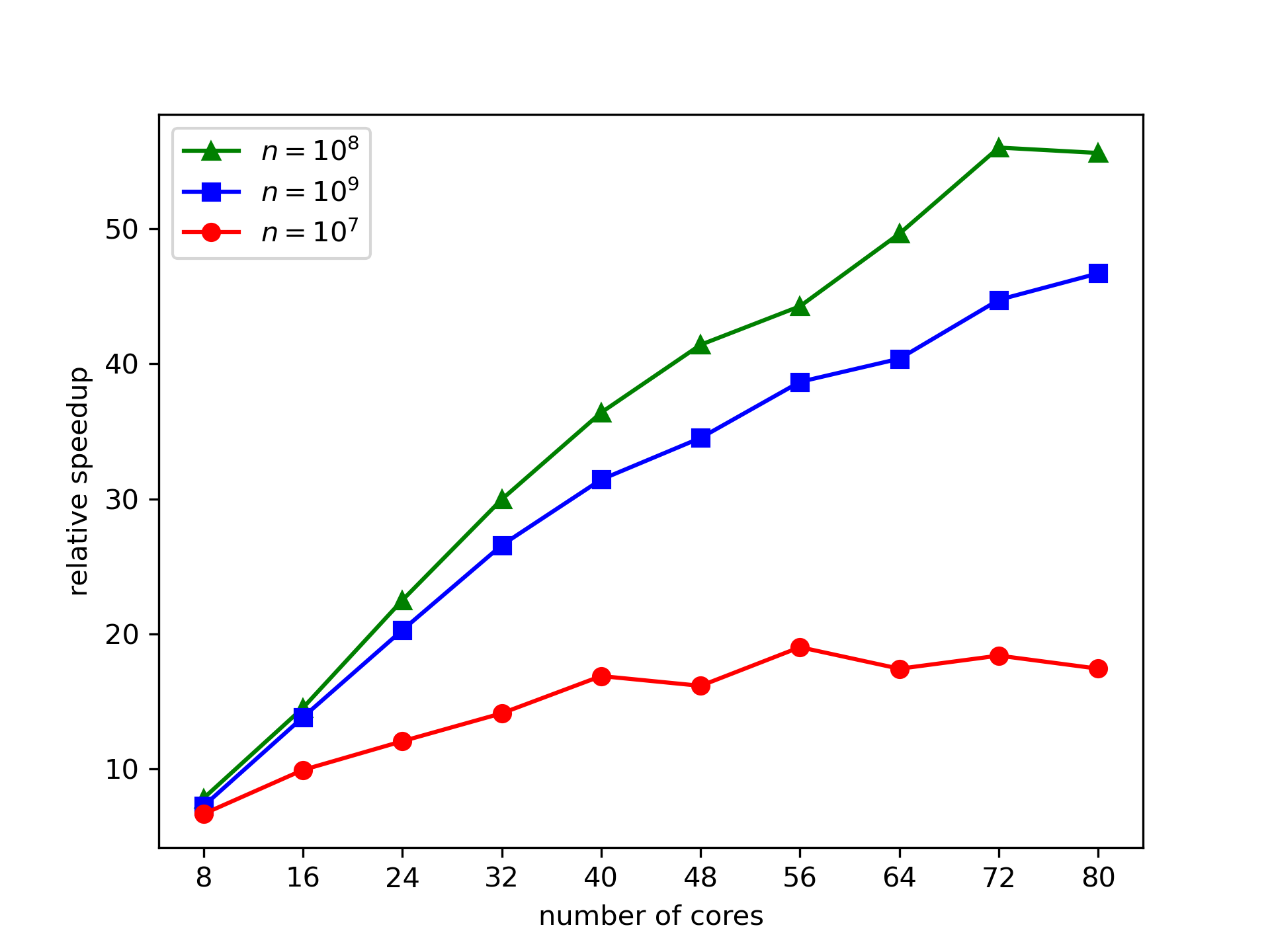}
    \end{minipage}
    \hfill
    \begin{minipage}[b]{0.45\textwidth}
        \centering
        \includegraphics[width=\textwidth]{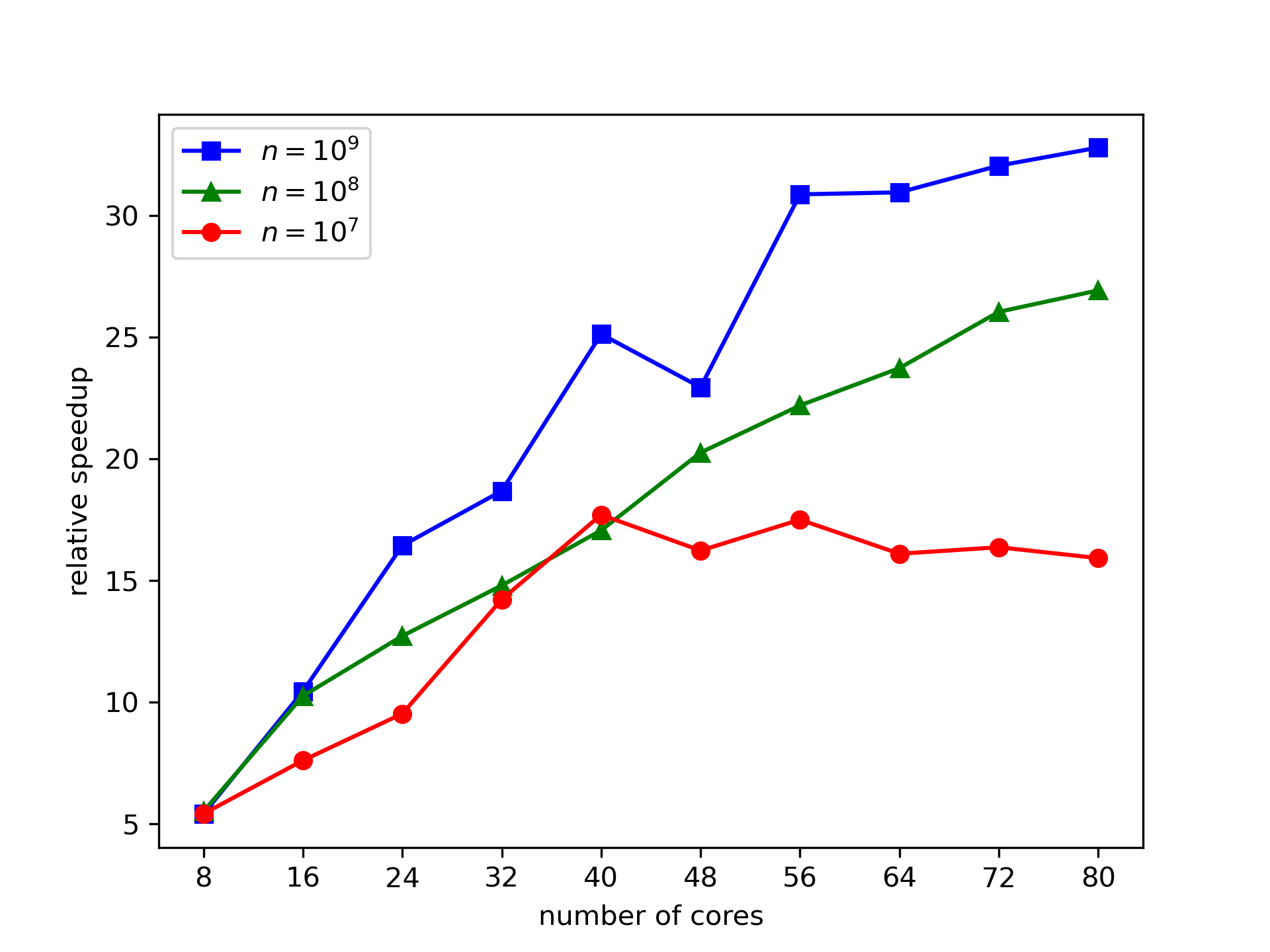}
    \end{minipage}
    \caption{Relative speedup vs cores in simplex projection. Each line represents a different size of input vector, and 4 graphs represent 4 different projection methods: Parallel Sort and Scan (left top), Parallel Sort and Partial Scan (right top), Parallel Michelot (left below), and Parallel Condat (right below).}
    \label{fig:mlen}
\end{figure}

Instances in Figures~\ref{fig:mlen_comp} and ~\ref{fig:mlen} have varying input sizes of $n = 10^7, 10^8, 10^9$, with $d_i$ are drawn i.i.d. from $N(0,1)$.  This demonstrates that the speedup per cores is a function of problem size.  At $10^7$, parallel Condat tails off in absolute speedup around 40 cores (where communication costs become marginally higher than overall gains), while consistently increasing speedups are observed up to 80 cores on the $10^9$-sized instances.  Similar patterns are observed for all algorithms in the relative speedups. For a fixed number of cores, larger instances yield larger partition sizes; hence the subvector projection problem given to each core tends to reduce more of the original vector, producing the observed effect for our parallel methods. More severe tailoff effects are observed in the Scan and Sort algorithms, which use an entirely different parallelization scheme.

For additional experiments varying $b$, please see Appendix E.2. 

\subsubsection{$\ell_{1}$ ball}
\begin{figure}[!htbp]
    \centering
    \begin{minipage}[b]{0.45\textwidth}
         \centering
         \includegraphics[width=\textwidth]{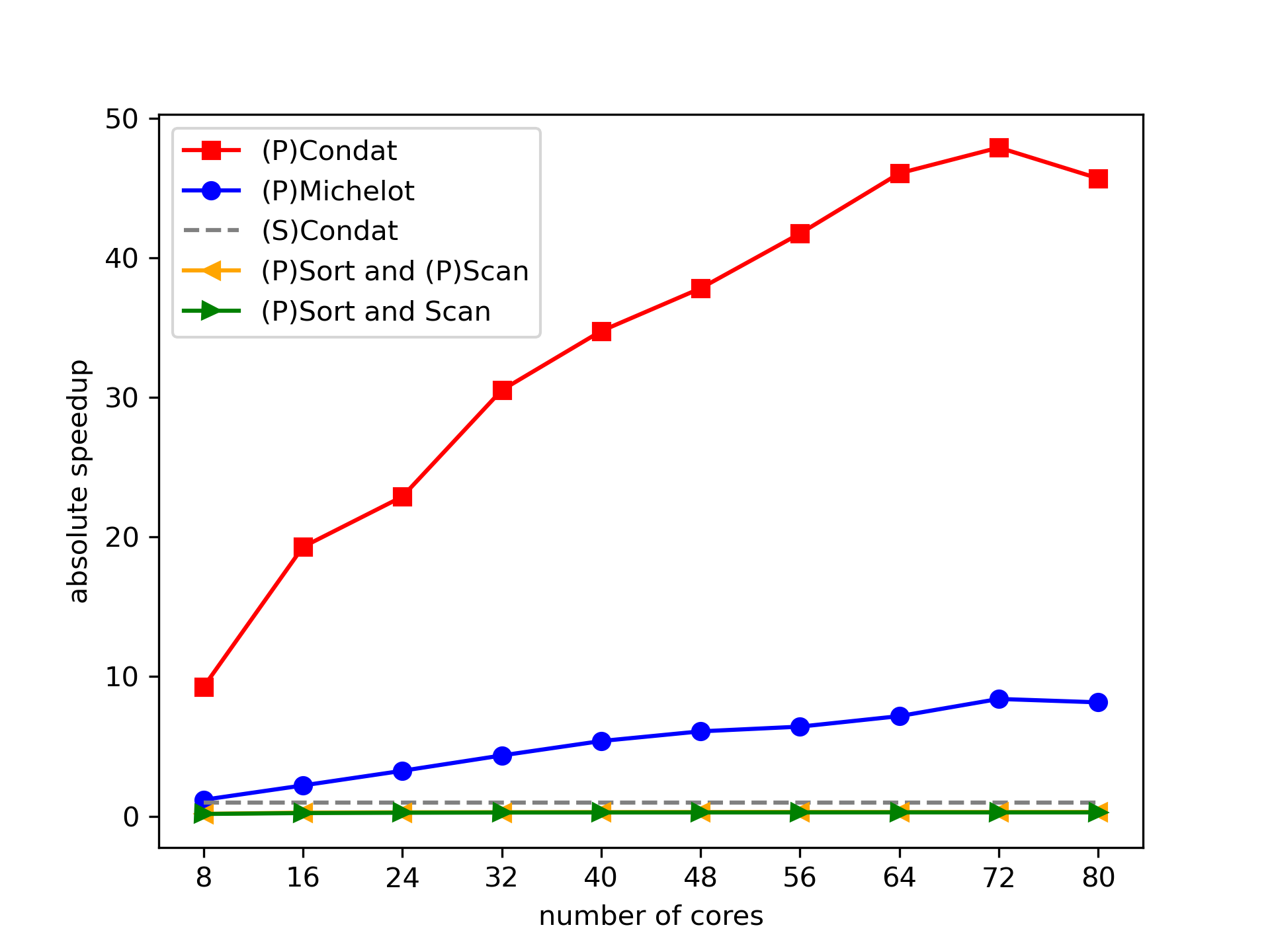}
    \end{minipage}
    \hfill
    \begin{minipage}[b]{0.45\textwidth}
         \centering
         \includegraphics[width=\textwidth]{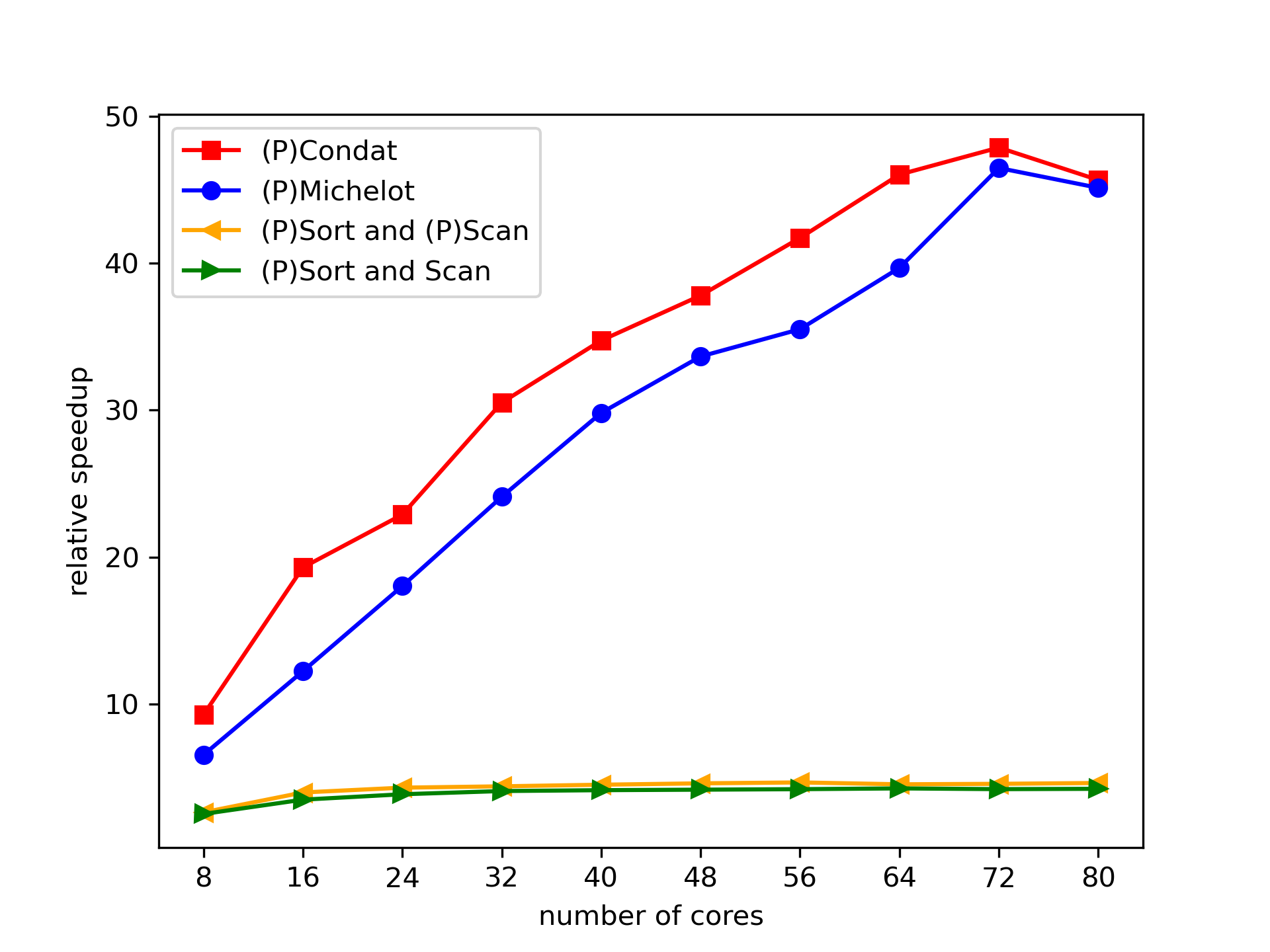}
    \end{minipage}
    \caption{Speedup vs cores in $\ell_1$ ball projection.  Each line represents a different projection method.}
    \label{fig:l1ball}
\end{figure}

We conduct $\ell_1$ ball projection experiments with problem size $n=10^{8}$. Inputs $d_i$ are drawn i.i.d. from $N(0,1)$. Algorithms were implemented as described in Sec~\ref{sec:ell1desc}. Results are shown in Figure~\ref{fig:l1ball}. Similar to the standard projection onto Simplex problems, our parallel Condat implementation attains considerably superior results over the benchmark, with nearly 50X speedup.

\subsubsection{Weighted Simplex and Weighted $\ell_{1}$ ball}
We have conducted additional experiments using the weighted versions of simplex and $\ell_{1}$ ball projections.  These have been placed in the online supplement, as results are similar. A description of the algorithms are given in Appendix B.2; pseudocode in Appendix D; and experimental results in Appendix E.3.  

\subsubsection{Parity Polytope}
\begin{figure}[!htbp]
    \centering
    \begin{minipage}[b]{0.45\textwidth}
         \centering
         \includegraphics[width=\textwidth]{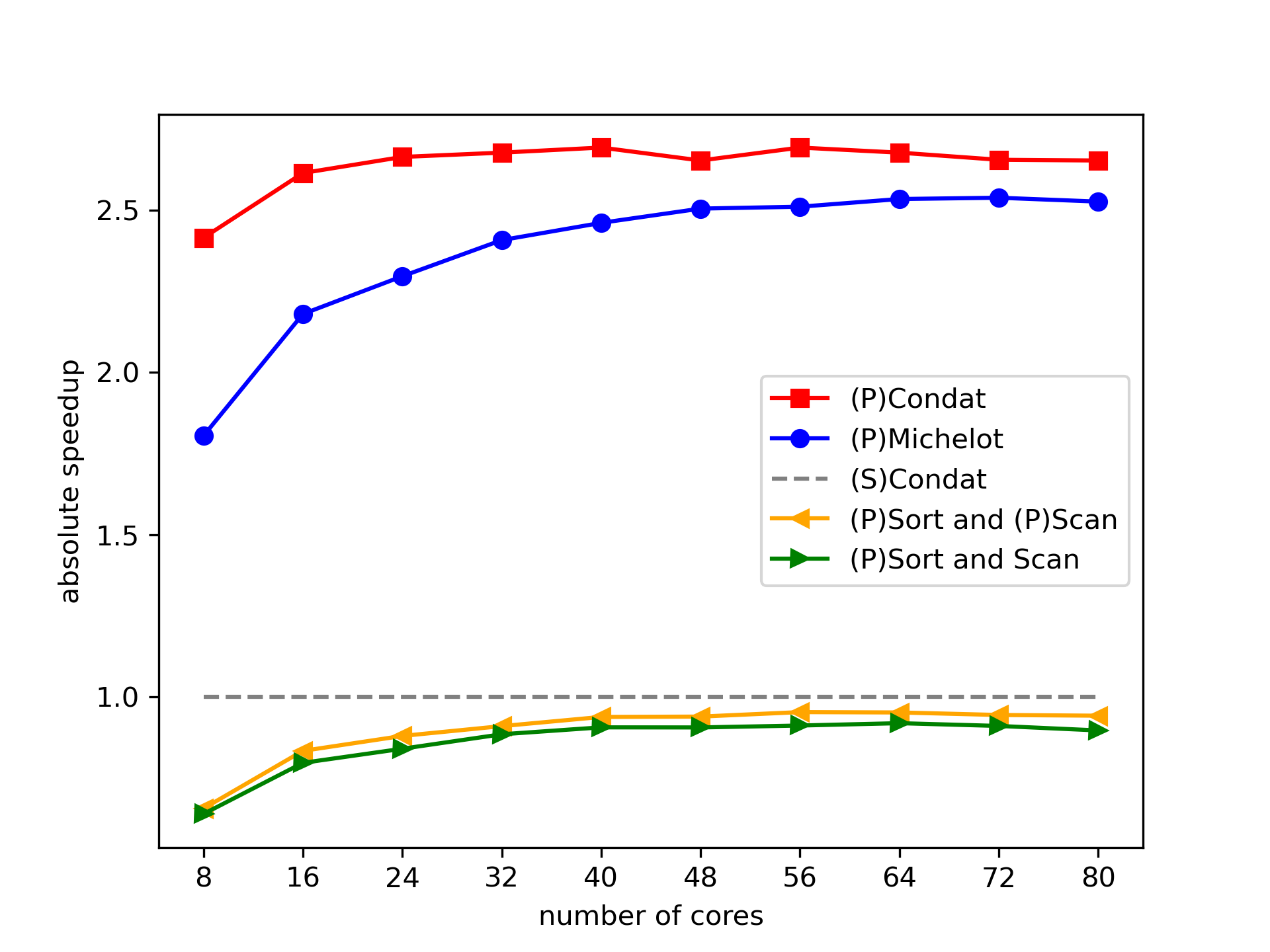}
    \end{minipage}
    \hfill
    \begin{minipage}[b]{0.45\textwidth}
         \centering
         \includegraphics[width=\textwidth]{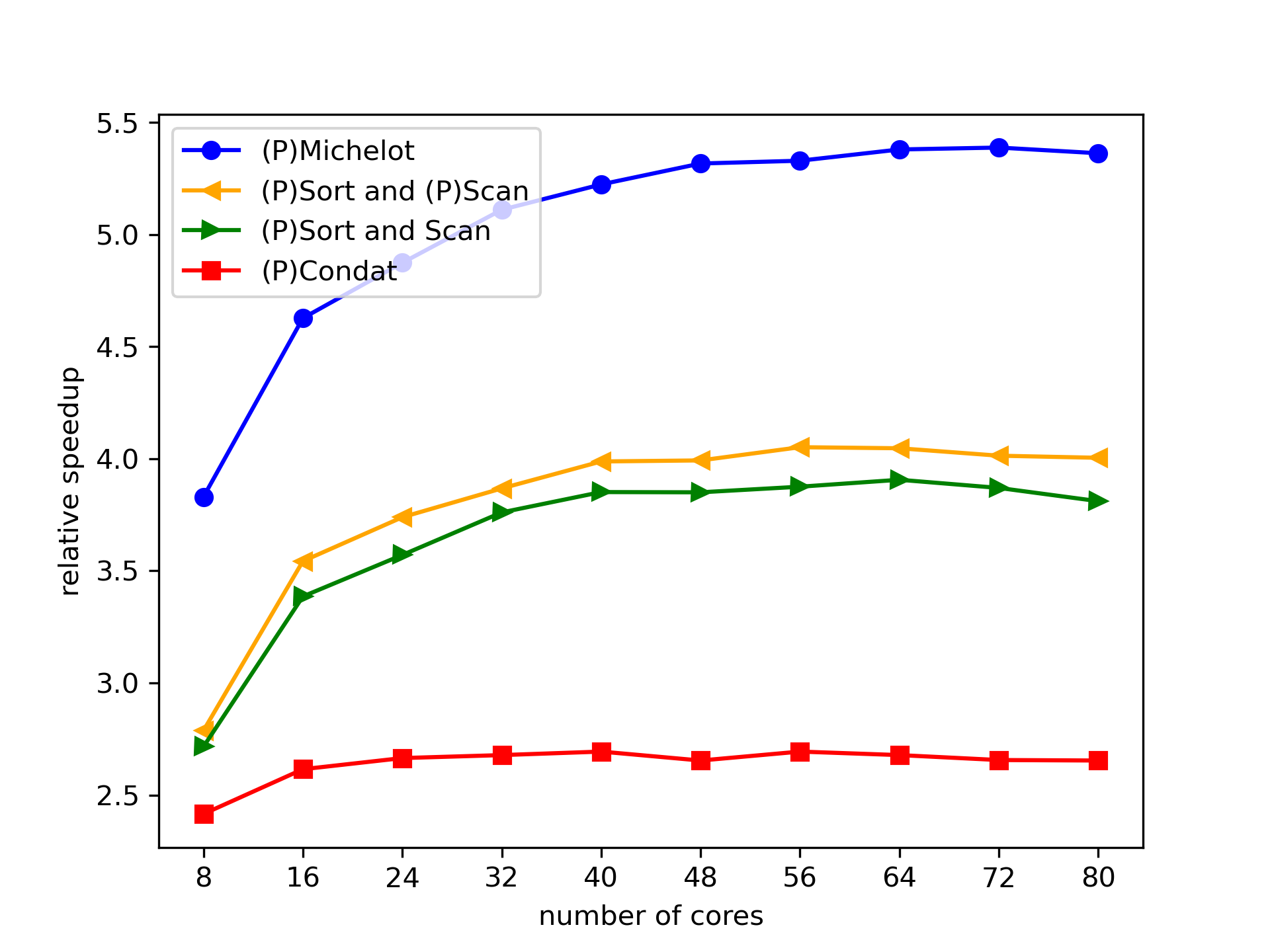}
    \end{minipage}
    \hfill
    \begin{minipage}[b]{0.45\textwidth}
         \centering
         \includegraphics[width=\textwidth]{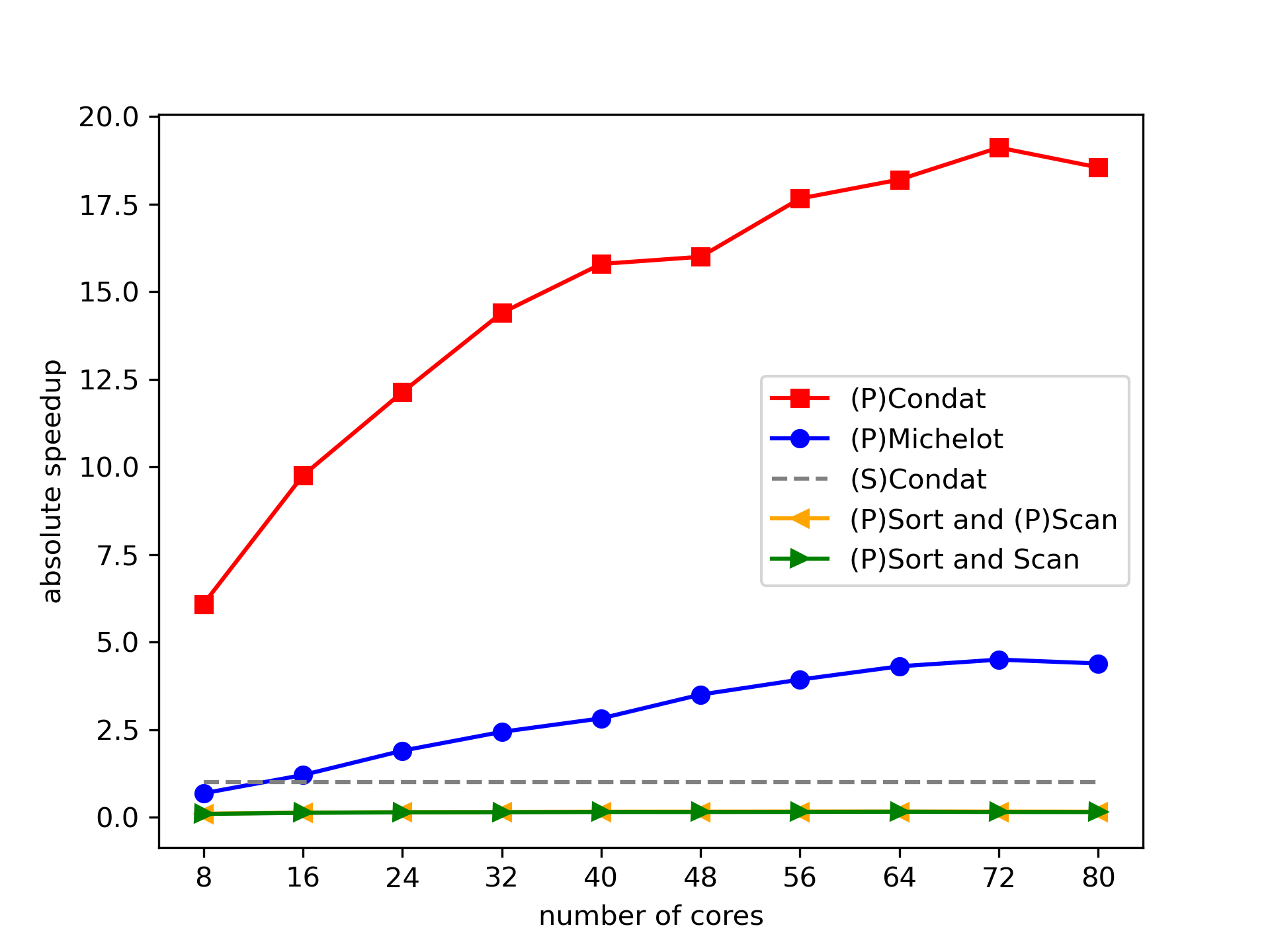}
    \end{minipage}
    \hfill
    \begin{minipage}[b]{0.45\textwidth}
         \centering
         \includegraphics[width=\textwidth]{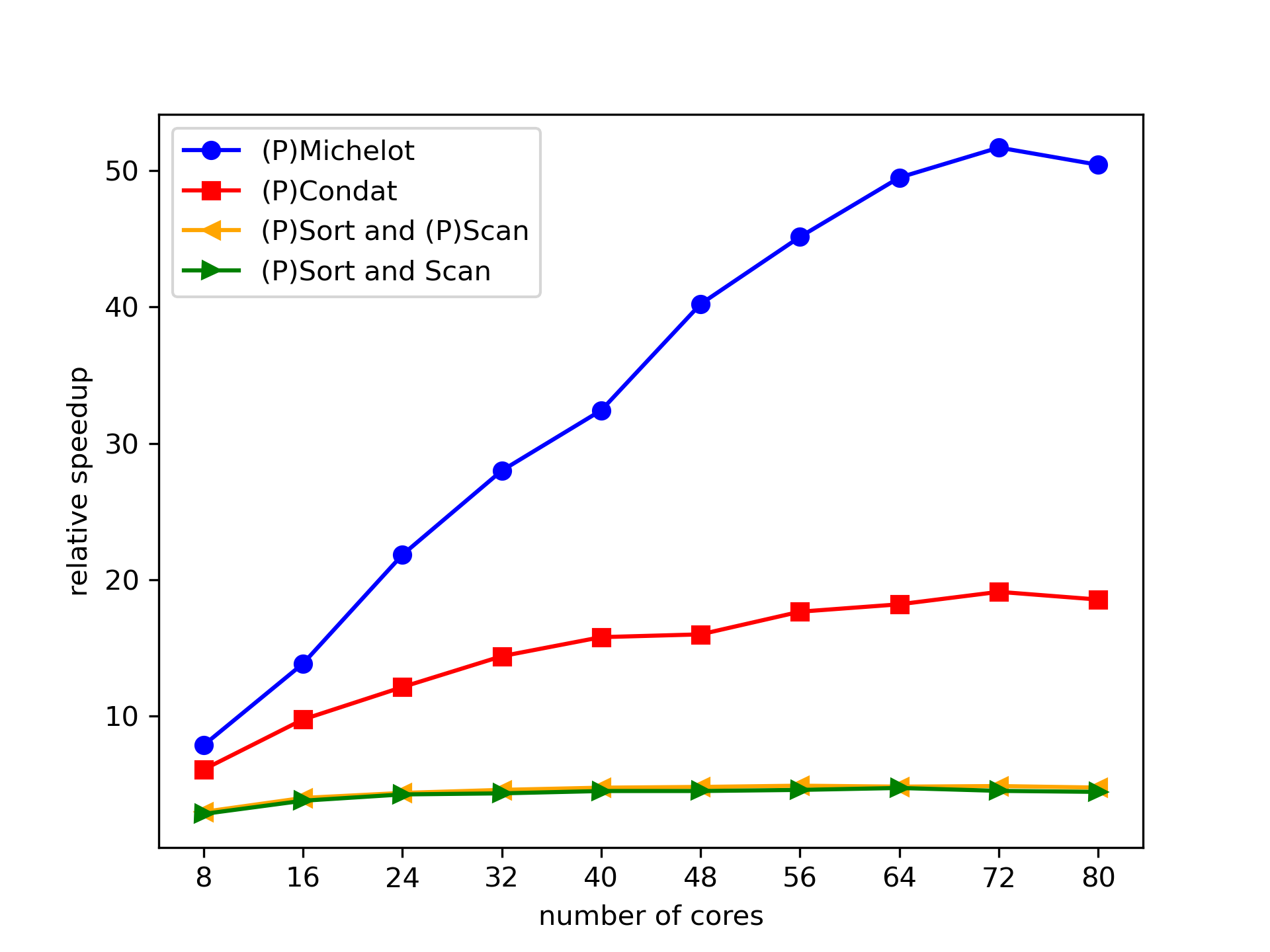}
    \end{minipage}
    \caption{Speedup vs cores in projection onto a parity polytope. Each line represents a different projection method, and 4 graphs represent 4 different experiments: absolute speedup and relative speedup in parity polytope projection (top line), absolute and relative speedup only for the simplex projection component of parity polytope projection (below line).}
    \label{fig:parity_simplex_proj}
\end{figure}

We conduct parity polytope projection experiments with a problem size setting of $n=10^{8}-1$. Inputs $d_i$ are drawn i.i.d. from $U[1,2]$. Algorithms were implemented as described in Sec~\ref{sec:cpppdesc}. Results are shown in Figure~\ref{fig:parity_simplex_proj}. Our parallel Condat has worse relative speedup on projections compared to parallel Michelot, but overall this still results in higher absolute speedups since the baseline serial Condat runs quickly.  With up to around 20x absolute speedup in simplex projection subroutines from parallel Condat's, we report an overall absolute speedup of up to around 2.75x for parity polytope projections. We note that the overall absolute speedup tails off more quickly; this is an expected effect from Amdahl's law \cite{amdahl1967validity}: diminishing returns due to partial parallelization of the algorithm.

\subsubsection{Lasso on Real-World Data}
We selected a dataset from a paper implementing a Lasso method \citep{realdata}: \emph{kdd2010} (named as \emph{kdda} in the cited paper), and its updated version \emph{kdd2012}; both of them can be found in LIBSVM data sets \citep{libsvm}. There are $n=20,216,830$ features in \emph{kdd2010} and $n=54,686,452$ features in \emph{kdd2012}. 

\begin{figure}[!htbp]
    \centering
    \begin{minipage}[b]{0.45\textwidth}
         \centering
         \includegraphics[width=\textwidth]{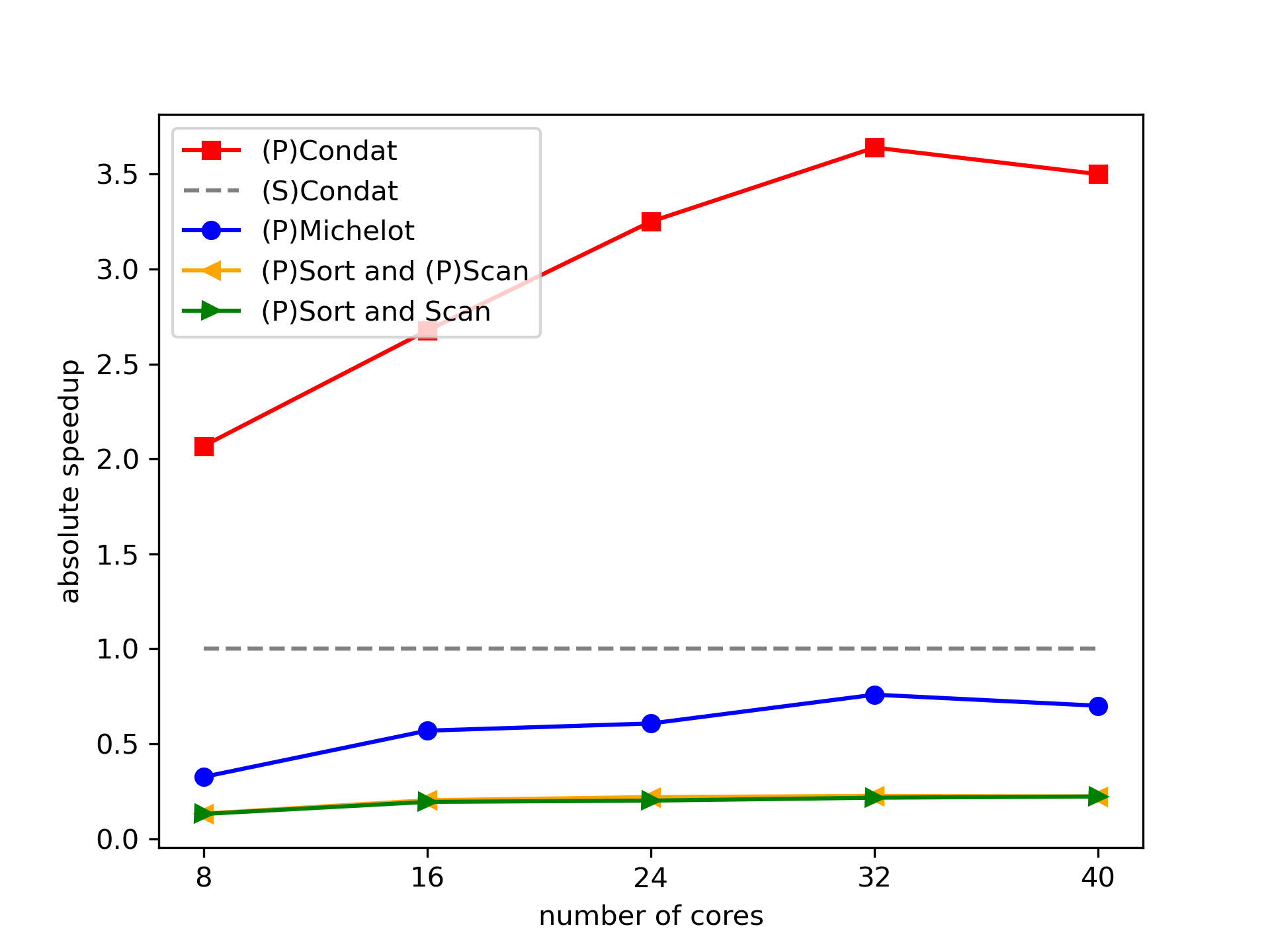}
    \end{minipage}
    \hfill
    \begin{minipage}[b]{0.45\textwidth}
         \centering
         \includegraphics[width=\textwidth]{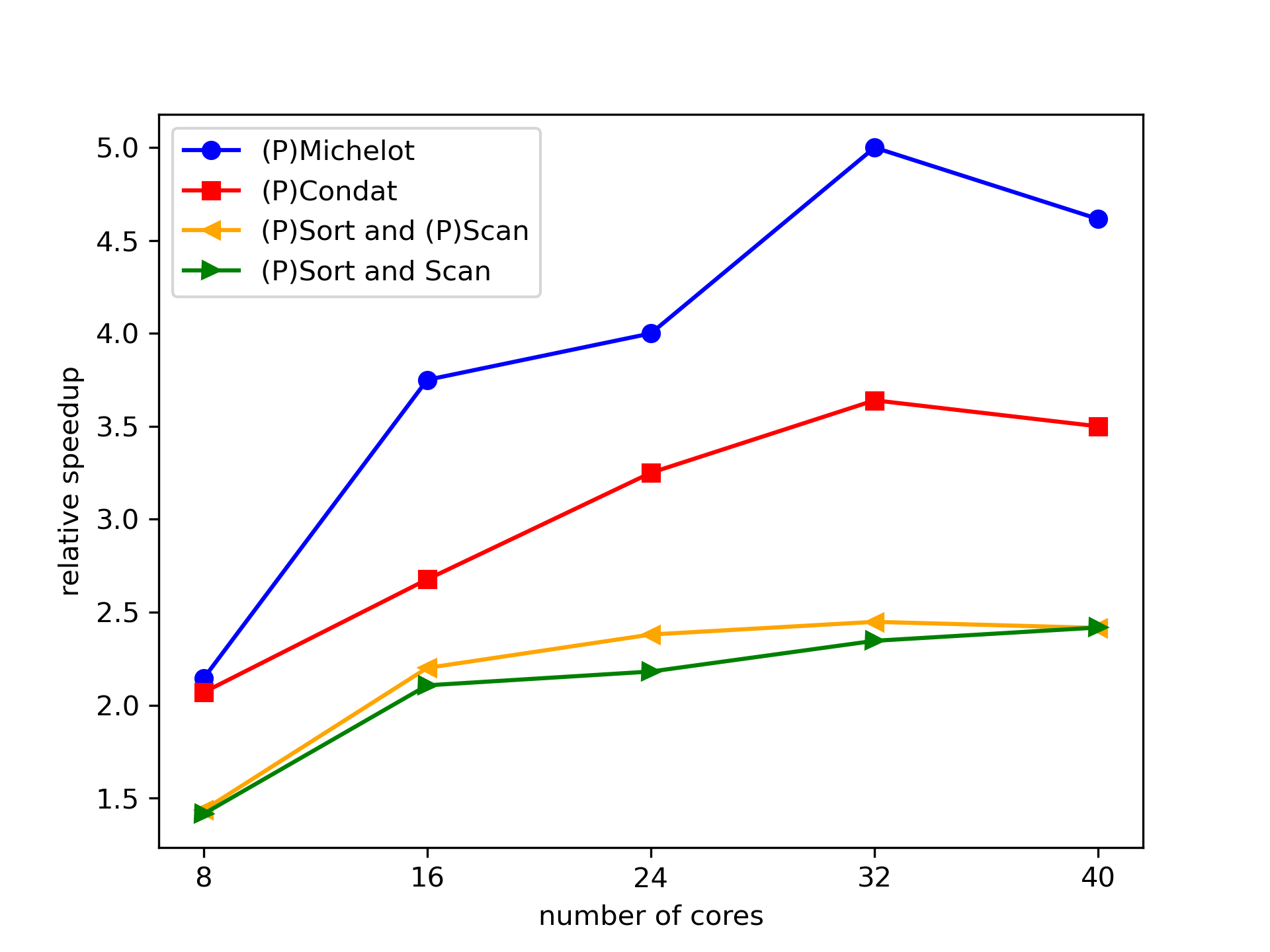}
    \end{minipage}
    \caption{Speedup vs cores of $\ell_1$ ball Projection in Lasso on \emph{kdd2010}. Each line represents a different projection method.}
    \label{fig:kdd10}
\end{figure}
\begin{figure}[!htbp]
    \centering
    \begin{minipage}[b]{0.45\textwidth}
         \centering
         \includegraphics[width=\textwidth]{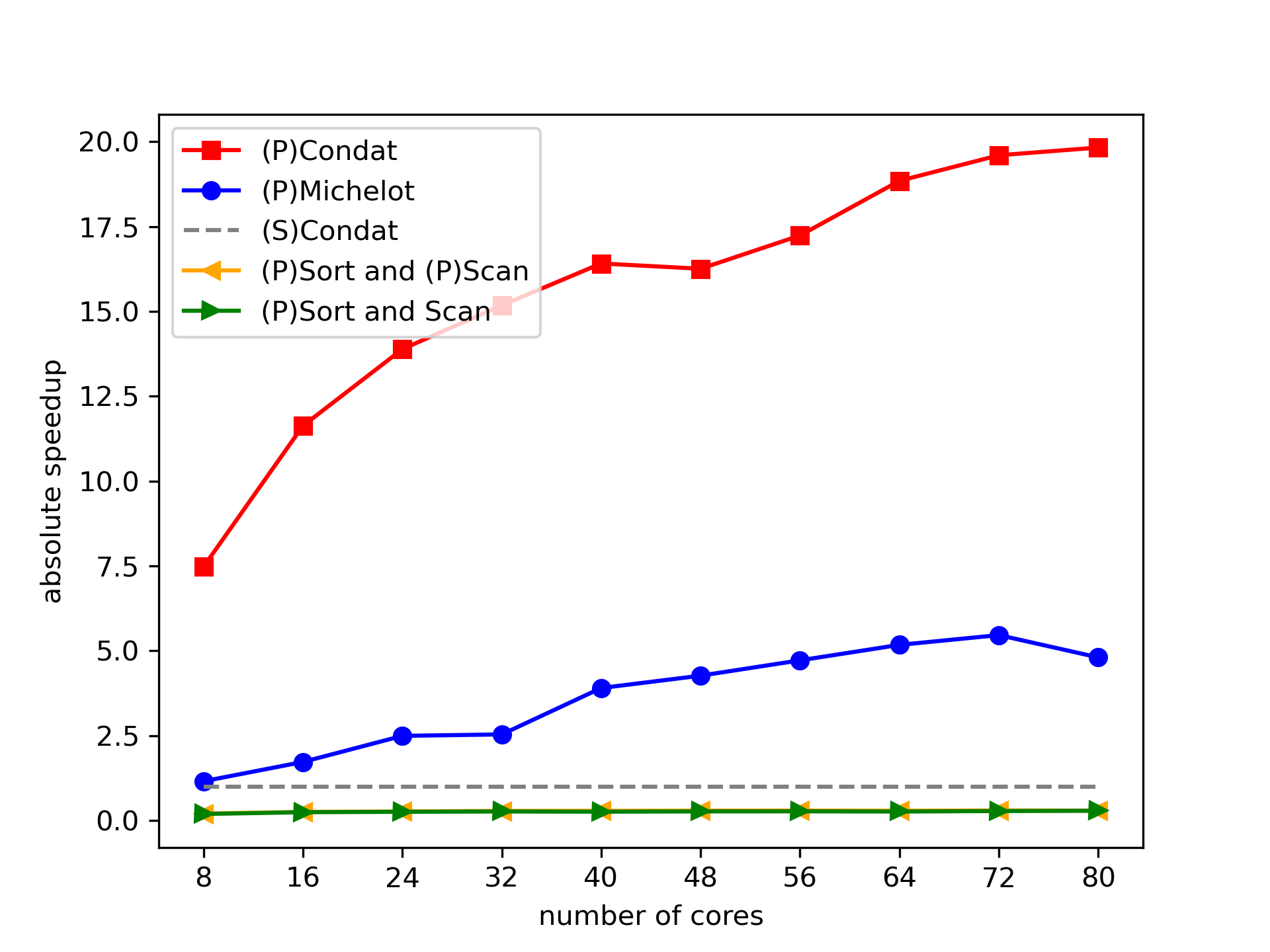}
    \end{minipage}
    \hfill
     \begin{minipage}[b]{0.45\textwidth}
         \centering
         \includegraphics[width=\textwidth]{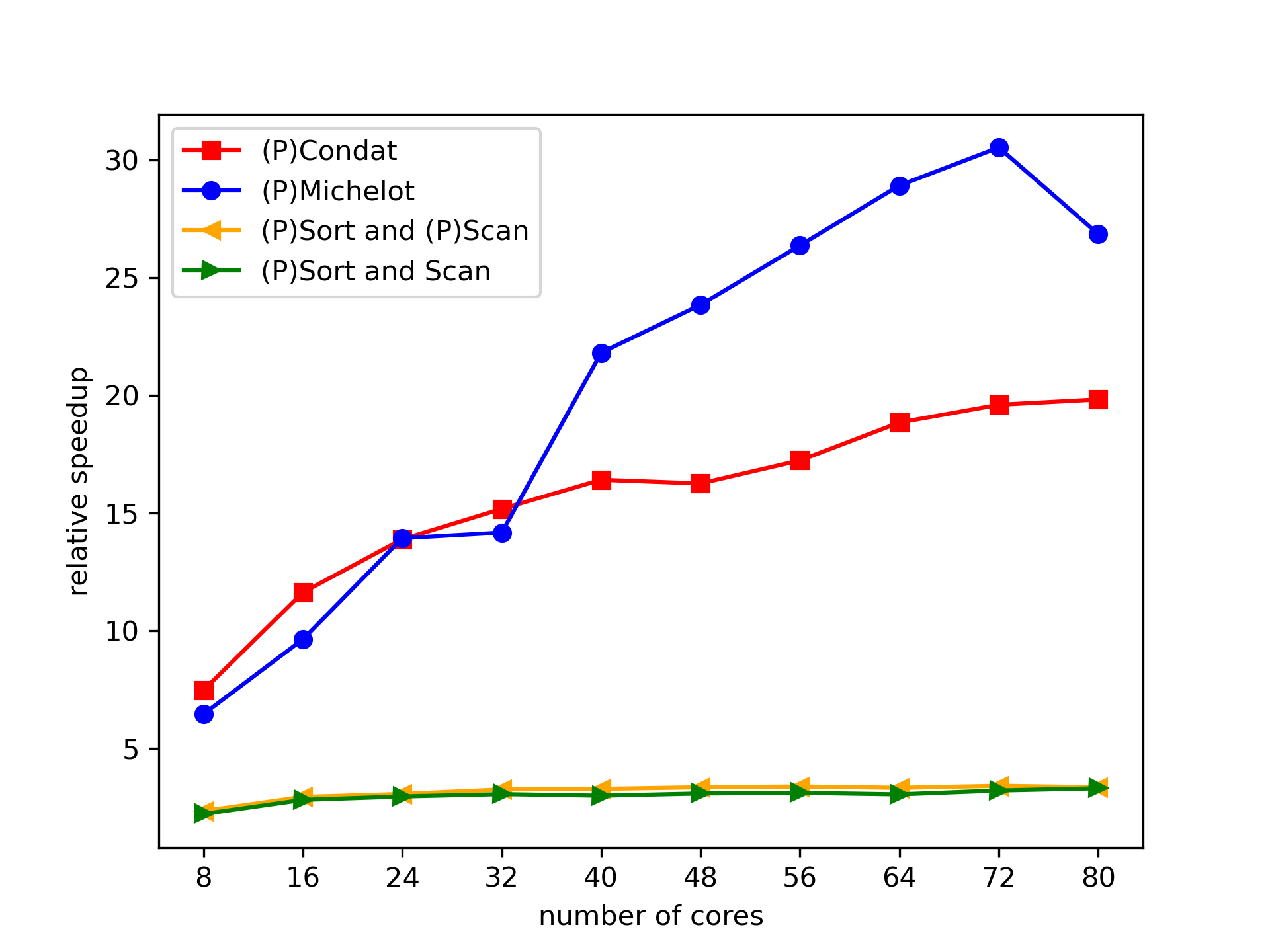}
    \end{minipage}
    \caption{Speedup vs cores of $\ell_1$ ball Projection in Lasso on \emph{kdd2012}. Each line represents a different projection method.}
    \label{fig:kdd12}
\end{figure}

We implemented PGD with Mini-batch \citep[Sec. 12.5]{minibatch}, and measure the runtime for subroutine of projecting onto $\ell_1$ ball, which is
\[x_{t+1} = \mbox{proj}_{\mathcal{B}_1}(x_t-\alpha \cdot 2\tilde{A}^{T}(\tilde{A}x_t-\beta)),\]
where the initial point $x_0$ is drawn \texttt{i.i.d.} from sparse $U[0,1]$ with $0.5$ sparse rate, $\tilde{A}\in \mathbb{R}^{m\times n}$ with sample size $m$ and $n$ features, and $\beta$ includes labels for samples; moreover, $\alpha=0.05$. Moreover, we set $m=128$ samples in each iteration. We  measure the runtime for $\ell_1$ ball projections for the first $10$ iterations. In later iterations the projected vector $x$ will be dense, at which point it is better to use serial projection methods (see Appendix E.5). Runtimes are measured by \texttt{time\_ns()} Function), and the absolute speedup and relative speedup are reported in Figure~\ref{fig:kdd10} (for \emph{kdd2010}) and Figure~\ref{fig:kdd12} (for \emph{kdd2012}). Considerably more speedup is obtained for \emph{kdd2012}, which may be due to problem size since the projection input vectors are more than twice the size of inputs from \emph{kdd2010}. 

\subsubsection{Discussion}
We observed consistent patterns of performance across a wide range of test instances, varying in size, distributions, and underlying projections. In terms of relative speedups, our parallelization scheme was surprisingly at least as effective as that of Sort and Scan. Parallel sort and parallel scan are well-studied problems in parallel computing, so we expected \emph{a priori} for such speedups to be a strict upper bound on what our method could achieve. Thus our method is not simply benefiting from its compatibility with more advanced serial projection algorithms---the distributed method itself appears to be highly effective in practice.

\section{Conclusion}
\label{sec:conclusion}
We proposed a distributed preprocessing method for projection onto a simplex. Our method distributes subvector projection problems across processors in order to reduce the candidate index set on large-scale instances. One advantage of our method is that it is compatible with all major serial projection algorithms. Empirical results demonstrate, across a wide range of simulated distributions and real-world instances, that our parallelization of well-known serial algorithms are comparable and at times superior in relative speedups to highly developed and well-studied parallelization schemes for sorting and scanning. Moreover, the sort-and-scan serial approach involves substantially more work than e.g. Condat's method; hence our parallelization scheme provides considerable absolute speedups in our experiments versus the state of the art.

The effectiveness depends on the sparsity of the projection; in the case of large-scale problems with i.i.d. inputs and $b \in o(n)$, we can expect high levels of sparsity. A wide range of large-scale computational experiments demonstrates the consistent benefits of our method, which can be combined with any serial projection algorithms. Our experiments on real-world data suggest that significant sparsity can be exploited even when such distributional assumptions could be violated.  We also note that, due to Proposition~\ref{prop: define tau} and Corollary~\ref{cor: tau & t}, highly dense projections occur when $\tau$ has a low value relative to the entries of $d$.  Now, iterative (serial) methods as Michelot's and Condat's can be interpreted as starting with a pivot value that is a lower bound on $\tau$ and increasing the pivot iteratively until the true value $\tau$ is attained.  Hence, when our distributed method performs poorly due to density of the projection, the problem can simply be solved in serial using a small number of iterations, and vice-versa. This might be expected, for instance, when the input vector $d$ itself is sparse (see Appendix E.5). 

\section{Acknowledgements}
This work was funded in part by the Office of Naval Research under grant N00014-23-1-2632. We also thank an anonymous reviewer for various code optimization suggestions.

\newpage
\bibliography{ref.bib} 
\bibliographystyle{informs2014.bst}

\newpage
\appendix
\section{Mathematical Proofs}
\textbf{Corollary 1}
    For any $t_1,t_2 \in \mathbb{R}$ such that $t_1 < \tau < t_2$, we have
    \[f(t_1) > f(\tau) = 0 > f(t_2). \]

\begin{proof}
    By definition,
    \begin{alignat*}{3}
    f(t) &=&& \frac{\sum_{i\in \mathcal{I}_t}d_{i}-b}{|\mathcal{I}_t|}-t,\\
    &=&& \frac{\sum_{i\in \mathcal{I}_t}(d_{i}-t)-b}{|\mathcal{I}_t|},\\
    &=&& \frac{\sum_{i=1}^n \max\{d_{i}-t,0\}-b}{|\mathcal{I}_t|}.\\
    \end{alignat*}
	Observe that $g(t) := \sum_{i=1}^n \max\{d_{i}-t,0\}-b$ is a strictly decreasing function for $t\leq \max_i\{d_i\}$, and $g(t) = -b$ for $t\geq \max_i\{d_i\}$.  Furthermore, from Proposition 1, we have that $\tau$ is the unique value such that $g(\tau) = 0$; moreover,  since $b>0$ then $\tau < \max_i\{d_i\}$. Thus $g(t_1)>g(\tau)=0>g(t_2)$, which implies $f(t_1)>0, f(\tau)=0$, and $f(t_2)<0$.
\end{proof}
We assume uniformly distributed inputs, $d_{1},\dots,d_{n}$ are $\mathrm{i.i.d} \sim U[l,u]$, and we have\\
\textbf{Lemma 1}
\[E[\frac{\sum_{i\in\mathcal{I}}d_{i}-b}{|\mathcal{I}|}] \rightarrow \frac{l+u}{2}\ \mathrm{as}\ n \rightarrow \infty\]
    with a sublinear convergence rate, where $\mathcal{I}:=\{1,...,n\}$.

\begin{proof}
	Observe that $E[d_{i}]=\frac{l+u}{2}$. Since $d_1,...d_n$ are i.i.d, then we have
	\begin{equation}\label{eq:le1}
	    E[\frac{\sum_{i\in\mathcal{I}}d_{i}-b}{|\mathcal{I}|}]=\frac{l+u}{2}-\frac{b}{n}.
	\end{equation}
	Thus $E[\frac{\sum_{i\in\mathcal{I}}d_{i}-b}{|\mathcal{I}|}]$ converges to $\frac{l+u}{2}$ sublinearly at the rate of 
 \[\lim_{n\rightarrow \infty} [\frac{l+u}{2}-\frac{b}{n+1}]/[\frac{l+u}{2}-\frac{b}{n}]=1.\]
\end{proof}
~\\
Let $X|t$ denote the conditional variable such that $P_{X|t}(x)=P(X=x|X>t)$.
~\\
We assume uniformly distributed inputs, $d_{1},\dots,d_{n}$ are $\mathrm{i.i.d} \sim U[l,u]$, and then\\
~\\
\textbf{Proposition 2} Michelot's method has an average runtime of $O(n)$.
\begin{proof}
	 Let $\delta_i$ be the number of elements that Algorithm 3 (from the main body) removes from the (candidate) active set $\mathcal{I}_p^{(i)}$ in iteration $i$ of the do-while loop, and let $T$ be the total number of iterations. 
	\begin{equation}\label{eq:le2_1}
	    \sum_{i=1}^{T}\delta_i=n-|\mathcal{I}_\tau|,
	\end{equation}
	where $\mathcal{I}_\tau$ is the active index set of the projection.
	Let $p^{(i)}$ be the $i$th pivot (line 3 in Algorithm 3 of the main body) with
	\[p^{(i)}=\frac{\sum_{j\in\mathcal{I}_{p}^{(i)}}d_j-b}{|\mathcal{I}_{p}^{(i)}|},\]
	and define $p^{(0)} := l$.

Now, from Proposition 9 we have that all $d_j$ with $j\in \mathcal{I}_{p}^{(i)}$ are i.i.d.$\sim U[p^{(i-1)},u]$. Thus $E[d_j|d_j \in \mathcal{I}_p^{(i)}]=\frac{p^{(i-1)}+u}{2}$, and so
	\begin{equation}\label{eq:le2_2}
	    E[p^{(i)}|\mathcal{I}_{p}^{(i)},p^{(i-1)}]=E[\frac{\sum_{j\in\mathcal{I}_{p}^{(i)}}d_j-b}{|\mathcal{I}_{p}^{(i)}|}\ |\ \mathcal{I}_{p}^{(i)},p^{(i-1)}] = \frac{p^{(i-1)}+u}{2}-\frac{b}{|\mathcal{I}_{p}^{(i)}|}.
	\end{equation}
	In iteration $i$, all elements in the range $[p^{(i-1)},p^{(i)}]$ are removed, and again the i.i.d. uniform property is preserved by Proposition 9, so
	\[E[\delta_{i}|\mathcal{I}_p^{(i)},p^{(i)},p^{(i-1)}] = |\mathcal{I}_p^{(i)}|\frac{p^{(i)}-p^{(i-1)}}{u-p^{(i-1)}}=|\mathcal{I}_p^{(i)}|(\frac{1}{2}-\frac{(u+p^{(i-1)})/2-p^{(i)}}{u-p^{(i-1)}});\]
	by Law of Total Expectation,
	\begin{alignat*}{3}
	    E[\delta_i|\mathcal{I}_p^{(i)},p^{(i-1)}]
	    =&E[E[\delta_{i}|\mathcal{I}_p^{(i)},p^{(i)},p^{(i-1)}]]\\
	    =&|\mathcal{I}_p^{(i)}|(\frac{1}{2}-\frac{(u+p^{(i-1)})/2-E[p^{(i)}|\mathcal{I}_p^{(i)},p^{(i-1)}]}{u-p^{(i-1)}}).
	\end{alignat*}
	Replacing the right-hand side expectation using Equation~(\ref{eq:le2_2}),
	\begin{equation*}
	    \begin{alignedat}{3}
	     &E[\delta_i\ |\ \mathcal{I}_p^{(i)},p^{(i-1)}]\\
	     =&|\mathcal{I}_p^{(i)}|(\frac{1}{2}-\frac{(u+p^{(i-1)})/2-(u+p^{(i-1)})/2+b/|\mathcal{I}_p^{(i)}|}{u-p^{(i-1)}})\\
	     =&\frac{|\mathcal{I}_p^{(i)}|}{2}-\frac{b}{u-p^{(i-1)}}.
	    \end{alignedat}
	\end{equation*}
	Using the Law of Total Expectation again,
	\begin{equation}\label{eq:p2_2}
	    E[\delta_i\ |\ p^{(i-1)}]=E[E[\delta_i\ |\ \mathcal{I}_p^{(i)},p^{(i-1)}]]=\frac{E[|\mathcal{I}_p^{(i)}|]}{2}-\frac{b}{u-p^{(i-1)}}.
	\end{equation}
	Let $\sigma := \sqrt{\frac{2b}{u-l}}+1$. We will now consider cases.  Observe that $|\mathcal{I}_\tau| \leq n$, i.e. the active set is a subset of the ground set. The following two cases exhaust the possibilities of where $\sigma \sqrt{n}$ lies: either $|\mathcal{I}_\tau|<\sigma\cdot \sqrt{n}<n$ (Case 1); otherwise,  either $|\mathcal{I}_\tau|\leq n \leq \sigma\cdot \sqrt{n}$ or $\sigma\cdot \sqrt{n} \leq |\mathcal{I}_\tau|\leq n $ (Case 2). \\
 \emph{Case 1:  $|\mathcal{I}_\tau|<\sigma\cdot \sqrt{n}<n$}\\
 Let $\xi$ be an iteration such that $1\leq \xi<T$ satisfies:
	\[|\mathcal{I}_p^{(\xi)}|\geq \sigma\cdot\sqrt{n} > |\mathcal{I}_p^{(\xi+1)}|> ... > |\mathcal{I}_p^{(T)}| = |\mathcal{I}_\tau|;\]
    Then,
	\begin{equation}\label{eq:p2_8}
	    E[\sum_{i=1}^T|\mathcal{I}_p^{(i)}|] = E[\sum_{i=1}^{\xi}|\mathcal{I}_p^{(i)}|] + E[\sum_{i=\xi}^T|\mathcal{I}_p^{(i)}|]
	\end{equation}
	Since $|\mathcal{I}_\tau|<\sigma\cdot\sqrt{n}$, then from $|\mathcal{I}_p^{(\xi+1)}|<\sigma\cdot\sqrt{n}$ Michelot's Method will stop within at most $\sigma\cdot\sqrt{n}$ iterations since at least one element is removed in each iteration; thus $T-\xi<\sigma\cdot\sqrt{n}$. Since $|\mathcal{I}_p^{(i)}|< |\mathcal{I}_p^{(\xi+1)}|<\sigma\cdot\sqrt{n}$ for $i\geq \xi+1$,
	\begin{equation}\label{eq:p2_3}
	    E[\sum_{i=\xi}^T|\mathcal{I}_p^{(i)}|]\leq E[\sum_{i=\xi}^T \sigma\cdot\sqrt{n}]< \sigma^2\cdot n \in O(n).
	\end{equation}
	From Equation~(\ref{eq:p2_2}),
	\[E[\sum_{i=1}^{\xi}\delta_i\ |\ p^{(1)},...,p^{(\xi)}, \xi]=\frac{E[\sum_{i=1}^{\xi}|\mathcal{I}_p^{(i)}|\ |\ \xi]}{2}-\sum_{i=1}^{\xi}\frac{b}{u-p^{(i-1)}}.\]
	Then, from the Law of Total Expectation,
	\begin{equation}
 \label{eq:eqsix}
	    \begin{aligned}
	        E[\sum_{i=1}^{\xi} \delta_i] =& E[E[\sum_{i=1}^{\xi}\delta_i\ |\ p^{(1)},...,p^{(\xi)}, \xi]]\\
	        =& \frac{E[\sum_{i=1}^{\xi} |\mathcal{I}_p^{(i)}|]}{2}-E[\sum_{i=1}^{\xi}\frac{b}{u-p^{(i-1)}}].
	    \end{aligned}
	\end{equation}
	Michelot's method always maintains a nonempty active index set (decreasing per iteration), and so $\sum_{i=1}^{\xi}\delta_i = n - |\mathcal{I}_p^{(\xi)}| < n $. Then, continuing from~(\ref{eq:eqsix}),
	\[\frac{E[\sum_{i=1}^{\xi} |\mathcal{I}_p^{(i)}|]}{2}-E[\sum_{i=1}^{\xi}\frac{b}{u-p^{(i-1)}}]\leq n,\]
	\begin{equation}\label{eq:p2_4}
	    \begin{aligned}
	        \implies \frac{E[\sum_{i=1}^{\xi} |\mathcal{I}_p^{(i)}|]}{2}\leq& n + E[\sum_{i=1}^{\xi}\frac{b}{u-p^{(i-1)}}]\\
	    \end{aligned}
	\end{equation}
	\begin{claim}
    $E[\sum_{i=1}^{\xi}\frac{b}{u-p^{(i-1)}}]\in O(n)$.
    \end{claim}
	\begin{claimproof}
	    \begin{equation}\label{eq:p2_5}
	        \begin{aligned}
	        E[\sum_{i=1}^{\xi}\frac{b}{u-p^{(i-1)}}|\xi] = &\sum_{i=1}^{\xi} E[\frac{b}{u-p^{(i-1)}}]\\
	        \leq &\sum_{i=1}^{\xi}\frac{b}{E[u-p^{(i-1)]}}\quad (\mbox{Jensen's Inequality}).
	        \end{aligned}
	    \end{equation}
	    Since $\mathcal{I}_p^{(i)}=\{j\in\mathcal{I}\ |\ d_j>p^{(i-1)}\}$ and $d_1,...,d_n$ are $i.i.d.\sim U[l,u]$, 
	    \[E[|\mathcal{I}_p^{(i)}|]=n\cdot\frac{E[u-p^{(i-1)}]}{u-l}.\]
	    So for $i=1,...,\xi$, since $|\mathcal{I}_p^{(i)}|>|\mathcal{I}_p^{(i+1)}|$,
	    \[n\cdot\frac{E[u-p^{(i-1)}]}{u-l}= E[|\mathcal{I}_p^{(i)}|]> E[|\mathcal{I}_p^{(\xi)}|]\geq \sigma\cdot\sqrt{n};\]
	    \begin{equation}\label{eq:p2_7}
	        \implies E[u-p^{(i-1)}]>\sigma\cdot \frac{u-l}{\sqrt{n}}.
	    \end{equation}
	    Substituting into Inequality~(\ref{eq:p2_5}),
	    \[E[\sum_{i=1}^{\xi}\frac{b}{u-p^{(i-1)}}|\xi]\leq \sum_{i=1}^{\xi}\frac{\sqrt{n}b}{\sigma\cdot(u-l)}=\frac{\xi\sqrt{n}b}{\sigma\cdot(u-l)};\]
	    using Law of Total Expectation,
	    \begin{equation}\label{eq:p2_6}
	        E[\sum_{i=1}^{\xi}\frac{b}{u-p^{(i-1)}}]=E[E[\sum_{i=1}^{\xi}\frac{b}{u-p^{(i-1)}}|\xi]]\leq E[\xi]\frac{\sqrt{n}b}{\sigma\cdot(u-l)}.
	    \end{equation}
	    So then it remains to show $E[\xi]\in O(\sqrt{n})$. From Equation~(\ref{eq:le2_2}), for any iteration $i$
        \[u-E[p^{(i)}|\mathcal{I}_{p}^{(i)},p^{(i-1)}] = u-\frac{p^{(i-1)}+u}{2}+\frac{b}{|\mathcal{I}_{p}^{(i)}|},\]
    \[\implies E[u-p^{(i)}|\mathcal{I}_{p}^{(i)},p^{(i-1)}] = \frac{u-p^{(i-1)}}{2}+\frac{b}{|\mathcal{I}_{p}^{(i)}|};\]
    thus, using the Law of Total Expectation,
    \begin{equation}
    \label{eq:le2_3}
        E[u-p^{(i)}] = \frac{E[u-p^{(i-1)}]}{2}+E[\frac{b}{|\mathcal{I}_p^{(i)}|}].
    \end{equation}
    Applying Equation~(\ref{eq:le2_3}) recursively, starting with the base case $E[u-p^{(0)}]=u-l$,
    \begin{equation*}
        \begin{aligned}
            E[u-p^{(i)}] =& \frac{u-l}{2^i}+\sum_{j=1}^{i}E[\frac{b}{2^{i-j}\cdot|\mathcal{I}_p^{(j)}|}],\\
            \implies E[p^{(i)}] =& u - \frac{u-l}{2^i}-\sum_{j=1}^{i}E[\frac{b}{2^{i-j}\cdot|\mathcal{I}_p^{(j)}|}],\\
            \implies E[p^{(\xi-1)} | \xi] =& u - \frac{u-l}{2^{\xi-1}}-E[\sum_{j=1}^{\xi-1}\frac{b}{2^{\xi-1-j}\cdot|\mathcal{I}_p^{(j)}|}| \xi].
        \end{aligned}
    \end{equation*}
    Using the Law of Total Expectation,
    \[E[p^{(\xi-1)}] = E[E[p^{(\xi-1)}|\xi]] = u - E[\frac{u-l}{2^{\xi-1}}]-E[\sum_{j=1}^{\xi-1}\frac{b}{2^{\xi-1-j}\cdot|\mathcal{I}_p^{(j)}|}]\]
    \begin{equation*}
        \begin{aligned}
            \implies E[p^{(\xi-1)}] \geq&  u - E[\frac{u-l}{2^{\xi-1}}]-E[\sum_{j=1}^{\xi-1}\frac{b}{2^{\xi-1-j}\sigma\cdot\sqrt{n}}]\quad (\mbox{since }|\mathcal{I}_p^{(i)}|\geq \sigma\cdot\sqrt{n}\ \mbox{for } i\leq \xi)\\
            \geq& u - E[\frac{u-l}{2^{\xi-1}}]-\frac{2b}{\sigma\cdot\sqrt{n}}\\
            \geq & u - \frac{u-l}{E[2^{\xi-1}]}-\frac{2b}{\sigma\cdot\sqrt{n}} \quad (\mbox{Jensen's inequality})\\
            \implies \frac{u-l}{E[2^{\xi-1}]}\geq &E[u-p^{(\xi-1)}]-\frac{2b}{\sigma\cdot\sqrt{n}}
        \end{aligned}
    \end{equation*}
    From Inequality~(\ref{eq:p2_7}),
    $E[u-p^{(\xi-1)}]>\sigma\cdot\frac{u-l}{\sqrt{n}}$; thus
    \[\frac{u-l}{E[2^{\xi-1}]}>\sigma\cdot\frac{u-l}{\sqrt{n}}-\frac{2b}{\sigma\cdot\sqrt{n}}=\frac{1}{\sigma\cdot\sqrt{n}}[\sigma^2\cdot (u-l)-2b].\]
    Since $\sigma = \sqrt{\frac{2b}{u-l}}+1$, $\sigma > \sqrt{\frac{2b}{u-l}}$; thus $\sigma^2\cdot (u-l)-2b>0$; as a result,
    \[E[2^{\xi-1}]<\frac{\sigma\sqrt{n}(u-l)}{\sigma^2\cdot (u-l)-2b};\]
    \[\implies 2^{E[\xi]-1}\leq E[2^{\xi-1}]<\frac{\sigma\sqrt{n}(u-l)}{\sigma^2\cdot (u-l)-2b}\quad (\mbox{Jensen's Inequality});\]
    \[\implies E[\xi]\leq \log (\frac{\sigma\sqrt{n}(u-l)}{\sigma^2\cdot (u-l)-2b})+1.\]
    Thus $E[\xi]\in O(\sqrt{n})$. Substituting into Inequality~(\ref{eq:p2_6})
    \[E[\sum_{i=1}^{\xi}\frac{b}{u-p^{(i-1)}}]\leq (\log (\frac{\sigma\sqrt{n}(u-l)}{\sigma^2\cdot (u-l)-2b})+1)\cdot\frac{\sqrt{n}b}{\sigma\cdot(u-l)}\in O(n).\]
	\end{claimproof}
	
	The Claim, together with Inequality~(\ref{eq:p2_4}) establish that
	\[E[\sum_{i=1}^{\xi} |\mathcal{I}_p^{(i)}|]\in O(n).\]
	Altogether with Inequalities~(\ref{eq:p2_8}) and~(\ref{eq:p2_3}), we have $E[\sum_{i=1}^T|\mathcal{I}_p^{(i)}|] \in O(n)$ in \emph{Case 1}.\\
	\emph{Case 2: $n\geq |\mathcal{I}_\tau|\geq\sigma\cdot \sqrt{n}$, or $\sigma\cdot \sqrt{n}\geq n\geq |\mathcal{I}_\tau|$}\\
	If $n\geq |\mathcal{I}_\tau|\geq\sigma\cdot \sqrt{n}$, let $\xi := T$. Since \[|\mathcal{I}_p^{(\xi)}|=|\mathcal{I}_p^{(T)}|=|\mathcal{I}_\tau|\geq \sigma\cdot \sqrt{n},\]
	the Claim from \emph{Case 1} and Equation~(\ref{eq:p2_4}) hold for $\xi = T$; thus
	\[E[\sum_{i=1}^{T}\frac{b}{u-p^{(i-1)}}]\in O(n),\]
	\[\frac{E[\sum_{i=1}^{T} |\mathcal{I}_p^{(i)}|]}{2}\leq n + E[\sum_{i=1}^{T}\frac{b}{u-p^{(i-1)}}];\]
    \[\implies E[\sum_{i=1}^{T} |\mathcal{I}_p^{(i)}|] \in O(n).\]
    
    If $\sigma\cdot \sqrt{n}\geq n\geq |\mathcal{I}_\tau|$, then $\sigma \geq \sqrt{n}$, which implies $n\leq \sigma^2$. However, $u,l,b$ are given independently of $n$ and so $\sigma$ is fixed.  Thus for asymptotic analysis, $n\leq \sigma^2$ does not hold for sufficiently large $n$.
    
    Together with \emph{Case 1} and \emph{Case 2}, $E[\sum_{i=1}^{T} |\mathcal{I}_p^{(i)}|] \in O(n)$. Hence $O(n)$ operations are used for scanning/prefix-sum.  All other operations, i.e. assigning $\mathcal{I}$ and $\mathcal{I}_p$, are within a constant factor of the scanning operations.
\end{proof}
~\\
We assume uniformly distributed inputs, $d_{1},\dots,d_{n}$ are $\mathrm{i.i.d} \sim U[l,u]$, and we have\\
\\
\textbf{Lemma 2}
    Filter provides a pivot $p$ such that $\tau \geq p\geq \frac{\sum_{i\in\mathcal{I}}d_i-b}{|\mathcal{I}|}$.
\begin{proof}
    The upper bound $p \leq \tau$ is given by construction of $p$ (see \cite[Section 3, Paragraph 2]{condat2016}).
    
    We can establish the lower bound on $p$ by considering the sequence $p^{(1)}\leq...\leq p^{(n)}\in\mathbb{R}$, which represents the initial as well as subsequent (intermediate) values of $p$ from the first outer for-loop on line 2 of presented as Algorithm 4 (from the main body), and their corresponding index sets $\mathcal{I}_p^{(1)}\subseteq...\subseteq\mathcal{I}_p^{(n)}$.
    Filter initializes with $p^{(1)}:=d_1-b$ and $\mathcal{I}_p^{(1)}:=\{1\}$. For $i=2,...,n$, if $d_i>p^{(i-1)}$, \begin{alignat*}{1}
    p^{(i)}:=p^{(i-1)}+(d_i-p^{(i-1)})/(|\mathcal{I}_p^{(i-1)}|+1),\ \mathcal{I}_p^{(i)}:=\mathcal{I}_p^{(i-1)}\cup\{i\};
    \end{alignat*}
    otherwise $p^{(i)}:=p^{(i-1)}$, and $\mathcal{I}_p^{(i)}:=\mathcal{I}_p^{(i-1)}$. Then it can be shown that  $p^{(i)}=(\sum_{j\in\mathcal{I}_p^{(i)}}d_j-b)/|\mathcal{I}_p^{(i)}|$ (see \cite[Section 3, Paragraph 2]{condat2016}), and $p\geq p^{(n)}$ (see \cite[Section 3, Paragraph 5]{condat2016}). Now in terms of $p^{(n)}$ we may write
    \begin{alignat*}{3}
        \frac{\sum_{i\in\mathcal{I}}d_i-b}{|\mathcal{I}|}=&\frac{\sum_{i\in\mathcal{I}_p^{(n)}}d_i-b+\sum_{i\in\mathcal{I}\backslash\mathcal{I}_p^{(n)}}d_i}{|\mathcal{I}|}\\
        =&\frac{(\sum_{i\in\mathcal{I}_p^{(n)}}d_i-b)|\mathcal{I}_p^{(n)}|}{|\mathcal{I}||\mathcal{I}_p^{(n)}|}+\frac{\sum_{i\in\mathcal{I}\backslash\mathcal{I}_p^{(n)}}d_i}{|\mathcal{I}|}\\
        =&\frac{|\mathcal{I}_p^{(n)}|}{|\mathcal{I}|}p^{(n)}+\frac{\sum_{i\in\mathcal{I}\backslash\mathcal{I}_p^{(n)}}d_i}{|\mathcal{I}|}\\
        =& p^{(n)}+\frac{\sum_{i\in\mathcal{I}\backslash\mathcal{I}_p^{(n)}}d_i-(|\mathcal{I}|-|\mathcal{I}_p^{(n)}|)p^{(n)}}{|\mathcal{I}|}\\
        =&p^{(n)}+\frac{\sum_{i\in\mathcal{I}\backslash\mathcal{I}_p^{(n)}}(d_i-p^{(n)})}{|\mathcal{I}|}.
    \end{alignat*}
    For any $i\in \mathcal{I}\backslash\mathcal{I}_p^{(n)}$, since $\mathcal{I}_p^{(i)}\subseteq \mathcal{I}_p^{(n)}$ then $i\not\in\mathcal{I}_p^{(i)}$. By construction of $p^{(i-1)}$ and $\mathcal{I}_p^{(i)}$, we have $d_i\leq p^{(i-1)}\leq p^{(n)}$. Thus, $\sum_{i\in\mathcal{I}\backslash\mathcal{I}_p^{(n)}}(d_i- p^{(n)})\leq 0$, and $p^{(n)}\geq (\sum_{i\in\mathcal{I}}d_i-b)/|\mathcal{I}|$. So $p\geq p^{(n)}\geq (\sum_{i\in\mathcal{I}}d_i-b)/|\mathcal{I}|$.
\end{proof}
~\\
Now we introduce some notation in order to compare subsequent iterations of Condat's method with iterations of Michelot's method. Let $t_C\leq n$ and $t_M \leq n$ be the total number of iterations taken by Condat's method and Michelot's method (respectively) on a given instance. Let $\mathcal{I}^C_0,...,\mathcal{I}^C_n$ be the active index sets per iteration for Condat's method with corresponding pivots $p^C_0,...,p^C_n$. Likewise, we denote the index sets and pivots of Michelot's method as $\mathcal{I}^M_0,...,\mathcal{I}^M_n$ and $p^M_0,...,p^M_n$, respectively. If $t_C< n$ then we set $p^C_{t_C}=p^C_{t_C+1}=...=p^C_n=\tau$ and $\mathcal{I}^C_{t_C}=\mathcal{I}^C_{t_C+1}=...=\mathcal{I}^C_n=\mathcal{I}_\tau$; likewise for Michelot's algorithm.
~\\
~\\
\textbf{Lemma 3}
$\mathcal{I}_i^C\subseteq \mathcal{I}_i^M$, and $p^C_i\geq p^M_i$ for $i=0,...,n$.
\begin{proof}
We will prove this by induction. For the base case, $\mathcal{I}_0^{C}$ is obtained by Filter. So $\mathcal{I}_0^{C}\subseteq \mathcal{I}=\mathcal{I}_0^M$. Moreover, from Lemma 2, $p_0^C\geq p_0^M$.

Now for any iteration $i\geq 1$, suppose $\mathcal{I}^C_i\subseteq \mathcal{I}^M_i$, and $p_i^C\geq p_i^M$. From line 5 in Algorithm 3 (from the main body), $\mathcal{I}^M_{i+1}:=\{j\in\mathcal{I}^M_i:d_j>p_i^M\}$. From Condat \cite[Section 3, Paragraph 3]{condat2016}, Condat's method uses a dynamic pivot between $p^C_i$ to $p^C_{i+1}$ to remove inactive entries that would otherwise remain in Michelot's method. Therefore, $\mathcal{I}^C_{i+1}\subseteq \{j\in\mathcal{I}^C_i:d_j>p_i^C\}\subseteq \mathcal{I}^M_{i+1}$, and moreover for any $j\in \mathcal{I}^M_{i+1}\backslash\mathcal{I}^C_{i+1}$, we have that $d_j\leq p^C_{i+1}$. Now observe that
\begin{alignat*}{3}
        &p^{C}_{i+1}-p^{M}_{i+1}\\
        =&\frac{\sum_{j\in\mathcal{I}^{C}_{i+1}}d_j-b}{|\mathcal{I}^{C}_{i+1}|}-\frac{\sum_{j\in\mathcal{I}^{M}_{i+1}}d_j-b}{|\mathcal{I}^{M}_{i+1}|}\\
        =&\frac{\sum_{j\in\mathcal{I}^{C}_{i+1}}d_j-b}{|\mathcal{I}^{C}_{i+1}|}-\frac{\sum_{j\in\mathcal{I}^{C}_{i+1}}d_j+\sum_{j\in\mathcal{I}^{M}_{i+1}\backslash\mathcal{I}^{C}_{i+1}}d_j-b}{|\mathcal{I}^{C}_{i+1}|+|\mathcal{I}^{M}_{i+1}\backslash\mathcal{I}^{C}_{i+1}|}
        \\
        =&\frac{\sum_{j\in\mathcal{I}_{i+1}^{M}\backslash\mathcal{I}_{i+1}^{C}}(p_{i+1}^{C}-d_j)}{|\mathcal{I}_{i+1}^{C}|+|\mathcal{I}_{i+1}^{M}\backslash\mathcal{I}_{i+1}^{C}|}\geq 0,
    \end{alignat*}
    and so $p^{C}_{i+1}\geq p^{M}_{i+1}$.
\end{proof}
~\\
\textbf{Corollary 2}
$\sum_{i=1}^{t_C}|\mathcal{I}_i^C|\leq \sum_{i=1}^{t_M}|\mathcal{I}_i^M|$.
\begin{proof}
Observe that both algorithms remove elements (without replacement) from their candidate active sets $\mathcal{I}^C_i,\mathcal{I}^M_i$ at every iteration; moreover, they terminate with the pivot value $\tau$ and so   $\mathcal{I}^C_{t_C}=\mathcal{I}^M_{t_M}=\mathcal{I}_\tau$.  So, together with Lemma 3, we have for  $i=0,...,n$ that $\mathcal{I}_\tau\subseteq \mathcal{I}_i^C\subseteq \mathcal{I}_i^M$. So $\mathcal{I}_{t_M}^M=\mathcal{I}_\tau$ implies $\mathcal{I}_{t_M}^C=\mathcal{I}_\tau$, and so $t_C\leq t_M$. Therefore $\sum_{i=1}^{t_C}|\mathcal{I}_i^C|\leq \sum_{i=1}^{t_M}|\mathcal{I}_i^M|$.
\end{proof}
~\\
\textbf{Lemma 4}
The worst-case runtime of Filter is $O(n)$.
\begin{proof}
Since $\mathcal{I}_w\subseteq \mathcal{I}$ at any iteration, Filter will scan at most $2|\mathcal{I}|$ entries; including $O(1)$ operations to  update $p$. 
\end{proof}

We assume uniformly distributed inputs, $d_{1},\dots,d_{n}$ are $\mathrm{i.i.d} \sim U[l,u]$, and we have\\
~\\
\textbf{Proposition 3}
	Condat's method has an average runtime of $O(n)$.
\begin{proof}
	Filter takes $O(n)$ operations from Lemma 4. From Corollary 2, the total operations spent on scanning in Condat's method is less than (or equal to) the average $O(n)$ operations for Michelot's method (established in Proposition 2); hence Condat's average runtime is $O(n)$.
\end{proof}

We assume uniformly distributed inputs, $d_{1},\dots,d_{n}$ are $\mathrm{i.i.d} \sim U[l,u]$, and we have\\
~\\
\textbf{Theorem 1}
    $E[|\mathcal{I}_\tau|]  < \sqrt{\frac{2b(n+1)}{u-l}+\frac{1}{4}}+\frac{1}{2}$.
\begin{proof}
    Sort $d$ such that $d_{\pi_1}\geq d_{\pi_2}\geq...\geq d_{\pi_n}$. Thus for a given order statistic, (see e.g.  \cite[p. 63]{orderstatistic}),
    \[E[d_{\pi_i}]=u-\frac{i}{n+1}(u-l).\]
    Define $N:=|\mathcal{I}_\tau|$ for ease of presentation. From Corollary 1,
    \[\frac{\sum_{i=1}^{N}d_{\pi_i}-b}{N}=\tau;\]
    and, together with $d_{\pi_{N+1}}\leq \tau < d_{\pi_N}$ (by definition of $\mathcal{I}_\tau$), we have 
    \[N\cdot d_{\pi_{N+1}}\leq \sum_{i=1}^{N}d_{\pi_i}-b< N\cdot d_{\pi_{N}},\]
    \begin{equation}\label{neq:EN}
        \implies E[\sum_{i=1}^{N}d_{\pi_i}-b]< E[N\cdot d_{\pi_{N}}].
    \end{equation}
    Furthermore, from the Law of Total Expectation,
    \[E[N\cdot d_{\pi_{N}}]=uE[N]-E[N^2]\frac{u-l}{n+1},\]
    \[E[\sum_{i=1}^{N}d_{\pi_i}-b]=uE[N]-E[N(N+1)]\frac{u-l}{2(n+1)}-b.\]
    Substituting into~(\ref{neq:EN}) yields
    \begin{equation*}
        E[N^2]-E[N]<\frac{2b(n+1)}{u-l},
    \end{equation*}
    Since $E^2[N]\leq E[N^2]$, then
    \[E^2[N]-E[N]\leq E[N^2]-E[N]< \frac{2b(n+1)}{u-l},\]
    \begin{equation}\label{eq: E[I]}
        \implies E[|\mathcal{I}_\tau|] = E[N] < \sqrt{\frac{2b(n+1)}{u-l}+\frac{1}{4}}+\frac{1}{2}\in O(\sqrt{n}).
    \end{equation}
\end{proof}
~\\
\textbf{Lemma 5}
Suppose $d_1,...,d_n$ are $\mathrm{i.i.d.}$ from an arbitrary distribution $X$, with PDF $f_X$ and CDF  $F_X$. Let $\epsilon$ such that $0<\epsilon < 1$ be some positive number, and $t\in\mathbb{R}$ be such that $1-F_X(t)=\epsilon$.  Then i) $|\mathcal{I}_t|\rightarrow \infty$ as $n\rightarrow \infty$ and ii) $P(\frac{|\mathcal{I}_t|}{n}\leq \epsilon)=1\ \mbox{as}\ n\rightarrow \infty$.
\begin{proof}
Since $f_X$ is a density function, the CDF $F_X$ is absolutely continuous (see e.g. \citep[Page 59, Definition 2.1]{cdf}); thus for any $0<\epsilon < 1$, there exists $t\in \mathbb{R}$ such that $F_x(t) = 1-\epsilon$.

For $i=1,...,n$, define the indicator variable  	
\[\delta_{i}:=\begin{cases} 1, & \mbox{if}\ d_{i}> t\\ 0, &\mathrm{otherwise} \end{cases}\] 
and let $S_{n}:=\sum_{i=1}^{n}\delta_{i}$. So $S_{n}=|\mathcal{I}_t|$, and we can show that it is binomially distributed. 
	
\begin{claim} $S_n\sim B(n,\epsilon)$.
\end{claim}
\begin{claimproof}
Observe that $P(\delta_i=1)=P(d_i>t)=\epsilon$, and $P(\delta_i=0)=P(d_i\leq t)=1-\epsilon$; thus $\delta_i\sim  \mbox{Bernoulli}(\epsilon)$. Moreover, $d_1,...,d_n$ are independent, and so $\delta_1,...,\delta_n$ are i.i.d. $\mbox{Bernoulli}(\epsilon)$; consequently, $S_n:=\sum_{i=1}^{n}\delta_i \sim B(n,\epsilon)$. 
\end{claimproof}

    So, as $n\rightarrow \infty$, we can apply the Central Limit Theorem (see e.g. \cite{centrallimit}):
	\begin{equation}
	\label{eq:norm distribution}
	    S_{n}^{*}:=\frac{S_{n}-n\epsilon}{\sqrt{n\epsilon(1-\epsilon)}}\sim N(0,1).
	\end{equation} 
    Denote $\Phi(q)=\int_{-\infty}^{q}\frac{1}{\sqrt{2\pi}}e^{-\frac{x^2}{2}}\,dx$ to be the CDF of the standard normal distribution.  Consider (for any $q \in \mathbb{R}$) the right-tail probability
    
	\begin{equation} \label{eq: z score}
			P(\frac{S_{n}-n\epsilon}{\sqrt{n\epsilon(1-\epsilon)}}\geq -q)=1-\Phi(-q)=\Phi(q),
	\end{equation}
	\[P(S_{n} \geq n\epsilon-q\sqrt{n\epsilon (1-\epsilon)})=\Phi(q),\]
	\begin{equation}\label{eq: range of It}
	    \implies P(|\mathcal{I}_{t}|\geq \lfloor-q\sqrt{n\epsilon(1-\epsilon)}+n\epsilon\rfloor) \geq \Phi(q),
	\end{equation}
	where $\lfloor . \rfloor$ is the floor function.  

    Setting $q=\sqrt{2\log n}$ yields 
	\[P(|\mathcal{I}_t|\geq \lfloor-\sqrt{2(n\log n)\epsilon(1-\epsilon)}+n\epsilon\rfloor)\geq \Phi(\sqrt{2\log n}).\]
	Consider the right-hand side, $\lim_{n\rightarrow \infty} \Phi(\sqrt{2\log n})$. Since $\Phi$ is a CDF, it is monotonically increasing, continuous, and converges to $1$; thus  $\lim_{n\rightarrow\infty}\Phi(\sqrt{2\log n})=1$. So as $n$ approaches infinity,
    \begin{equation}\label{eq:It_order}
        P(|\mathcal{I}_t|\geq \lfloor-\sqrt{2(n\log n)\epsilon(1-\epsilon)}+n\epsilon\rfloor)=1,
    \end{equation}
    which establishes condition (i).

    Now consider the left-tail probability,
    \[P(\frac{S_{n}-n\epsilon}{\sqrt{n\epsilon(1-\epsilon)}}\leq q)=\Phi(q)\]
    \[\implies P(|\mathcal{I}_{t}|\leq \lceil q\sqrt{n\epsilon(1-\epsilon)}+n\epsilon\rceil) \geq \Phi(q).\]
    Again setting $q=\sqrt{2\log n}$, we have that $\lim_{n\rightarrow \infty}\lceil q\sqrt{n\epsilon(1-\epsilon)}+n\epsilon\rceil/n=\epsilon$. So as $n$ approaches infinity,
    \[P(\frac{|\mathcal{I}_t|}{n}\leq \epsilon)=1,\]
    which establishes condition (ii).
\end{proof}
~\\
\textbf{Theorem 2}
	Suppose $d_1,...,d_n$ are $\mathrm{i.i.d.}$ from an arbitrary distribution $X$, with PDF $f_X$ and CDF  $F_X$. Then, for any $\epsilon>0$, $P(\frac{|\mathcal{I}_\tau|}{n}\leq \epsilon)=1$ as $n\rightarrow \infty$.
\begin{proof}
    \emph{Case 1}: $\epsilon \geq 1$\\
    Since $\mathcal{I}_\tau \subseteq \mathcal{I}$ and $|\mathcal{I}| = n$, $\frac{|\mathcal{I}_\tau|}{n}\leq 1$. Hence, $P(\frac{|\mathcal{I}_\tau|}{n}\leq \epsilon) = 1$.\\ 
    \emph{Case 2}: $0<\epsilon<1$\\
    Since $f_X$ is a density function, $F_X$ is absolutely continuous. So there exists $t\in\mathbb{R}$ be such that $1-F_X(t)=\epsilon$. We shall first establish that $P(\tau> t)=1$ as $n\rightarrow \infty$.
    
    From Corollary 1, $\tau> t \iff f(t)> 0$, and so

	\begin{alignat}{3}
	    	P(\tau> t)=& P(f(t)> 0)\nonumber\\
			=&P(\frac{\sum_{i\in \mathcal{I}_{t}}d_i-b}{|\mathcal{I}_t|}>t)\nonumber\\
			=& 1-P(\sum_{i\in \mathcal{I}_{t}}d_i\leq |\mathcal{I}_t|\cdot t+b) \nonumber\\
			\label{eq:f(t)>0}
			=& 1-P(\sum_{i\in \mathcal{I}_{t}}d_i-E[\sum_{i\in \mathcal{I}_{t}}d_i]\leq |\mathcal{I}_t|\cdot t+b-E[\sum_{i\in \mathcal{I}_{t}}d_i]).
	\end{alignat}
	
  Now observe that, for any  $i\in\mathcal{I}_t$, $d_i$ can be treated as a conditional variable: $d_i|_{d_i>t}$. Since all such $d_i$ are i.i.d, we may denote the (shared) expected value $\mu:=E[d_i|d_i>t]$ and variance as $\sigma^2:=\mbox{Var}[d_i|d_i>t]$. Moreover, by definition of $\mathcal{I}_t$ we have
  \begin{equation}
      \label{eq: mu_t}
      E[d_i| i\in\mathcal{I}_t]=\mu>t.
  \end{equation}

   Together with condition (i) of Lemma~5, this implies that the right-hand side of the probability in (\ref{eq:f(t)>0}) is negative as $n \rightarrow \infty$:
   \begin{equation}\label{eq:corollarycore}
       |\mathcal{I}_t|\cdot t+b-E[\sum_{i\in\mathcal{I}_t}d_i]=|\mathcal{I}_t|\cdot (t-\mu)+b<0.
   \end{equation}
	It follows, continuing from Equation~(\ref{eq:f(t)>0}), that 
	\begin{equation*}
		\begin{aligned}
		P(f(t)>0)&\geq 1-P(|\sum_{i\in\mathcal{I}_t}d_i-E[\sum_{i\in\mathcal{I}_t}d_i]|\geq E[\sum_{i\in\mathcal{I}_t}d_i]-|\mathcal{I}_t|\cdot t-b)\\
			&\geq 1-\frac{\texttt{Var}(\sum_{i\in\mathcal{I}_t}d_i)}{(E[\sum_{i\in\mathcal{I}_t}d_i]-|\mathcal{I}_t|\cdot t-b)^{2}} \ \  \mbox{(Chebyshev's inequality)}\\
			&= 1-\frac{\sigma^2|\mathcal{I}_t|}{(\mu|\mathcal{I}_t|-t|\mathcal{I}_t|-b)^{2}}.
		\end{aligned}
	\end{equation*}
	Thus we have the desired result:
	\begin{equation}
	\label{eq:tau>t}
	    P(\tau>t)\geq 1-\frac{\sigma^2|\mathcal{I}_t|}{(\mu|\mathcal{I}_t|-t|\mathcal{I}_t|-b)^{2}} \rightarrow 1\ \mbox{as}\ n \rightarrow \infty.
	\end{equation}

    Now observe that $\tau>t$ implies $\mathcal{I}_\tau\subseteq \mathcal{I}_t$, and subsequently $|\mathcal{I}_\tau|\leq \mathcal{I}_t$. Together with  condition (ii) from Lemma 5  we have
    \[P(\frac{|\mathcal{I}_\tau|}{n}\leq \epsilon)\geq P(\frac{|\mathcal{I}_t|}{n}\leq \epsilon)=1,\ \mbox{as}\ n\rightarrow \infty.\]
\end{proof}
~\\
\textbf{Corollary 2}
For projection of $d\in\mathbb{R}^n$ onto a simplex $\Delta_b$, if $b\in o(n)$, the conclusion from Theorem 2 keeps true.
\begin{proof}
From Equation~(\ref{eq:It_order}), $|\mathcal{I}_t| \in \Theta(n)$. If $b\in o(n)$, since Equation~(\ref{eq: mu_t}) implies $t-\mu < 0$, Equation~(\ref{eq:corollarycore}) keeps true. As a result, Theorem 2 keeps true.
\end{proof}
~\\
\textbf{Proposition 4}
	Let $\hat{d}$ be a subvector of $d$ with $m\leq n$ entries; moreover, without loss of generality suppose the subvector contains the first $m$ entries. Let $\hat v^*$ be the projection of $\hat{d}$ onto the simplex $\hat \Delta := \{v\in\mathbb{R}^m \ |\  \sum_{i=1}^m v_i = b, v\geq 0 \}$, and $\hat{\tau}$ be the corresponding pivot value. Then, $\tau\geq \hat{\tau}$. Consequently, for $1\leq i\leq m$ we have that $\hat v_i^*=0 \implies v_i^* =0$.
\begin{proof}
	Define two index sets,
	\[\mathcal{I}_{\hat{\tau}}:=\{i=1,...,n\ |\ d_{i}>\hat{\tau}, d_i \in d\};\]
	\[\hat{\mathcal{I}}_{\hat{\tau}}:=\{i=1,...,m\ |\ d_{i}>\hat{\tau}, d_i \in \hat{d}\}.\]
	As $\hat{d}$ is a subvector of $d$, we have $\hat{\mathcal{I}}_{\hat{\tau}}\subseteq \mathcal{I}_{\hat{\tau}}$; thus, 
	\[\sum_{i\in \mathcal{I}_{\hat{\tau}}}(d_{i}-\hat{\tau})\geq \sum_{i\in \hat{\mathcal{I}}_{\hat{\tau}}}(d_{i}-\hat{\tau})=b,\]
	\[\implies \frac{\sum_{i\in \mathcal{I}_{\hat{\tau}}}d_{i}-b}{|\mathcal{I}_{\hat{\tau}}|} \geq \hat{\tau}.\]
	From Corollary 1, $\tau\geq \hat{\tau}$; from Proposition 1 it thus follows that $\mathcal{I}_{\hat \tau} \supset \mathcal{I}_\tau$.
\end{proof}

We assume uniformly distributed inputs, $d_{1},\dots,d_{n}$ are $\mathrm{i.i.d} \sim U[l,u]$, and we have\\
~\\
\textbf{Proposition 5}
	Parallel Pivot and Partition with either the median, random, or Michelot's pivot rule, has an average runtime of $O(\frac{n}{k}+\sqrt{kn})$. 
\begin{proof}
   The algorithm starts by distributing projections. Pivot and Partition has linear runtime on average with any of the stated pivot rules, and so for the $i$th core with input $d^i$, whose size is $O(\frac{n}{k})$, we have an average runtime of $O(\frac{n}{k})$. Now from Theorem 1, each core returns in expectation at most $O(\sqrt{\frac{n}{k}})$ active terms. So the reduced input $\hat v$ (line 3 of Algorithm 8 from the main body) will have at most $O(\sqrt{kn})$ entries on average; thus the final projection will incur an expected number of operations in $O(\sqrt{kn})$.
\end{proof}

We assume uniformly distributed inputs, $d_{1},\dots,d_{n}$ are $\mathrm{i.i.d} \sim U[l,u]$, and we have\\
~\\
\textbf{Proposition 6}
	Let $\mathcal{I}_p$ be the output of $\mathrm{Filter}(d,b)$. Then $E[|\mathcal{I}_p|] \in O(n^{\frac{2}{3}})$.
\begin{proof}
   We assume that the \textbf{if} in line 5 from Algorithm 4 (from the main body) does not trigger; this is a conservative assumption as otherwise more elements would be removed from $\mathcal{I}_p$, reducing the number of iterations.
    
    Let $p^{(j)}$ be the $j$th pivot with $p^{(1)}:=d_1-b$ (from line 1 in Algorithm 4 from the main body), and subsequent pivots corresponding to the for-loop iterations of line 2. Whenever Filter finds  $d_i>p^{(j)}$, it updates $p$ as follows:
    \begin{equation}\label{eq:updatepivot}
        p^{(j+1)}:=p^{(j)}+\frac{d_i-p^{(j)}}{j+1}.
    \end{equation}
    Moreover, when $d_i>p^{(j)}$ is found, $d_i\sim U[p^{(j)},u]$; thus, from the Law of Total Expectation, $E[d_i]=E[E[d_i|p^{(j)}]]=(u+E[p^{(j)}])/2$. Together with~(\ref{eq:updatepivot}),
    \begin{alignat*}{3}
        E[p^{(j+1)}]=&E[p^{(j)}+\frac{d_i-p^{(j)}}{j+1}]\\
        =&E[p^{(j)}]+\frac{E[d_i]-E[p^{(j)}]}{j+1}\\
        =&E[p^{(j)}]+\frac{(u+E[p^{(j)}])/2-E[p^{(j)}]}{j+1}\\
        =&\frac{(2j+1)E[p^{(j)}]+u}{2j+2};\\
    \implies E[p^{(j+1)}]-u=&(E[p^{(j)}]-u)\frac{2j+1}{2j+2}.
    \end{alignat*}
    
    Using the initial value $E[p^{(1)}]=E[E[p^{(1)}|d_1]]=E[d_1]-b=\frac{u+l}{2}-b$, we can obtain a closed-form representation for the recursive formula:
    \begin{equation}
	E[p^{(j)}]-u=-2(b+\frac{u-l}{2})\prod_{i=0}^{j-1}\frac{2i+1}{2i+2}.
	\label{eq:expectedpivot}
    \end{equation}
    
    Now let $L_j$ denote the number of terms Filter scans after calculating $p^{(j)}$ and  before finding some $d_i>p^{(j)}$. Since Filter scans $n$ terms in total (from initialize and $n-1$ calls to line 4), then
    \[\sum_{j=1}^{|\mathcal{I}_p|-1}L_j\leq n-1 < \sum_{j=1}^{|\mathcal{I}_p|}L_j;\]
    \begin{equation}\label{neq:ELj}
        \implies 1+\sum_{j=1}^{|\mathcal{I}_p|-1}E[L_j]\leq n<1+\sum_{j=1}^{|\mathcal{I}_p|}E[L_j].
    \end{equation}
    We can show $L_j$ has a geometric distribution as follows.
    \begin{claim} $L_j\sim \mbox{Geo}(\frac{u-p^{(j)}}{u-l})$.
	\end{claim}
\begin{claimproof}
For each term $d_i\sim U[l,u]$, we have
\[P(d_i>p^{(j)})=1-(p^{(j)}-l)/(u-l)=(u-p^{(j)})/(u-l).\] 
So $d_i>p^{(j)}$ can be interpreted as a Bernoulli trial. Hence $L_j$ is distributed with $\mbox{Geo}(\frac{u-p^{(j)}}{u-l})$. 
\end{claimproof}

    Now, applying Jensen’s Inequality to $E[L_j]$, 
    \[E[L_j]=E[\frac{u-l}{u-p^{(j)}}]\leq \frac{u-l}{u-E[p^{(j)}]},\]
    \[\implies \ln(E[L_j])\leq E[\ln(L_j)]=\ln(u-l)-\ln(u-E[p^{(j)}]),\]
     and together with (\ref{eq:expectedpivot}) we have
	\begin{equation}\label{eq:ELj2}
	    \begin{aligned}
	        \ln(E[L_j])\leq & \ln(\frac{u-l}{2b+u-l})+\sum_{i=0}^{j-1}\ln(\frac{2i+2}{2i+1})\\
	        =& \ln(\frac{u-l}{2b+u-l})+\sum_{i=1}^{j}\ln(\frac{2i}{2i-1})
	    \end{aligned}
	\end{equation}
	Now observe that
	\[\lim_{i\rightarrow \infty}\frac{\ln(\frac{2i}{2i-1})}{\frac{1}{i}}=\lim_{i\rightarrow \infty}\frac{\frac{2i-1}{2i}\frac{4i-2-4i}{(2i-1)^{2}}}{-\frac{1}{i^{2}}}=\lim_{i\rightarrow \infty}\frac{2i^{2}}{(2i-1)(2i)}=\frac{1}{2},\]
	where the first equality is from L'H\^{o}pital's rule.  Thus $E[L_j]\in \Theta(\mbox{exp}(\sum_{i=1}^{j}\frac{1}{2i}))$.
	Furthermore, we have the classical bound on the harmonic series:
	\[\sum_{i=1}^{j}\frac{1}{i}=\ln(j)+\gamma+\frac{1}{2j}\leq \ln(j)+1,\]
	where $\gamma$ is the Euler-Mascheroni constant, see e.g. \cite{EulerConstant}; thus,
	\[e^{\sum_{i=1}^{j}\frac{1}{2i}}=\sqrt{j}+e^{\frac{\gamma}{2}+\frac{1}{4j}}, \]
	which implies $E[L_{j}]\in O(\sqrt{j})$. Moreover (see e.g. \cite[Section 1.2.7]{harmonic}),
	\[\sum_{j=1}^{|\mathcal{I}_p|}\sqrt{j}= \frac{2}{3}|\mathcal{I}_p|\sqrt{|\mathcal{I}_p|+\frac{3}{2}}+o(\sqrt{|\mathcal{I}_p|}),\]
	which implies $1+\sum_{j=1}^{|\mathcal{I}_p|-1}E[L_j] \in O(|\mathcal{I}_p|^{\frac{3}{2}})$. From (\ref{neq:ELj}) it follows that $|\mathcal{I}_p|\in O(n^{\frac{2}{3}})$. 
\end{proof}

We assume uniformly distributed inputs, $d_{1},\dots,d_{n}$ are $\mathrm{i.i.d} \sim U[l,u]$, and we have\\
~\\
\textbf{Proposition 7}
	Parallel Condat's method has an average complexity $O(\frac{n}{k}+\sqrt[3]{kn^{2}})$.
\begin{proof}
    In Distributed Filter (Algorithm 9 from the main body), each core is given input $\mathcal{I}_i$ with $|\mathcal{I}_i|\in O(\frac{n}{k})$. (serial) Filter has linear runtime from  Lemma 4, and the \textbf{for} loop in line 5 of Algorithm 9 (from the main body) will scan at most $|\mathcal{I}_i|$ terms. Thus distributed Filter runtime is in $O(\frac{n}{k})$.

    From Proposition 6, $E[|\mathcal{I}_i|]\in O((n/k)^{\frac{2}{3}})$. Since the output of Distributed Filter is $\mathcal{I}_p=\cup_{i=1}^{k}\mathcal{I}_i$, we have that (given $d$ is i.i.d.) $E[|\mathcal{I}_p|]\in O(\sqrt[3]{kn^2})$.

    
    Parallel Condat's method takes the input from Distributed Filter and applies serial Condat's method (Algorithm 10 from the main body, lines 2-7), excluding the serial Filter.  From Proposition 3, Condat's method has average linear runtime, so this application of serial Condat's method has average complexity $O(\sqrt[3]{kn^2})$. 
\end{proof}
~\\
\textbf{Lemma 6}
If $X\sim U[l,u]$ and $l<t<u$, then $X|t \sim U[t,u]$.
\begin{proof}
 The CDF of $X$ is $F_{X}(x)=P(X\leq x)=(x-l)/(u-l)$. Then, 
    \[F_{X|t}(x)=P(X\leq x|X>t)=\frac{P(X\leq x,X>t)}{P(X>t)}=(\frac{x-t}{u-l})/(\frac{u-t}{u-l})=\frac{x-t}{u-t},\]
    which implies $X|t\sim U[t,u]$.
\end{proof}
~\\
\textbf{Lemma 7}
If $X,Y$ are independent random variables, then $X|t,Y|t$ are independent.
\begin{proof}
$X,Y$ are independent and so $P(X\leq x,Y\leq y)=P(X\leq x)P(Y\leq y)$ for any $x,y \in \mathbb{R}$. So considering the joint conditional probability,
    \begin{alignat*}{3}
        &P(X\leq x,Y\leq y|X>t,Y>t)\\
        &=\frac{P(X\leq x,Y\leq y,X>t,Y>t)}{P(X>t,Y>t)}\\
        &=\frac{P(X\leq x,X>t)}{P(X>t)}\frac{P(Y\leq y,Y>t)}{P(Y>t)}\\
        &=P(X\leq x|X>t)P(Y\leq y|y>t).
    \end{alignat*}
\end{proof}    
~\\
\textbf{Proposition 9}
If $d_1,...,d_n$ i.i.d. $\sim U[l,u]$ and $l<t<u$, then $\{d_i\ |\ i\in\mathcal{I}_t\}$ i.i.d. $\sim U[t,u]$.
\begin{proof}
From Lemma 6, for any $i\in\mathcal{I}_t$, $d_i \sim U[t,u]$. From Lemma 7, for any $i,j\in\mathcal{I}_t$, $i\neq j$, $d_i,d_j$ are conditionally independent. So, $\{d_i\ |\ i\in\mathcal{I}_t\}$ i.i.d. $\sim U[t,u]$.
\end{proof}

\newpage

\section{Algorithm Descriptions}
\subsection{Bucket Method}
\label{sec:bucket}
Pivot and Partition selects one pivot in each iteration to partition $d$ and applies this recursively in order to create sub-partitions in the manner of a binary search. The Bucket Method, developed by \cite{bucket2020}, can be interpreted as a modification that uses multiple pivots and partitions (buckets) per iteration. 

\begin{algorithm}[!ht]
\SetAlgoLined
\LinesNumbered
\SetKwInput{Input}{Input}
\SetKwInput{Output}{Output}
\Input{vector $d=(d_{1},\cdots,d_{n})$, scaling factor $b$, bucket number $c$, maximum number of iterations $T$.}
\Output{projection $v^{*}$.}
Set $\mathcal{I}=\{1,...,n\}$, $\mathcal{I}_\tau=\emptyset$\;
\For{$t=1:T$}{
    Set $\mathcal{I}_1,...,\mathcal{I}_c:=\emptyset$\;
    \For{$j = 1:c$}{
	    Set $p_j:= (\max_{i\in\mathcal{I}}\{d_i\}-\min_{i\in\mathcal{I}}\{d_i\})\cdot (c-j)/c+\min_{i\in\mathcal{I}}\{d_i\}$\;
	    \For{$i\in\mathcal{I}:d_i\geq p_j$}{
	        Set $\mathcal{I}_{j}:=\mathcal{I}_j\cup\{i\}$, $\mathcal{I}:=\mathcal{I}\backslash\{i\}$
	    }
	}
	\For{$j=1:c$}{
	    Set $p=\frac{\sum_{i \in \mathcal{I}_\tau}d_{i}+\sum_{i \in\mathcal{I}_j}d_{i}-b}{|\mathcal{I}_\tau|+|\mathcal{I}_{j}|}$\;
	    \uIf{$p\geq p_{j}$}{
	        Set $\mathcal{I}:=\mathcal{I}_j$\;
		    Break the inner loop\;
	    }\uElseIf{$j<c$ $\&$ $p> \max_{i\in\mathcal{I}_{j+1}}\{d_i\}$}{
	        Set $\mathcal{I}_\tau:=\mathcal{I}_\tau\cup \mathcal{I}_j$\;
		    Break the outer loop\;
	    }\Else{
	        $\mathcal{I}_\tau:=\mathcal{I}_\tau\cup \mathcal{I}_j$\;
	    }
    }
}
Set $\tau:=\frac{\sum_{i \in \mathcal{I}_\tau}d_i-b}{|\mathcal{I}_\tau|}$\;
Set $v_{i}^{*}:=\mbox{max}\{d_{i}-\tau,0\}$ for all $1\leq i \leq n$\;
\textbf{return} $v^{*}:=(v_{1}^{*},\cdots,v_{n}^{*})$.
\caption{Bucket method}
\label{Alg:bucket}
\end{algorithm}

The algorithm, presented as Algorithm~\ref{Alg:bucket} is initialized with tuning parameters $T$, the maximum number of iterations, and $c$, the number of buckets with which to subdivide the data. In each iteration the algorithm partitions the problem into the buckets $\mathcal{I}_j$ with the inner for loop of line 4, and then calculates corresponding pivot values in the inner for loop of line 10.

The tuning parameters can be determined as follows. Suppose we want the algorithm to find a (final) pivot $\bar \tau$ within some absolute numerical tolerance $D$ of the true pivot $\tau$, i.e.  such that $|\bar \tau - \tau| \leq D$.  This can be ensured (see \cite{bucket2020}) by setting
\[T=\log_c\frac{R}{D},\]
where $R:=\max_{i\in\mathcal{I}}\{d_i\}-\min_{i\in\mathcal{I}}\{d_i\}$ denotes the range of $d$. Perez et al. \cite{bucket2020} prove the worst-case complexity is $O((n+c)\log_c (R/D))$.

We assume uniformly distributed inputs, $d_{1},\dots,d_{n}$ are $\mathrm{i.i.d} \sim U[l,u]$, and we have\\
~\\
\textbf{Proposition 8}
    The Bucket method has an average runtime of $O(cn)$.
\begin{proof}
    Let $\mathcal{I}^{(t)}$ denote the index set $\mathcal{I}$ at the start of iteration $t$ in the outer \textbf{for} loop (line 2), and $\mathcal{I}_j^{(t)}$ denote the index set of the $j$th bucket, $\mathcal{I}_j$, at the end of the first inner \textbf{for} loop (line 4).      
    
    For a given outer \textbf{for} loop iteration $t$ (line 2), the first inner \textbf{for} loop (line 4) uses $O(c|\mathcal{I}^{(t)}|)$ operations. Note that the $\max$ and $\min$ on line 5 can be reused in each iteration, and the nested \textbf{for} loop on line 6 has $|\mathcal{I}^{(t)}|$ iterations.  The second inner \textbf{for} loop (line 10) also uses $O(c|\mathcal{I}^{(t)}|)$ operations.  In line 11, the first sum $\sum_{i \in \mathcal{I}_\tau}d_i$ can be updated dynamically (in the manner of a scan) as a cumulative sum as $\tau$ is updated in line 16 or 19, thus requiring a constant number of operations per iteration $j$. The second sum $\sum_{i \in \mathcal{I}_j}d_i$ is bounded above by $O(|\mathcal{I}^{(t)}|)$ since $\mathcal{I}_j \subseteq \mathcal{I}^{(t)}$. Thus each iteration $j$ of the outer \textbf{for} loop uses $O(c|\mathcal{I}^{(t)}|)$ operations.

    Since $d_1,...,d_n$ are i.i.d $\sim U[l,u]$, then from Proposition 9, the terms from each sub-partition are also i.i.d uniform. So for any $t=1,...,T$ and $j=1,...,c$, $E[|\mathcal{I}_j^{(t)}|]=E[|\mathcal{I}^{(t)}|]/c$. From line 13, $E[|\mathcal{I}^{(t+1)}|]=E[|\mathcal{I}^{(t)}|]/c$. Since $E[|\mathcal{I}^{(1)}|]=n$ then $E[|\mathcal{I}^{(t)}|]=n/c^{t-1}$; thus
    \begin{alignat*}{1}
        E[\sum_{t=1}^{T}|\mathcal{I}^{(t)}|]
        =\sum_{t=1}^{\log_c(R/D)}\frac{n}{c^{t-1}}
        =\frac{c}{c-1}n(1-\frac{D}{R}).
    \end{alignat*}
    Therefore, $E[\sum_{t=1}^{T}c\cdot|\mathcal{I}^{(t)}|] \in O(cn)$.
\end{proof}

\subsection{Projection onto a Weighted Simplex and a Weighted $\ell_1$ Ball}
\label{sec:weightedprob}
The weighted  simplex and the weighted $\ell_1$ ball are
\[\Delta_{w,b}:=\{v \in \mathbb{R}^{n}\ |\ \sum_{i=1}^{n}w_{i}v_{i}=b, v\geq 0\},\]
\[\mathcal{B}_{w,b}:=\{v \in \mathbb{R}^{n}\ |\ \sum_{i=1}^{n}w_{i}|v_{i}|\leq b\},\]
where $w>0$ is a weight vector, and $b>0$ is a scaling factor. \cite{wsimplex2020} show there is a unique $\tau \in \mathbb{R}$ such that
$v^{*}=\max\{d_{i}-w_{i}\tau,0\},\ \forall i=1,\cdots,n$,
where $v^{*}=\mbox{proj}_{\Delta_{w,b}}(d) \in \mathbb{R}^{n}$. Thus pivot-based methods for the unweighted simplex extend to the weighted simplex in a straightforward manner. We present weighted Michelot's method as Algorithm 3 (Appendix D), and weighted Filter as Algorithm 4 (Appendix D).

Our parallelization depends on the choice of serial method for projection onto simplex, since projection onto the weighted simplex requires direct modification rather than oracle calls to methods for the unweighted case. Sort and Scan for weighted simplex projection can be implemented with a parallel merge sort algorithm in Algorithm 5 (Appendix D). Our distributed structure can be applied to Michelot and Condat methods in a similar manner as with the unweighted case; these are presented respectively as Algorithms 6 and 7 (Appendix D). We note that projection onto $\mathcal{B}_{w,b}$ is linear-time reducible to projection onto $\Delta_{w,b}$ \cite[Equation (4)]{wsimplex2020}.
\newpage

\section{Distribution Examples}
\label{sec:examples}
Here we apply Theorem 2 to three examples:\\
	(A) Let $d_{i}\sim U[0,1]$, $n_{1}=10^{5}$, $n_{2}=10^{6}$, $t=0.95$.\\
	(B) Let $d_{i}\sim N[0,1]$, $n_{1}=10^{5}$, $n_{2}=10^{6}$, $t=1.65$.\\
	(C) Let $d_{i}\sim N[0,10^{-3}]$, $n_{1}=10^{5}$, $n_{2}=10^{6}$, $t=1.65\times\sqrt{10^{-3}}=0.05218$.\\
Observe that
\begin{subequations}
\begin{alignat}{3}
    P(\tau >t)\geq &P(\tau>t, |\mathcal{I}_p|\geq \lfloor-\sqrt{2(n\log n)p(1-p)}+np\rfloor)\\
    \label{eq: first term}
    =& P(\tau>t \ |\ |\mathcal{I}_p|\geq \lfloor-\sqrt{2(n\log n)p(1-p)}+np\rfloor)\\
    \label{eq: second term}
    \times &P(|\mathcal{I}_p|\geq \lfloor-\sqrt{2(n\log n)p(1-p)}+np\rfloor))
\end{alignat}
\end{subequations}
(\ref{eq: first term}) can be calculated by (\ref{eq:tau>t}), and (\ref{eq: second term}) can be calculated by (\ref{eq: range of It}).

For	(A), $p=1-F_{d}(t)=0.05$. From $n_{1}p=5000$ and $n_{2}p=50000$, we have $\sqrt{n_{1}p(1-p)}=217.9$ and $\sqrt{n_{2}p(1-p)}=68.9$. So, 
\[P(|\mathcal{I}_t^1|\in [4793,5207])=P(|\mathcal{I}_t^2|\in [49347,50654])=0.9973.\] 
Now, since $f_{d}^{*}(x)=\frac{1}{0.05}=20,\ \mbox{for}\ x\in [0.95,1]$, we have $E=0.975$ and $V=\frac{1}{4800}$. Applying Theorem 2 yields:
\[P(\tau_1>t) \geq 0.99723,\ \forall\ |\mathcal{I}_t^1|\in [4793,5207],\]
\[P(\tau_2>t)\geq 0.99729,\ \forall\ |\mathcal{I}_t^2|\in [49347,50654],\]
which implies the number of active elements in the projection should be less than $5\%$ of $n_1$ or $n_2$ with high probability.

For (B), $p=1-F_{d}(t)=0.05$; similar to the first example,
\[P(|\mathcal{I}_t^1|\in [4793,5207])=P(|\mathcal{I}_t^2|\in [49347,50654])=0.9973.\] 
Together with
\[f_{d}^{*}=\frac{1}{0.05}\frac{1}{\sqrt{2\pi}}e^{-\frac{x^{2}}{2}}=\frac{20}{\sqrt{2\pi}}e^{-\frac{x^{2}}{2}},\ \mbox{for} \ x\in[1.65,+\infty),\]
we can calculate $E$ and $V$:
\[E=\frac{20}{\sqrt{2\pi}}\int_{1.65}^{\infty}xe^{-\frac{x^{2}}{2}}\,dx =\frac{20}{\sqrt{2\pi}}e^{-\frac{1.65^{2}}{2}}=2.045,\]
\[E(x^{2})=\frac{20}{\sqrt{2\pi}}\int_{1.65}^{\infty}x^{2}e^{-\frac{x^{2}}{2}}\,dx=4.375,\]
\[\implies V=E(x^{2})-E^{2}=0.192.\]
Applying Theorem 2,
\[P(\tau_1>t)\geq 0.99704,\ \forall\ |\mathcal{I}_t^1|\in [4793,5207],\]
\[P_{2}(\tau>t)\geq 0.99728,\ \forall\ |\mathcal{I}_t^2|\in [49347,50654],\]
which imply $5\%$ of terms are active after projection with probability $>99\%$.

For (C), $p=1-F_{d}(t)=0.05$. Similar to the previous examples, \[P(|\mathcal{I}_t^1|\in [4793,5207])=P(|\mathcal{I}_t^2|\in [49347,50654])=0.9973.\]
Together with
\[f_{d}^{*}=\frac{20}{\sqrt{2\pi 10^{-3}}}e^{-\frac{x^{2}}{2\times 10^{-3}}},\ \forall \ x\in[0.05218,+\infty),\]
we can calculate $E$ and $V$ as follows,
\[E=\frac{20}{\sqrt{2\pi10^{-3}}}\int_{1.65\times\sqrt{10^{-3}}}^{\infty}xe^{-\frac{x^{2}}{2\times10^{-3}}}\,dx=0.03234,\]
\[E(x^{2})=\frac{20}{\sqrt{2\pi10^{-3}}}\int_{1.65\times\sqrt{10^{-3}}}^{\infty}x^{2}e^{-\frac{x^{2}}{2\times10^{-3}}}\,dx=0.002187,\]
\[\implies V=E(x^{2})-E^{2}=0.001142.\]
Applying Theorem 2, 
\[P(\tau_1>t)\geq 0.99671,\ \forall\ |\mathcal{I}_t^1|\in [4793,5207],\]
\[P(\tau_2>t)\geq 0.99724,\ \forall\ |\mathcal{I}_t^2|\in [49347,50654],\]
and so $5\%$ of terms are active after projection with probability $>99\%$.
\newpage

\section{Algorithm Pseudocode}
\begin{algorithm}[!htbp]
\SetAlgoLined
\LinesNumbered
\SetKwInput{Input}{Input}
\SetKwInput{Output}{Output}
\Input{vector $d=(d_{1},\cdots,d_{n})$}
\Output{projection $v^{*}$.}
 \For{$i=1:n$}{
    $f_{i}=\begin{cases}1, &\mbox{if}\ d_{i} \geq 0\\ 0, &\mbox{otherwise} \end{cases}$\;
 }
 \If{$1^{T}f$ is even}{
    Set $i^{*}=\mbox{argmin}_{i\in 1:n} |d_{i}|$\;
    Update $f_{i^{*}}=1-f_{i^{*}}$\;
 }
 \For{$i=1:n$}{
    Set $v_{i}=d_{i}(-1)^{f_{i}}$\;
 }
 \eIf{$1^{T}\mathrm{proj}_{[-\frac{1}{2},\frac{1}{2}]^{n}}(v)\geq 1-\frac{n}{2}$}{
    \textbf{return} $v^{*}=\mbox{proj}_{[-\frac{1}{2},\frac{1}{2}]^{n}}(d)$
 }{
    Set $v^{*}=\mbox{proj}_{\Delta-\frac{1}{2}}(d)$\;
    \For{$i=1:n$}{
        Update $v^{*}_{i}=v_{i}^{*}(-1)^{f_{i}}$\;
    }
    \textbf{return} $v^{*}$
 }
 \caption{Centered Parity Polytope Projection}
 \label{Alg:PP Project}
\end{algorithm}

\begin{algorithm}[!htbp]
\SetAlgoLined
\LinesNumbered
\SetKwRepeat{Do}{do}{while}
\SetKwInput{Input}{Input}
\SetKwInput{Output}{Output}
\Input{vector $d=(d_{1},\cdots,d_{n})$, scaling factor $b$, weight $w=\{w_1,...,w_n\}$.}
\Output{projection $v^{*}$.}
 Set $\mathcal{I}_p:=\{1,...,n\}$, $\mathcal{I}:=\emptyset$\;
 \Do{$|\mathcal{I}|> |\mathcal{I}_p|$}{
 	Set $\mathcal{I}:=|\mathcal{I}_p|$\;
 	Set $p:=\frac{\sum_{i\in \mathcal{I}_p}w_id_i-b}{\sum_{i\in \mathcal{I}_p}w_{i}^{2}}$\;
 	Set $\mathcal{I}_t=\{i \in \mathcal{I}_p\ | \ \frac{d_i}{w_i}>p\}$\;
 }
 Set $v_{i}^{*}:=\mbox{max}\{d_{i}-w_ip,0\}$ for all $1\leq i \leq n$\;
 \textbf{return} $v^{*}:=(v_{1}^{*},\cdots,v_{n}^{*})$. 
 \caption{Weighted Michelot's method}
 \label{Alg: wmichelot}
\end{algorithm}

\begin{algorithm}[!htbp]
\SetAlgoLined
\LinesNumbered
\SetKwInput{Input}{Input}
\SetKwInput{Output}{Output}
\Input{vector $d=(d_{1},\cdots,d_{n})$, scaling factor $b$, weight $w$.}
\Output{Index set $\mathcal{I}_p$.}
 Set $\mathcal{I}_p:=\{1\}$, $\mathcal{I}_w:=\emptyset$, $p=:\frac{w_{1}d_{1}-b}{w_{1}^{2}}$\;
 \For{$i=2:n$}{
 	\If{$\frac{d_{i}}{w_{i}}>p$}{
 		Set $p:=\frac{w_{i}d_{i}+\sum_{j\in \mathcal{I}_pw_{j}d_{j}}-b}{w_{i}^{2}+\sum_{j\in \mathcal{I}_pw_{j}^{2}}}$\;
 		\eIf{$p>\frac{w_{i}d_{i}-b}{w_{i}^{2}}$}{
 			Set $\mathcal{I}_p:=\mathcal{I}\cup \{i\}$\;
 		}{
 			Set $\mathcal{I}_w:=\mathcal{I}_w\cup\mathcal{I}_p$\;
 			Set $\mathcal{I}_p=\{i\}$, $p:=\frac{w_{i}d_{i}-b}{w_{i}^{2}}$\;
 		}
 	}
 }
 \If{$|\mathcal{I}_w| \neq 0$}{
 	\For{$i\in \mathcal{I}_w:\frac{d_i}{w_i} > p$}{
 		Set $\mathcal{I}_p:=\mathcal{I}_p\cup\{i\}$\;
 		Set $p:=\frac{w_{i}d_{i}+\sum_{j\in \mathcal{I}_pw_{j}d_{j}}-b}{w_{i}^{2}+\sum_{j\in \mathcal{I}_pw_{j}^{2}}}$\;
 	}
 }
 \textbf{return} $\mathcal{I}_p$. 
 \caption{Weighted Filter \cite{wsimplex2020}}
 \label{Alg: wfilter}
\end{algorithm}

\begin{algorithm}[!ht]
\SetAlgoLined
\LinesNumbered
\SetKwInput{Input}{Input}
\SetKwInput{Output}{Output}
\Input{vector $d=(d_{1},\cdots,d_{n})$, scaling factor $b$, weight $w$.}
\Output{projection $v^{*}$.}
 Set $z:=\{\frac{d_{i}}{w_{i}}\}$\;
 Parallel sort $z$ as $z_{(1)}\geq \cdots \geq z_{(n)}$, and apply this order to $d$ and $w$ \;
 Find $\kappa:=\mathrm{max}_{k=1,\cdots,n}\{\frac{\sum_{i=1}^{k}w_{i}d_{i}-b}{\sum_{i=1}^{k}w_{i}^{2}}\leq z_{k}\}$\;
 Set $\tau=\frac{\sum_{i=1}^{\kappa}w_{i}d_{i}-b}{\sum_{i=1}^{\kappa}w_{i}^{2}}$\;
 Parallel set $v_{i}^{*}:=\mbox{max}\{d_{i}-w_{i}\tau, 0\}$ for all $1\leq i\leq n$\;
 \textbf{return} $v^{*}=(v_{1}^{*},\cdots,v_{n}^{*})$. 
 \caption{Weighted Parallel Sort and Parallel Scan}
 \label{Alg: PWS & PWS}
\end{algorithm}

\begin{algorithm}[!htbp]
\SetAlgoLined
\LinesNumbered
\SetKwInput{Input}{Input}
\SetKwInput{Output}{Output}
\Input{vector $d=(d_{1},\cdots,d_{n})$, scaling factor $b$, weight $w$,  $k$ cores.}
\Output{projection $v^{*}$.}
 Partition $d$ into subvectors $d^1,...,d^k$ of dimension $\leq \frac{n}{k}$\;
 Set $v^i:=\texttt{Weighted\_Pivot\_Project}(v^i,w,b)$\quad \mbox{(distributed across cores $i=1,...,k$)}\;
 Set $\hat{v}:=\cup_{i=1}^{k}\{v_j^i\ |\ v_j^i>0\}$\;
 Set $v^*:=\texttt{Weighted\_Pivot\_Project}(\hat{v},w,b)$\;
 \textbf{return} $v^{*}$.
 \caption{Distributed Weighted Pivot and Project}
 \label{Alg: Dwpivot}
\end{algorithm}

\begin{algorithm}[!htbp]
\SetAlgoLined
\LinesNumbered
\SetKwRepeat{Do}{do}{while}
\SetKwInput{Input}{Input}
\SetKwInput{Output}{Output}
\Input{vector $d=(d_{1},\cdots,d_{n})$, scaling factor $b$, weight $w$, $k$ cores.}
\Output{projection $v^{*}$.}
 Partition $d$ into subvectors $d^1,...,d^k$ of dimension $\leq \frac{n}{k}$\;
 Set $v^i:=\texttt{Weighted\_Condat\_Project}(v^i,w,b)$\quad \mbox{(distributed across cores i=1,...,k)}\;
 Set $\mathcal{I}_p:=\cup_{i=1}^{k}\{j\ |\ v^i_j>0\}$\;
 Set $p:=\frac{\sum_{i\in \mathcal{I}_p}w_id_i-b}{\sum_{i\in\mathcal{I}_p}w_i^2}$, $\mathcal{I}:=\emptyset$\;
 \Do{$|\mathcal{I}|>|\mathcal{I}_p|$}{
 	Set $\mathcal{I}=\mathcal{I}_p$\;
 	\For{$i\in\mathcal{I}$}{
 	    \If{$\frac{d_i}{w_i}\leq p$}{
 	       Set $\mathcal{I}:=\mathcal{I}_p\backslash \{i\}$, $p:=\frac{\sum_{i\in \mathcal{I}}w_id_i-b}{\sum_{i\in\mathcal{I}}w_i^2}$\;
 	    }
 	}    
 }
 Set $v_{i}^{*}=\mbox{max}\{d_{i}-w_{i}p,0\}$ for all $1\leq i \leq n$\;
 \textbf{return} $v^{*}=(v_{1}^{*},\cdots,v_{n}^{*})$. 
 \caption{Distributed Weighted Condat}
 \label{Alg: Dwcondat}
\end{algorithm}
\phantom{a}
\newpage

\section{Additional Experiments}
All code and data can be found at: \href{https://github.com/foreverdyz/Parallel_Projection}{Github}\footnote{https://github.com/foreverdyz/Parallel\_Projection} or the IJOC repository \cite{code}

\subsection{Testing Theoretical Bounds}\label{sec:bounds}
\begin{figure}[!htbp]
\centering
\includegraphics[width=0.7\linewidth]{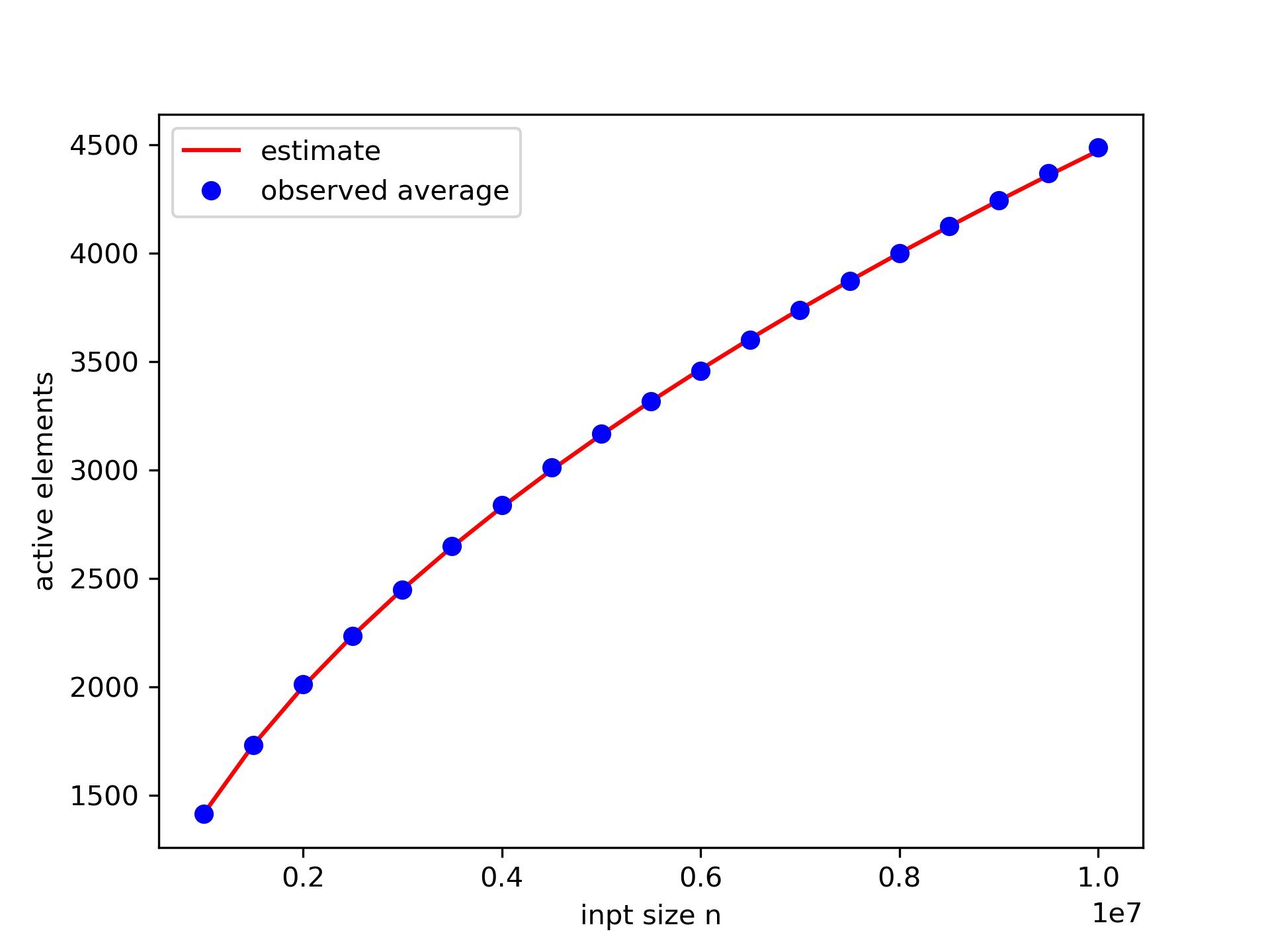}
\caption{Observed number of active terms vs bound value of $\sqrt{2n}$}
\label{fig: average_active}
\end{figure}

For Theorem 1, we calculate the average number of active elements in projecting a vector $d$, drawn i.i.d. from $U[0,1]$, onto the simplex $\Delta_{1}$, with $n$ between $10^{6}$ and $10^{7}$ and 10 trials per size. This empirical result is compared against the corresponding asymptotic bound of $\sqrt{2n}$ given by Theorem 1.  Results are shown in Figure \ref{fig: average_active}, and demonstrate that our asymptotic bound is rather accurate for small $n$.

\begin{figure}[!htbp]
\centering
\includegraphics[width=0.7\linewidth]{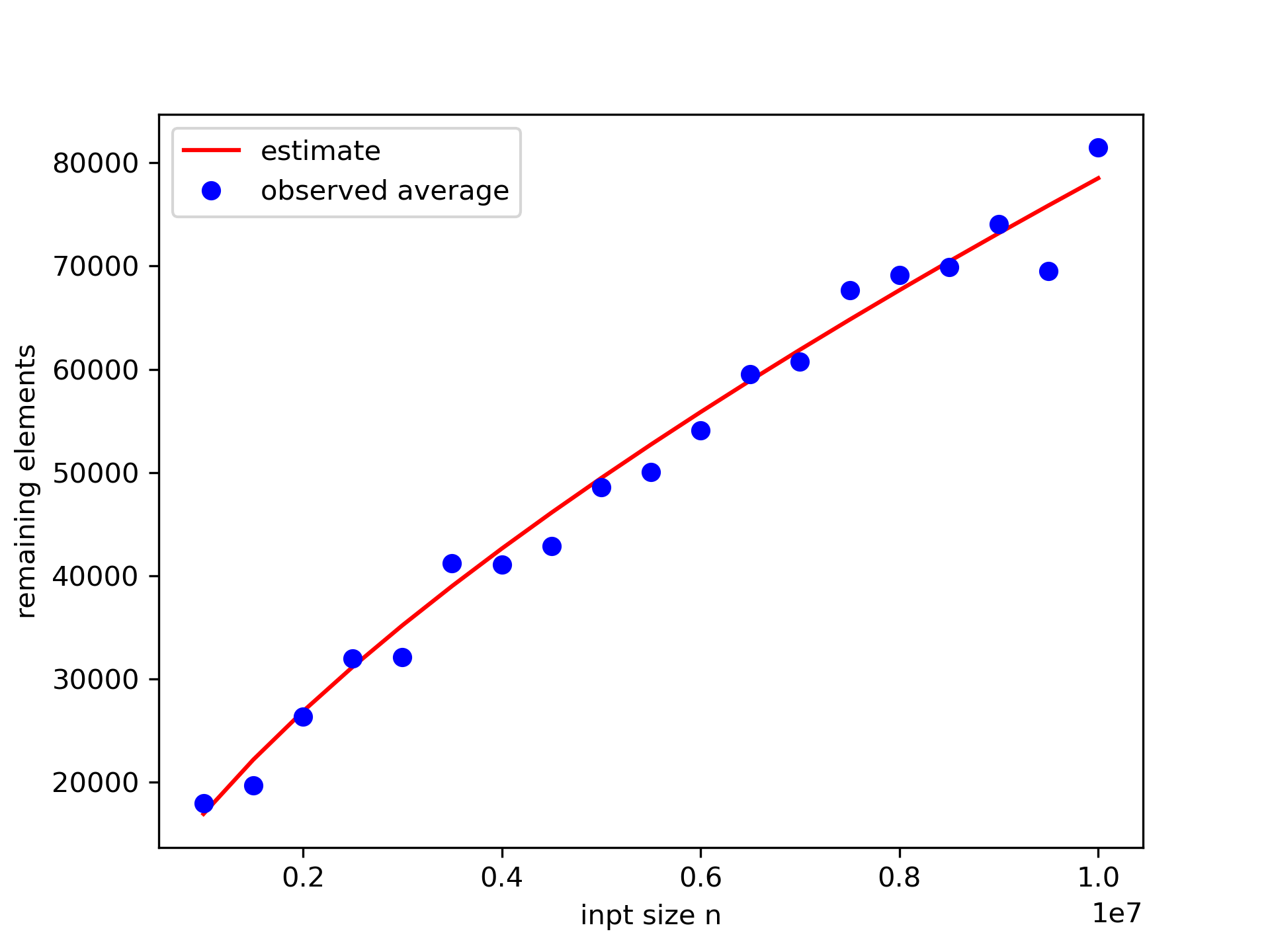}
\caption{Observed number of remaining terms after applying Filter vs bound value of $(2.2n)^{\frac{2}{3}}$}
\label{fig:filter_result}
\end{figure}

Similarly, for Proposition 6, we conduct the same experiments and compare the results against the function $(2.2n)^{\frac{2}{3}}$ in Figure \ref{fig:filter_result}, where the constant was found empirically.

For Lemma 1, we run Algorithm 3 (in the main paper) on $U[0,1]$ i.i.d. distributed inputs $d_i$ with scaling factor $b=1$, size $n=10^{6}$, and 100 trials. We compare the (average) remaining number of elements after each iteration of Michelot's method and against the geometric series with a ratio as $\frac{1}{2}$. We find the average number of remaining terms after each loop of Michelot's method is close to the corresponding value of the geometric series with a ratio as $\frac{1}{2}$ from $10^{6}$. So, the conclusion from Lemma 1, which claims that the Michelot method approximately discards half of the vector in each loop when the size is big, is accurate.

\begin{figure}[!htbp]
\centering
\includegraphics[width=0.7\linewidth]{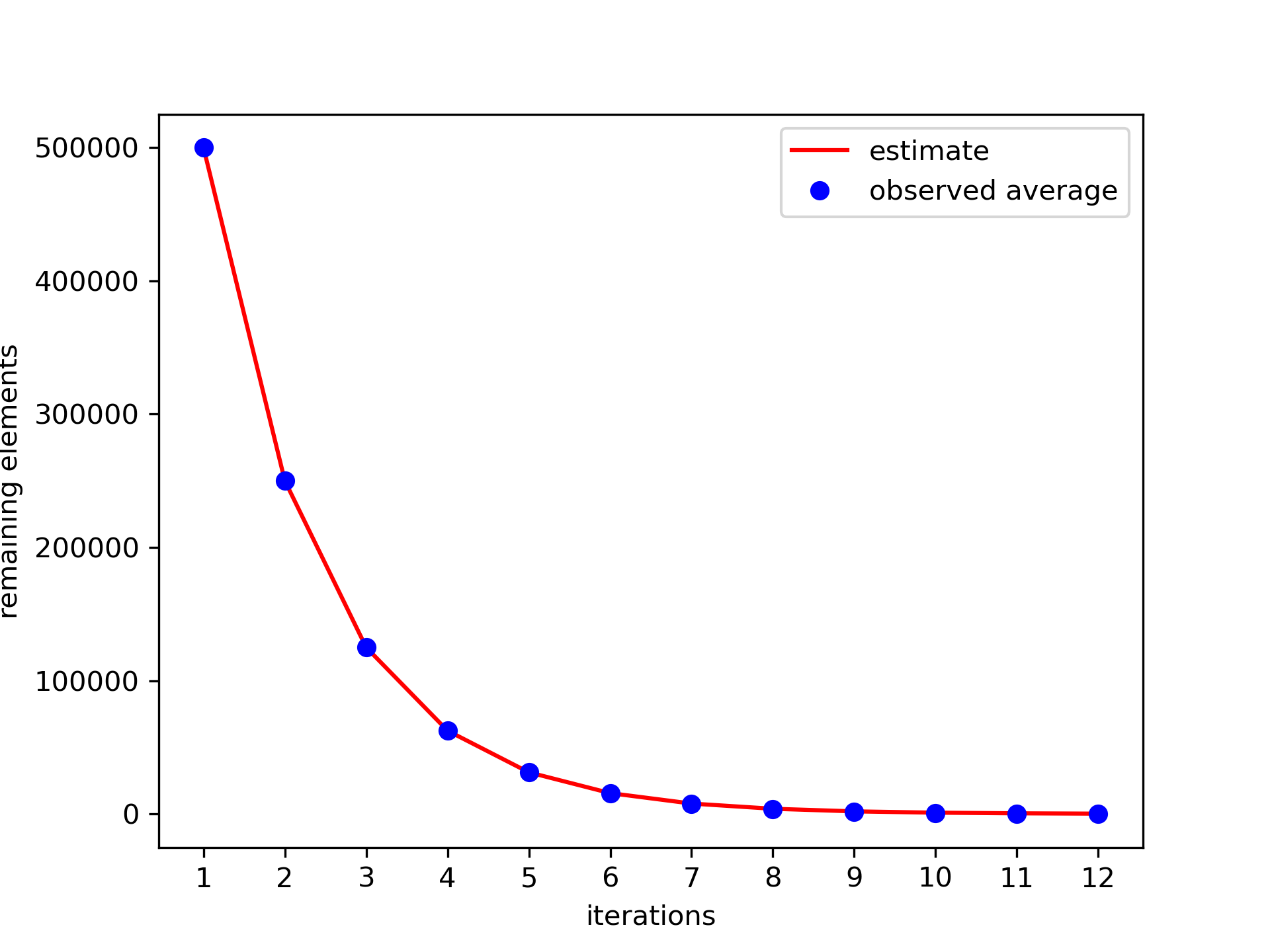}
\caption{Observed number of remaining terms after each Michelot iteration vs bound value of geometric series with a ratio of $\frac{1}{2}$}
\label{fig: michelot_loop}
\end{figure}

\subsection{Robustness Test}
\label{sec:robust}

\begin{figure}[!htbp]
    \centering
    \begin{minipage}[b]{0.45\textwidth}
         \centering
         \includegraphics[width=\textwidth]{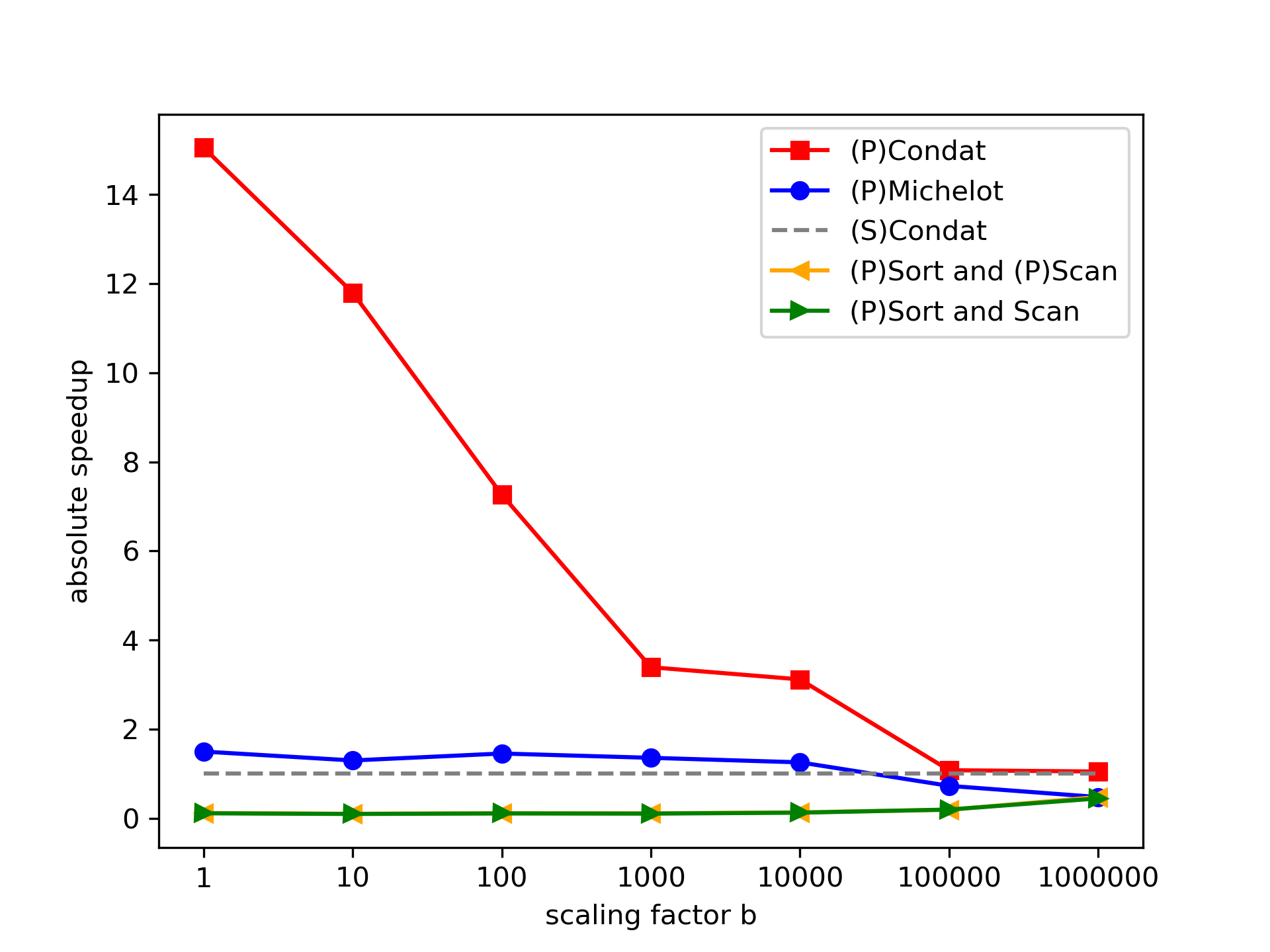}
    \end{minipage}
    \hfill
    \begin{minipage}[b]{0.45\textwidth}
         \centering
         \includegraphics[width=\textwidth]{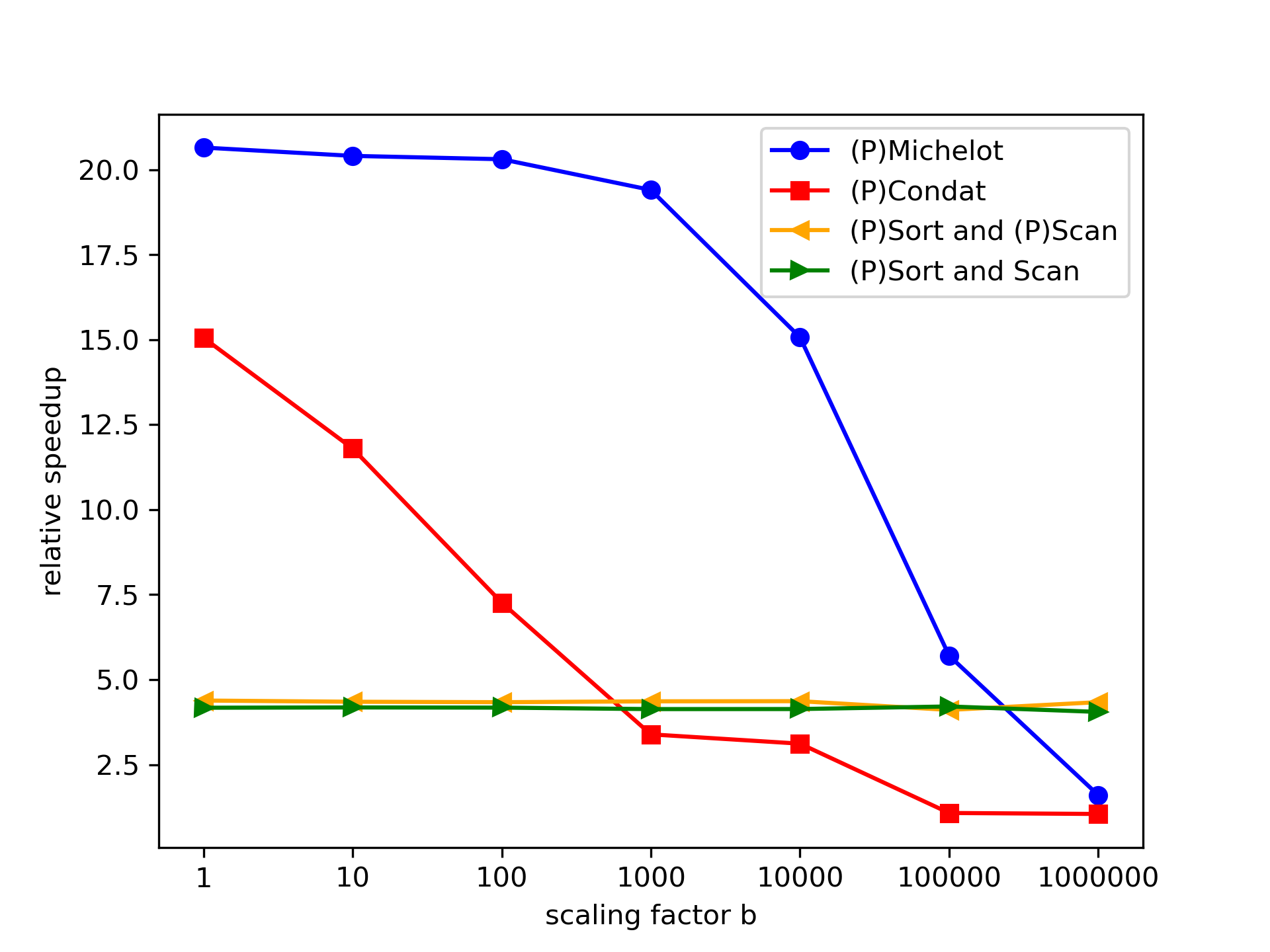}
    \end{minipage}
    \caption{Speedup vs b in Simplex Projection. Each line represents a different projection method, and the input distribution is $N(0,1)$.}
    \label{fig:b1}
\end{figure}
\begin{figure}[!htbp]
    \centering
    \begin{minipage}[b]{0.45\textwidth}
         \centering
         \includegraphics[width=\textwidth]{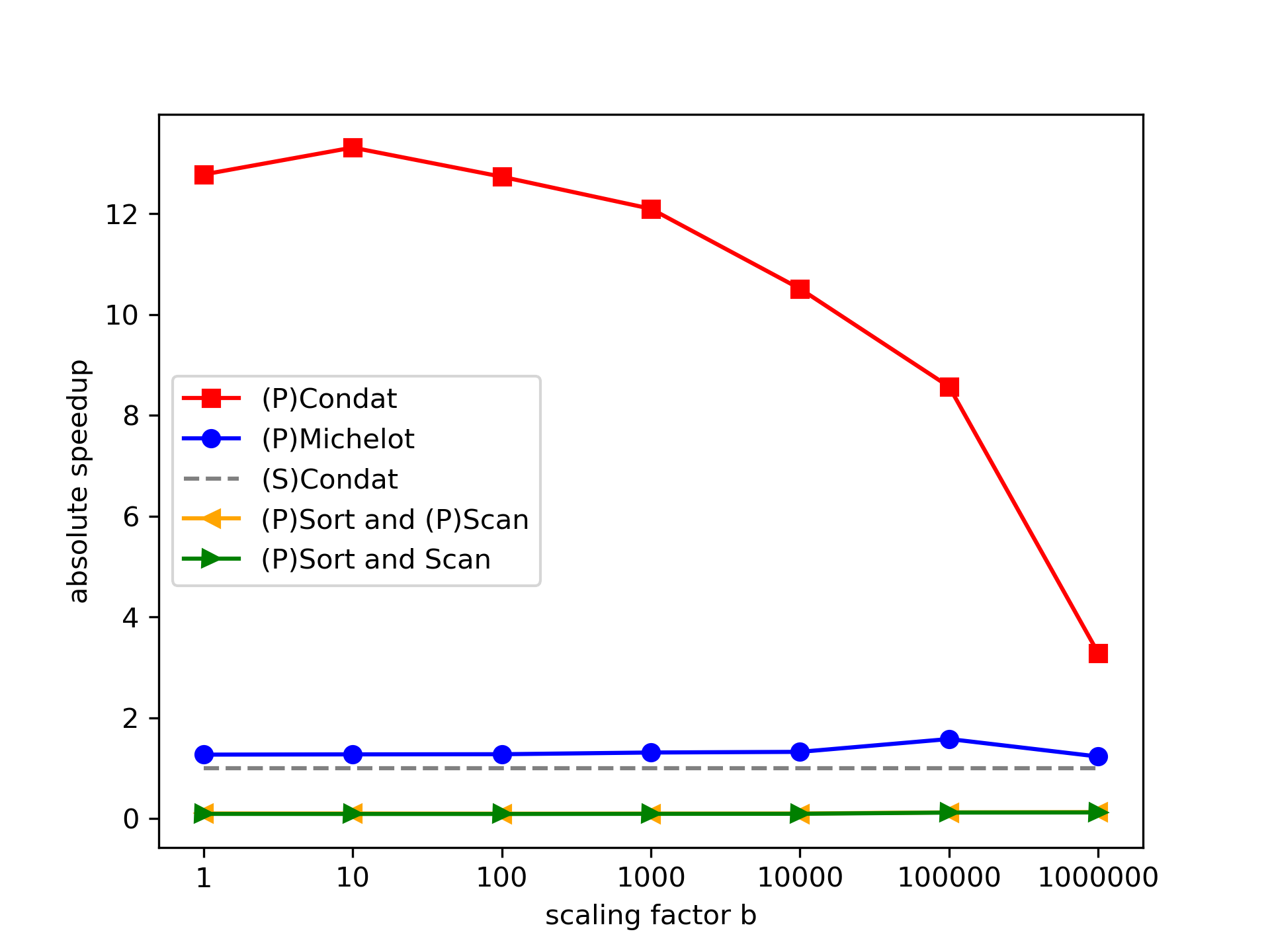}
    \end{minipage}
    \hfill
    \begin{minipage}[b]{0.45\textwidth}
         \centering
         \includegraphics[width=\textwidth]{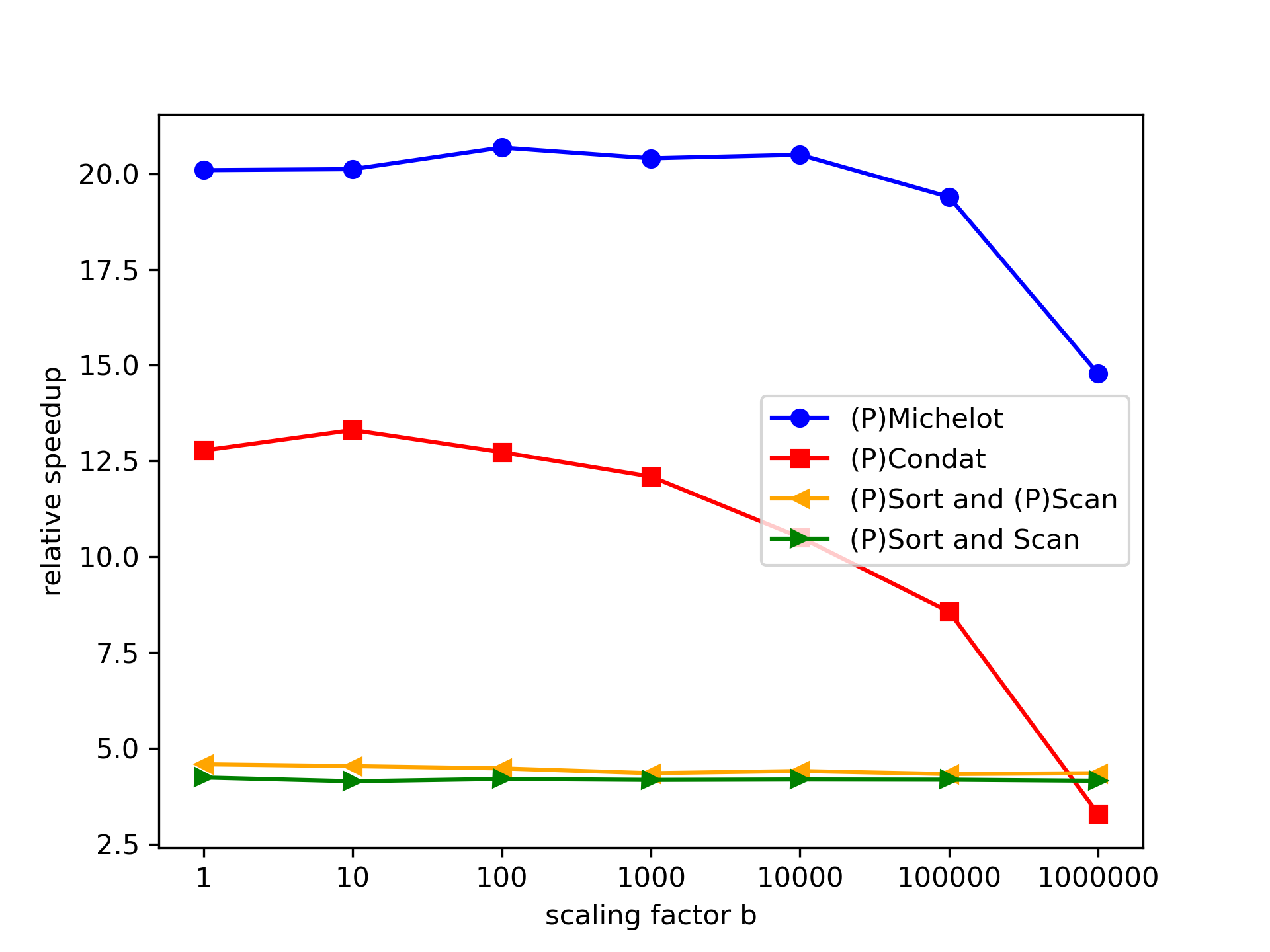}
    \end{minipage}
    \caption{Speedup vs b in Simplex Projection. Each line represents a different projection method, and the input distribution is $N(0,100)$.}
    \label{fig:b100}
\end{figure}

We created examples with a wide range of scaling factors $b = 1, 10, 10^2, 10^3$, $10^4, 10^5, 10^6$ and drawing $d_i\sim N(0,1)$ and $N(0, 100)$ with a size of $n = 10^8$ for serial methods and parallel methods with $40$ cores. Results are given in Figure~\ref{fig:b1} ($d_i\sim N(0,1)$) and Figure~\ref{fig:b100} ($d_i\sim N(0,100)$). 

\subsection{Weighted simplex, weighted $\ell_{1}$ ball}
\begin{figure}[!htbp]
    \centering
    \begin{minipage}[b]{0.45\textwidth}
         \centering
         \includegraphics[width=\textwidth]{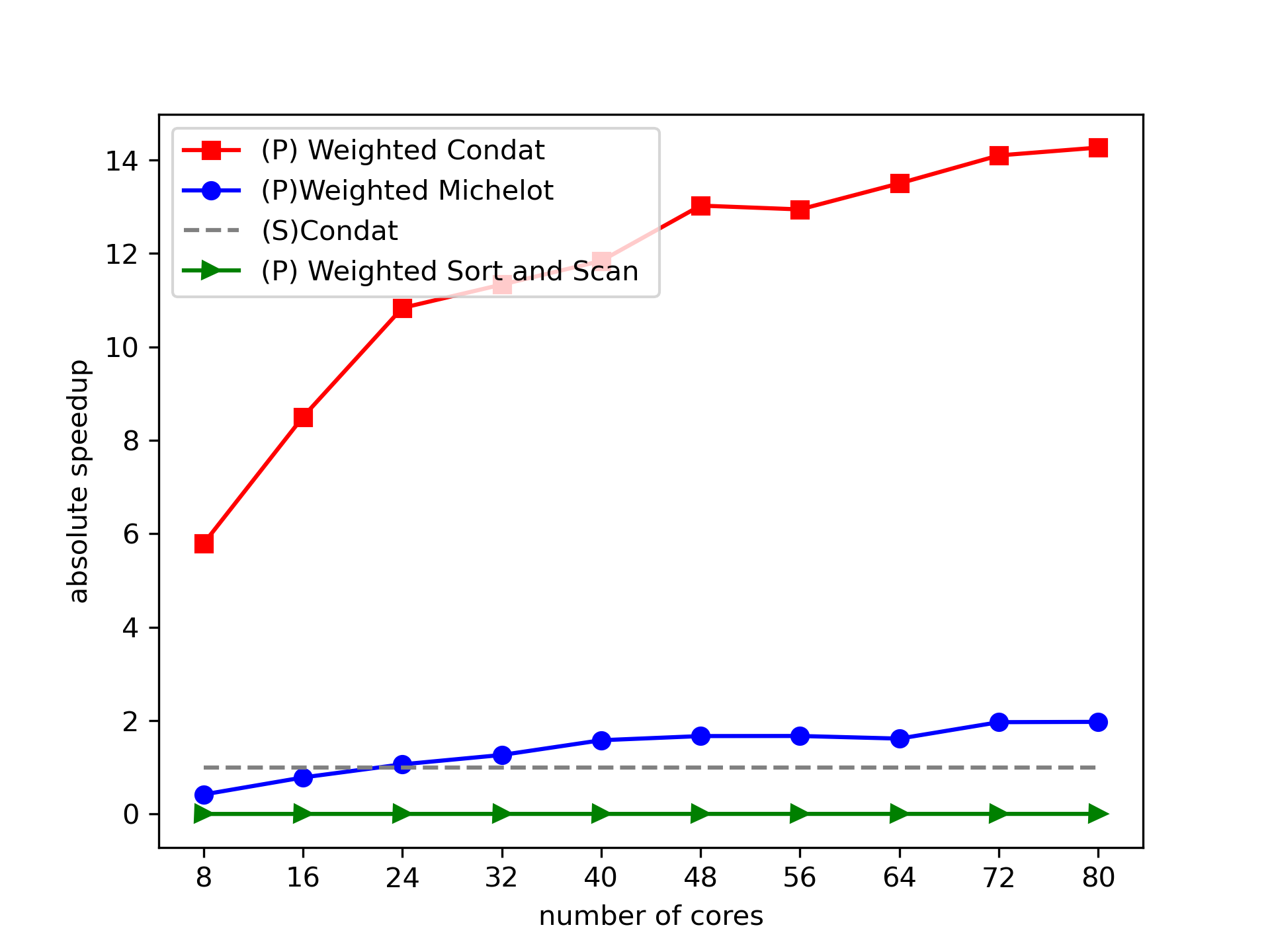}
    \end{minipage}
    \hfill
    \begin{minipage}[b]{0.45\textwidth}
         \centering
         \includegraphics[width=\textwidth]{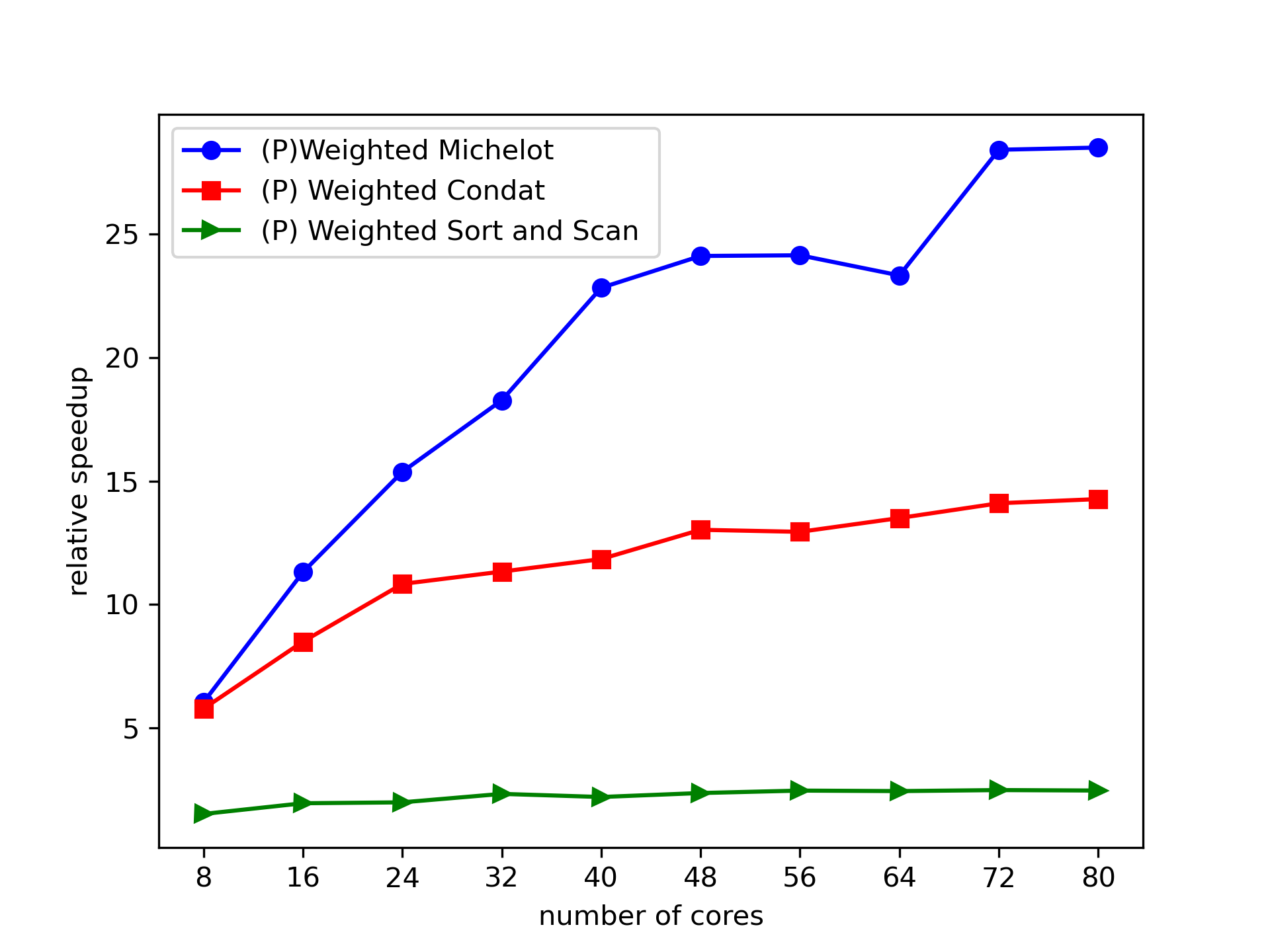}
    \end{minipage}
    \caption{Speedup vs cores in weighted simplex projection. Each line represents a different projection method.}
    \label{fig:wsimplex}
\end{figure}

\begin{figure}[!htbp]
    \centering
    \begin{minipage}[b]{0.45\textwidth}
         \centering
         \includegraphics[width=\textwidth]{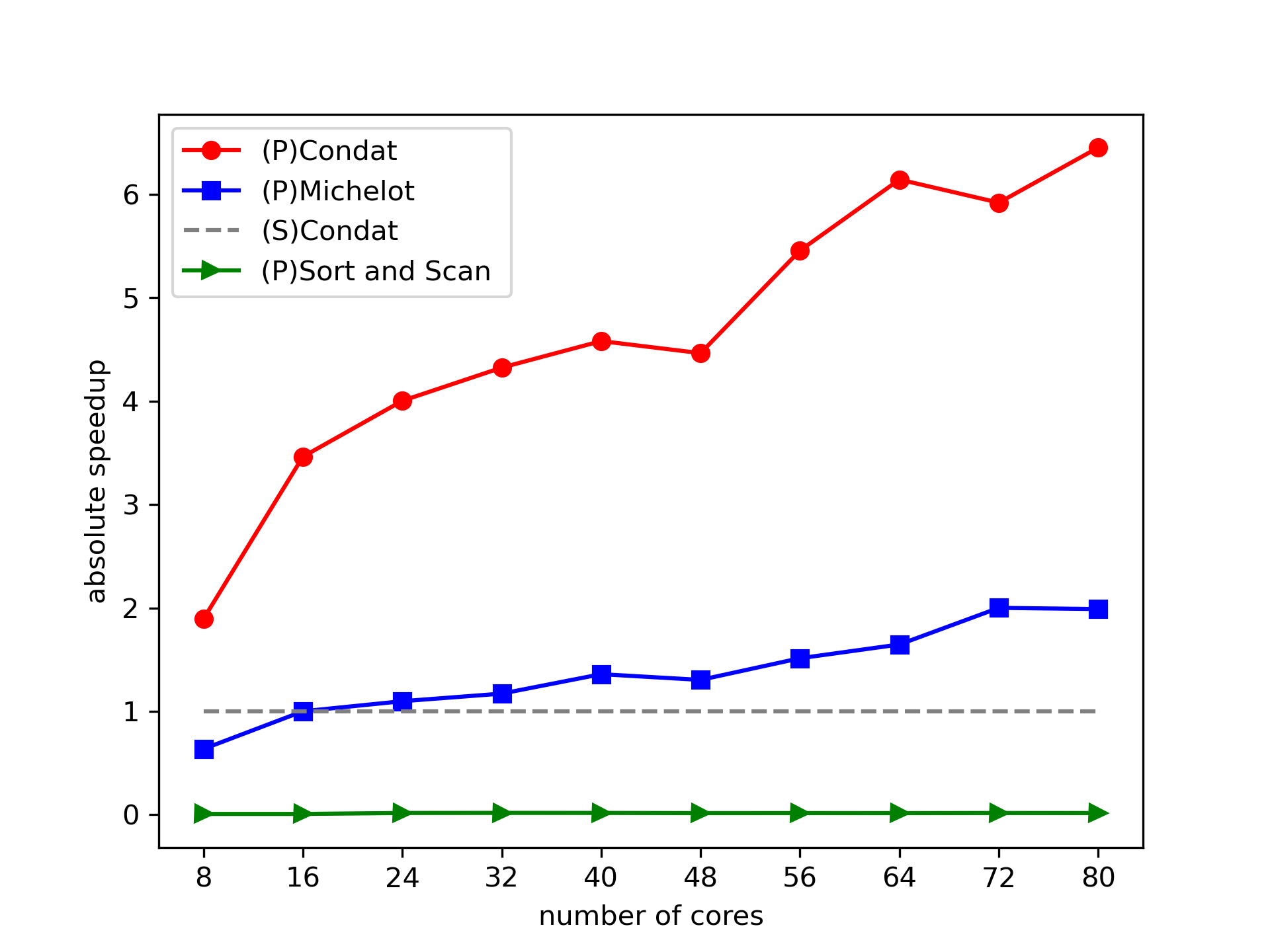}
    \end{minipage}
    \hfill
    \begin{minipage}[b]{0.45\textwidth}
         \centering
         \includegraphics[width=\textwidth]{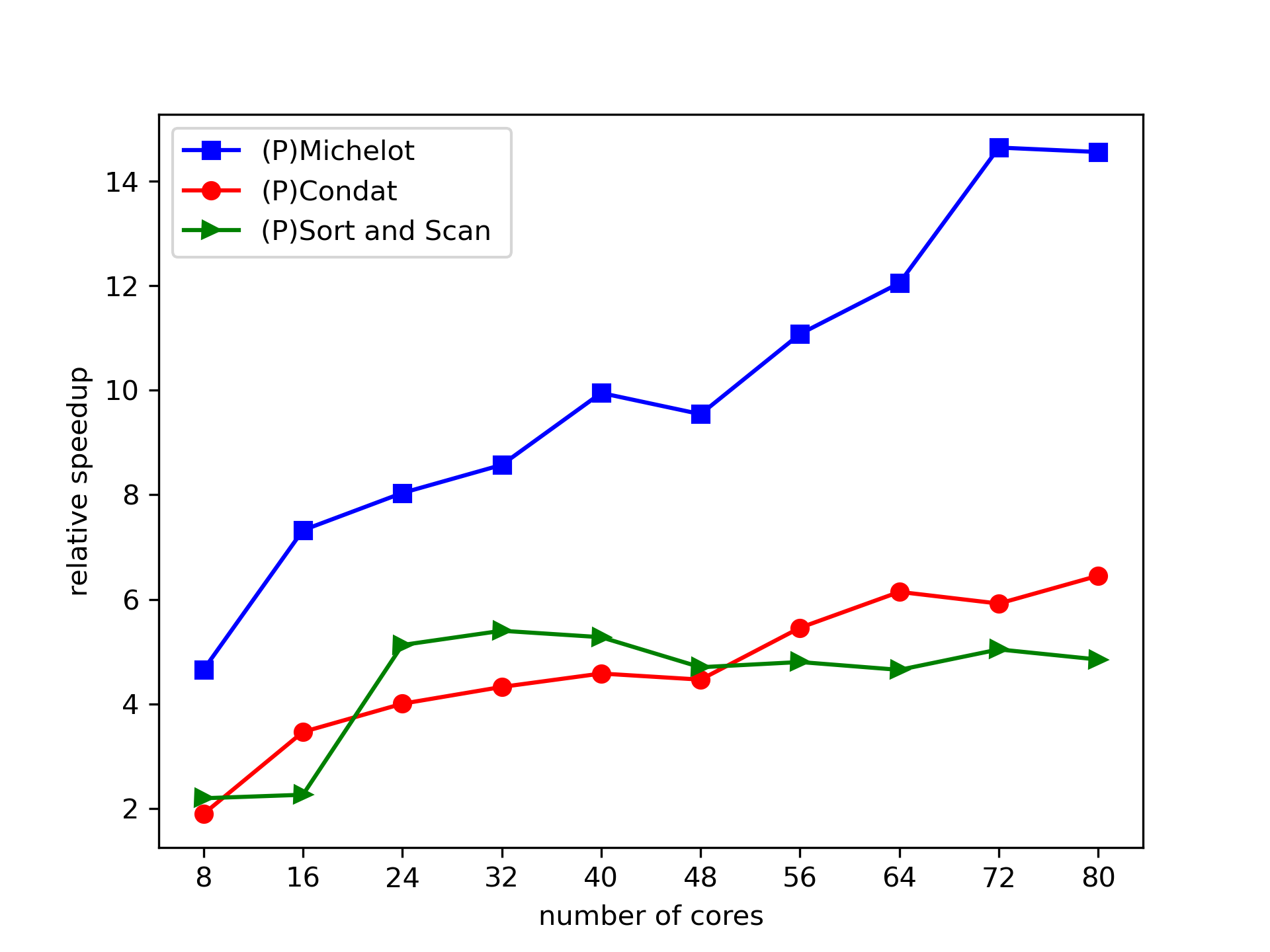}
    \end{minipage}
    \caption{Speedup vs cores in weighted $\ell_1$ ball projection. Each line represents a different projection method.}
    \label{fig:wl1ball}
\end{figure}

We conduct weighted simplex and weighted $\ell_1$ ball projection experiments with various methods to solve a problem with a problem size of $n=10^{8}-1$. Inputs $d_i$ are drawn i.i.d. from $N(0, 1)$ and weights $w_i$ are drawn i.i.d. from $U[0, 1]$. Algorithms were implemented as described in Sec~\ref{sec:weightedprob}. Results for weighted simplex are shown in Figure~\ref{fig:wsimplex} and results for weighted $\ell_1$ ball projection are shown in Fig.~\ref{fig:wl1ball}.  Slightly more modest speedups across all algorithms are observed in the weighted simplex compared to the unweighted experiments; nonetheless, the general pattern still holds, with parallel Condat attaining superior performance with up to 14x absolute speedup.  In the weighted $\ell_1$ ball projection, we observe that sort and scan has similar relative speedups to our parallel Condat's method.  However, the underlying serial Condat's method is considerably faster.

\subsection{Runtime Fairness Test}
We provide a runtime fairness test for serial methods (e.g. Sort and Scan method, Michelot's method, Condat's method) and their respective parallel implements. We restrict both serail and parallel methods to use only one core to solve a problem with a size of $n=10^8$ and a scaling factor $b = 1$. Inputs $d_i$ drawn i.i.d. from $u[0,1]$, and $b = 1$.  Results are provided in Table~\ref{tab:fair}.

\begin{table}[!htbp]
    \centering
    \begin{tabular}{lccc}
        \toprule[1pt]
        Method &Runtime\\
        \hline
        Sort + Scan &$\mathrm{1.037e+01}$\\
        (P)Sort + Scan &$\mathrm{1.571e+01}$\\
        (P)Sort + Partial Scan &$\mathrm{1.505e+01}$\\
        \hline
        Michelot &$\mathrm{4.030}$\\
        (P) Michelot &$\mathrm{4.229}$\\
        \hline
        Condat &$\mathrm{2.429e-01}$\\
        (P) Condat &$\mathrm{2.497e-01}$\\
        \bottomrule[1pt]
    \end{tabular}
    \caption{Runtime (s) for projection onto a simplex in fairness test}
    \label{tab:fair}
\end{table}

\subsection{Discussion on dense projections} \label{sec:explain}
To project a vector $d\in \mathbb{R}^n$ onto the $\ell_b$ ball, we first check if $\sum_{i=1}^n |d_i| \leq b$. If this condition holds, then $d$ is already within the $\ell_b$ ball. However, if $\sum_{i=1}^n |d_i| > b$, we project $|d| :=(|d_1|,...,|d_n|)$ onto the simplex with a scaling factor of $b$. As noted earlier, we have $\tau \geq \frac{\sum_{i=1}^n |d_i| - b}{n} >0$, which means that all zero terms in $d$ are inactive in the projection of $|d|$ onto the simplex. Therefore, we only need to project a subvector $(|d_i|)_{i:d_i\neq 0}$ onto the simplex. This is why when $d$ is sparse, it is probably better to use serial projection methods instead of their parallel counterparts.

\end{document}